\documentclass{gtart}

\def\ifplaintex{\expandafter\ifx\csname documentclass\endcsname\relax}

\ifplaintex 
\else
\fi

\expandafter\ifx\csname beginpicture\endcsname\relax
\expandafter\ifx\csname documentclass\endcsname\relax
\input pictex \else
\input prepictex \input pictex \input postpictex \fi\fi

\def\gt{{\mathsurround=0pt\it $\cal G\mskip-2mu$eometry \&\ 
$\cal T\!\!$opology}}        

\def\gtp{{\mathsurround=0pt\it $\cal G\mskip-2mu$eometry \&\ 
$\cal T\!\!$opology $\cal P\!$ublications}}  

\def\lognumber#1{\def\thelognumber{#1}}
\def\volumenumber#1{\def\thevolumenumber{#1}}
\def\papernumber#1{\def\thepapernumber{#1}}
\def\volumeyear#1{\def\thevolumeyear{#1}}

\def\pagenumbers#1#2{\def\startpage{#1}\def\finishpage{#2}}
\def\published#1{\def\publishdate{#1}}
\def\proposed#1{\def\theproposer{#1}}
\def\seconded#1{\def\theseconders{#1}}
\def\received#1{\def\receiveddate{#1}}

\def\accepted#1{\def\accepteddate{#1}}

\long\def\asciiabstract#1{\long\def\theasciiabstract{#1}}

\let\\\par\let\thelognumber\relax
\let\thevolumenumber\relax\let\thepapernumber\relax
\let\thevolumeyear\relax\let\thesamplenumber\relax\let\startpage\relax
\let\finishpage\relax\let\publishdate\relax\let\receiveddate\relax
\let\reviseddate\relax\let\accepteddate\relax\let\theasciititle\relax
\let\theasciiauthors\relax
\let\theasciiabstract\relax
\let\theasciiemail\relax\let\theshortauthors\relax\let\theshorttitle\relax

\long\def\maketitlep{   

\count0=\startpage

\gt\hfill      
\beginpicture
\setcoordinatesystem units <0.33truein, 0.33truein> point at 2.2 0.9
\setplotsymbol ({$\cal G$})
\plotsymbolspacing=9truept
\circulararc 315 degrees from 0 1 center at 0 0
\setplotsymbol ({$\cal T$})
\circulararc 315 degrees from 1 -1 center at 1 0
\endpicture
\break
{\small\ifx\thesamplenumber\relax 
Volume \else Sample
\fi\thevolumenumber\ (\thevolumeyear)
\startpage--\finishpage\nl
Published: \publishdate}
\vglue 0.5truein plus 0.4fil minus 0.1truein

{\parskip=0pt\leftskip 0pt plus 1fil\def\\{\par\smallskip}{\ifplaintex\large
\else\Large\fi\bf\thetitle}\par\medskip}   

\vglue 0pt plus 0.1fil 

{\parskip=0pt\leftskip 0pt plus 1fil\def\\{\par}{\sc\theauthors}
\par\medskip}

\vglue 0pt plus 0.1fil 

{\small\parskip=0pt\let\newline\\
{\leftskip 0pt plus 1fil\def\\{\par}{\sl\theaddress}\par}
\expandafter\ifx\theemail\relax    
\relax\else\vglue 5pt plus 0.02fil minus 2pt\def\\{\stdspace{\rm 
and}\stdspace} 
\cl{Email:\stdspace\tt\theemail}\fi
\ifx\theurl\relax                  
\relax\else\vglue 5pt plus 0.02fil minus 2pt\def\\{\stdspace{\rm 
and}\stdspace}
\cl{URL:\stdspace\tt\theurl}\fi\par}

\vglue 7pt plus 0.3fil minus 3pt

{\bf Abstract}
\vglue 5pt plus 0.1fil minus 2pt

\theabstract

\vglue 7pt plus 0.3fil minus 3pt

{\bf AMS Classification numbers}\quad Primary:\quad \theprimaryclass

Secondary:\quad \thesecondaryclass

\vglue 5pt plus 0.3fil minus 2pt

{\bf Keywords}\quad \thekeywords

\vglue 10pt plus 0.5fil minus 5pt

{\small  Proposed: \theproposer\hfill Received: \receiveddate\nl
Seconded: \theseconders\hfill 
\ifx\reviseddate\relax                         
Accepted: \accepteddate                        
\else
Revised: \reviseddate                          
\fi}
\eject
}       

\let\maketitlepage\maketitlep
\let\maketitle\maketitlepage

\font\phead=cmsl9 scaled 950
\font=cmsl9 scaled 1050
\font\pnum=cmbx10 scaled 913
\font\lnum=cmbx10 
\font\pfoot=cmsl9 scaled 950
\font=cmsl9 scaled 1050
\ifplaintex
\headline{\vbox to 0pt{\vskip -4.5mm\line{\small\phead\ifnum
\count0=\startpage ISSN 1364-0380 (on line)
1465-3060 (printed) \hfill {\pnum\folio}\else\ifodd\count0\def\\{ }%
\ifx\theshorttitle\relax\thetitle\else\theshorttitle\fi\hfill{\pnum\folio}
\else\def\\{ and }{\pnum\folio}\hfill\ifx\theshortauthors\relax\theauthors
\else\theshortauthors\fi\fi\fi}\vss}}
\footline{\vbox to 0pt{\vglue 0mm\line{\small\pfoot\ifnum\count0=\startpage
\copyright\ \gtp\hfill\else
\gt, Volume \thevolumenumber\ (\thevolumeyear)\hfill\fi}\vss
}}
\else
\makeatletter
\def\@oddhead{{\small\lhead\ifnum\count0=\startpage ISSN 1364-0380 (on line)
1465-3060 (printed) \hfill {\lnum\number\count0}\else\ifodd\count0
\def\\{ }\ifx\theshorttitle\relax \thetitle \else\theshorttitle\fi\hfill
{\lnum\number\count0}\else\def\\{ and }{\lnum\number\count0}
\hfill\ifx\theshortauthors\relax 
\theauthors\else\theshortauthors\fi\fi\fi}}\def\@evenhead{\@oddhead}
\def\@oddfoot{\small\lfoot\ifnum\count0=\startpage\copyright\ \gtp\hfill\else
\gt, Volume \thevolumenumber\ (\thevolumeyear)\hfill\fi}
\def\@evenfoot{\@oddfoot}
\makeatother
\fi

\newwrite\gtoutfile
\long\gdef\makeheadfile{  
{\def\\{, }\def\s{ }
\immediate\openout\gtoutfile head.xxx
\immediate\write\gtoutfile{To: math@arxiv.org}
\immediate\write\gtoutfile{Subject: put or rep NNNNN:pppp}
\immediate\write\gtoutfile{--text follows this line--}
\immediate\write\gtoutfile{Proxy-for: \ifx\theasciiauthors\relax
\theauthors\else\theasciiauthors\fi\s<\ifx\theasciiemail\relax\theemail\else\theasciiemail\fi>}
\immediate\write\gtoutfile{\noexpand\\}
\immediate\write\gtoutfile{Authors: \ifx\theasciiauthors\relax
\theauthors\else\theasciiauthors\fi}
\immediate\write\gtoutfile{Title: \ifx\theasciititle\relax
\thetitle\else\theasciititle\fi}
\immediate\write\gtoutfile{Subj-class: GT or SG or MG etc}
\immediate\write\gtoutfile{MSC-class: \theprimaryclass\ifx\thesecondaryclass\relax\else, \thesecondaryclass\fi}
\immediate\write\gtoutfile{Journal-ref: Geom. Topol. \thevolumenumber
(\thevolumeyear) \startpage-\finishpage}
\immediate\write\gtoutfile{Comments: Published by Geometry and Topology at}
\immediate\write\gtoutfile{\s\s http://www.maths.warwick.ac.uk/gt/GTVol\thevolumenumber/paper\thepapernumber.abs.html}
\immediate\write\gtoutfile{\noexpand\\}
\immediate\write\gtoutfile{}
\ifx\theasciiabstract\relax
\immediate\write\gtoutfile{\theabstract}\else
\immediate\write\gtoutfile{\theasciiabstract}\fi
\immediate\write\gtoutfile{}
\immediate\write\gtoutfile{\noexpand\\}
\immediate\write\gtoutfile{}
\immediate\closeout\gtoutfile}}  

\def\maketitlepage{\maketitlep\makeheadfile}
\let\maketitle\maketitlepage

\def\ifplaintex{\expandafter\ifx\csname documentclass\endcsname\relax}

\ifplaintex 
\else
\fi

\expandafter\ifx\csname beginpicture\endcsname\relax
\expandafter\ifx\csname documentclass\endcsname\relax
\input pictex \else
\input prepictex \input pictex \input postpictex \fi\fi

\def\gt{{\mathsurround=0pt\it $\cal G\mskip-2mu$eometry \&\ 
$\cal T\!\!$opology}}        

\def\gtp{{\mathsurround=0pt\it $\cal G\mskip-2mu$eometry \&\ 
$\cal T\!\!$opology $\cal P\!$ublications}}  

\def\lognumber#1{\def\thelognumber{#1}}
\def\volumenumber#1{\def\thevolumenumber{#1}}
\def\papernumber#1{\def\thepapernumber{#1}}
\def\volumeyear#1{\def\thevolumeyear{#1}}

\def\pagenumbers#1#2{\def\startpage{#1}\def\finishpage{#2}}
\def\published#1{\def\publishdate{#1}}
\def\proposed#1{\def\theproposer{#1}}
\def\seconded#1{\def\theseconders{#1}}
\def\received#1{\def\receiveddate{#1}}

\def\accepted#1{\def\accepteddate{#1}}

\long\def\asciiabstract#1{\long\def\theasciiabstract{#1}}

\let\\\par\let\thelognumber\relax
\let\thevolumenumber\relax\let\thepapernumber\relax
\let\thevolumeyear\relax\let\thesamplenumber\relax\let\startpage\relax
\let\finishpage\relax\let\publishdate\relax\let\receiveddate\relax
\let\reviseddate\relax\let\accepteddate\relax\let\theasciititle\relax
\let\theasciiauthors\relax
\let\theasciiabstract\relax
\let\theasciiemail\relax\let\theshortauthors\relax\let\theshorttitle\relax

\long\def\maketitlep{   

\count0=\startpage

\gt\hfill      
\beginpicture
\setcoordinatesystem units <0.33truein, 0.33truein> point at 2.2 0.9
\setplotsymbol ({$\cal G$})
\plotsymbolspacing=9truept
\circulararc 315 degrees from 0 1 center at 0 0
\setplotsymbol ({$\cal T$})
\circulararc 315 degrees from 1 -1 center at 1 0
\endpicture
\break
{\small\ifx\thesamplenumber\relax 
Volume \else Sample
\fi\thevolumenumber\ (\thevolumeyear)
\startpage--\finishpage\nl
Published: \publishdate}
\vglue 0.5truein plus 0.4fil minus 0.1truein

{\parskip=0pt\leftskip 0pt plus 1fil\def\\{\par\smallskip}{\ifplaintex\large
\else\Large\fi\bf\thetitle}\par\medskip}   

\vglue 0pt plus 0.1fil 

{\parskip=0pt\leftskip 0pt plus 1fil\def\\{\par}{\sc\theauthors}
\par\medskip}

\vglue 0pt plus 0.1fil 

{\small\parskip=0pt\let\newline\\
{\leftskip 0pt plus 1fil\def\\{\par}{\sl\theaddress}\par}
\expandafter\ifx\theemail\relax    
\relax\else\vglue 5pt plus 0.02fil minus 2pt\def\\{\stdspace{\rm 
and}\stdspace} 
\cl{Email:\stdspace\tt\theemail}\fi
\ifx\theurl\relax                  
\relax\else\vglue 5pt plus 0.02fil minus 2pt\def\\{\stdspace{\rm 
and}\stdspace}
\cl{URL:\stdspace\tt\theurl}\fi\par}

\vglue 7pt plus 0.3fil minus 3pt

{\bf Abstract}
\vglue 5pt plus 0.1fil minus 2pt

\theabstract

\vglue 7pt plus 0.3fil minus 3pt

{\bf AMS Classification numbers}\quad Primary:\quad \theprimaryclass

Secondary:\quad \thesecondaryclass

\vglue 5pt plus 0.3fil minus 2pt

{\bf Keywords}\quad \thekeywords

\vglue 10pt plus 0.5fil minus 5pt

{\small  Proposed: \theproposer\hfill Received: \receiveddate\nl
Seconded: \theseconders\hfill 
\ifx\reviseddate\relax                         
Accepted: \accepteddate                        
\else
Revised: \reviseddate                          
\fi}
\eject
}       

\let\maketitlepage\maketitlep
\let\maketitle\maketitlepage

\font\phead=cmsl9 scaled 950
\font=cmsl9 scaled 1050
\font\pnum=cmbx10 scaled 913
\font\lnum=cmbx10 
\font\pfoot=cmsl9 scaled 950
\font=cmsl9 scaled 1050
\ifplaintex
\headline{\vbox to 0pt{\vskip -4.5mm\line{\small\phead\ifnum
\count0=\startpage ISSN 1364-0380 (on line)
1465-3060 (printed) \hfill {\pnum\folio}\else\ifodd\count0\def\\{ }%
\ifx\theshorttitle\relax\thetitle\else\theshorttitle\fi\hfill{\pnum\folio}
\else\def\\{ and }{\pnum\folio}\hfill\ifx\theshortauthors\relax\theauthors
\else\theshortauthors\fi\fi\fi}\vss}}
\footline{\vbox to 0pt{\vglue 0mm\line{\small\pfoot\ifnum\count0=\startpage
\copyright\ \gtp\hfill\else
\gt, Volume \thevolumenumber\ (\thevolumeyear)\hfill\fi}\vss
}}
\else
\makeatletter
\def\@oddhead{{\small\lhead\ifnum\count0=\startpage ISSN 1364-0380 (on line)
1465-3060 (printed) \hfill {\lnum\number\count0}\else\ifodd\count0
\def\\{ }\ifx\theshorttitle\relax \thetitle \else\theshorttitle\fi\hfill
{\lnum\number\count0}\else\def\\{ and }{\lnum\number\count0}
\hfill\ifx\theshortauthors\relax 
\theauthors\else\theshortauthors\fi\fi\fi}}\def\@evenhead{\@oddhead}
\def\@oddfoot{\small\lfoot\ifnum\count0=\startpage\copyright\ \gtp\hfill\else
\gt, Volume \thevolumenumber\ (\thevolumeyear)\hfill\fi}
\def\@evenfoot{\@oddfoot}
\makeatother
\fi

\newwrite\gtoutfile
\long\gdef\makeheadfile{  
{\def\\{, }\def\s{ }
\immediate\openout\gtoutfile head.xxx
\immediate\write\gtoutfile{To: math@arxiv.org}
\immediate\write\gtoutfile{Subject: put or rep NNNNN:pppp}
\immediate\write\gtoutfile{--text follows this line--}
\immediate\write\gtoutfile{Proxy-for: \ifx\theasciiauthors\relax
\theauthors\else\theasciiauthors\fi\s<\ifx\theasciiemail\relax\theemail\else\theasciiemail\fi>}
\immediate\write\gtoutfile{\noexpand\\}
\immediate\write\gtoutfile{Authors: \ifx\theasciiauthors\relax
\theauthors\else\theasciiauthors\fi}
\immediate\write\gtoutfile{Title: \ifx\theasciititle\relax
\thetitle\else\theasciititle\fi}
\immediate\write\gtoutfile{Subj-class: GT or SG or MG etc}
\immediate\write\gtoutfile{MSC-class: \theprimaryclass\ifx\thesecondaryclass\relax\else, \thesecondaryclass\fi}
\immediate\write\gtoutfile{Journal-ref: Geom. Topol. \thevolumenumber
(\thevolumeyear) \startpage-\finishpage}
\immediate\write\gtoutfile{Comments: Published by Geometry and Topology at}
\immediate\write\gtoutfile{\s\s http://www.maths.warwick.ac.uk/gt/GTVol\thevolumenumber/paper\thepapernumber.abs.html}
\immediate\write\gtoutfile{\noexpand\\}
\immediate\write\gtoutfile{}
\ifx\theasciiabstract\relax
\immediate\write\gtoutfile{\theabstract}\else
\immediate\write\gtoutfile{\theasciiabstract}\fi
\immediate\write\gtoutfile{}
\immediate\write\gtoutfile{\noexpand\\}
\immediate\write\gtoutfile{}
\immediate\closeout\gtoutfile}}  

\def\maketitlepage{\maketitlep\makeheadfile}
\let\maketitle\maketitlepage

\def\ifplaintex{\expandafter\ifx\csname documentclass\endcsname\relax}

\ifplaintex 
\else
\fi

\expandafter\ifx\csname beginpicture\endcsname\relax
\expandafter\ifx\csname documentclass\endcsname\relax
\input pictex \else
\input prepictex \input pictex \input postpictex \fi\fi

\def\gt{{\mathsurround=0pt\it $\cal G\mskip-2mu$eometry \&\ 
$\cal T\!\!$opology}}        

\def\gtp{{\mathsurround=0pt\it $\cal G\mskip-2mu$eometry \&\ 
$\cal T\!\!$opology $\cal P\!$ublications}}  

\def\lognumber#1{\def\thelognumber{#1}}
\def\volumenumber#1{\def\thevolumenumber{#1}}
\def\papernumber#1{\def\thepapernumber{#1}}
\def\volumeyear#1{\def\thevolumeyear{#1}}

\def\pagenumbers#1#2{\def\startpage{#1}\def\finishpage{#2}}
\def\published#1{\def\publishdate{#1}}
\def\proposed#1{\def\theproposer{#1}}
\def\seconded#1{\def\theseconders{#1}}
\def\received#1{\def\receiveddate{#1}}

\def\accepted#1{\def\accepteddate{#1}}

\long\def\asciiabstract#1{\long\def\theasciiabstract{#1}}

\let\\\par\let\thelognumber\relax
\let\thevolumenumber\relax\let\thepapernumber\relax
\let\thevolumeyear\relax\let\thesamplenumber\relax\let\startpage\relax
\let\finishpage\relax\let\publishdate\relax\let\receiveddate\relax
\let\reviseddate\relax\let\accepteddate\relax\let\theasciititle\relax
\let\theasciiauthors\relax
\let\theasciiabstract\relax
\let\theasciiemail\relax\let\theshortauthors\relax\let\theshorttitle\relax

\long\def\maketitlep{   

\count0=\startpage

\gt\hfill      
\beginpicture
\setcoordinatesystem units <0.33truein, 0.33truein> point at 2.2 0.9
\setplotsymbol ({$\cal G$})
\plotsymbolspacing=9truept
\circulararc 315 degrees from 0 1 center at 0 0
\setplotsymbol ({$\cal T$})
\circulararc 315 degrees from 1 -1 center at 1 0
\endpicture
\break
{\small\ifx\thesamplenumber\relax 
Volume \else Sample
\fi\thevolumenumber\ (\thevolumeyear)
\startpage--\finishpage\nl
Published: \publishdate}
\vglue 0.5truein plus 0.4fil minus 0.1truein

{\parskip=0pt\leftskip 0pt plus 1fil\def\\{\par\smallskip}{\ifplaintex\large
\else\Large\fi\bf\thetitle}\par\medskip}   

\vglue 0pt plus 0.1fil 

{\parskip=0pt\leftskip 0pt plus 1fil\def\\{\par}{\sc\theauthors}
\par\medskip}

\vglue 0pt plus 0.1fil 

{\small\parskip=0pt\let\newline\\
{\leftskip 0pt plus 1fil\def\\{\par}{\sl\theaddress}\par}
\expandafter\ifx\theemail\relax    
\relax\else\vglue 5pt plus 0.02fil minus 2pt\def\\{\stdspace{\rm 
and}\stdspace} 
\cl{Email:\stdspace\tt\theemail}\fi
\ifx\theurl\relax                  
\relax\else\vglue 5pt plus 0.02fil minus 2pt\def\\{\stdspace{\rm 
and}\stdspace}
\cl{URL:\stdspace\tt\theurl}\fi\par}

\vglue 7pt plus 0.3fil minus 3pt

{\bf Abstract}
\vglue 5pt plus 0.1fil minus 2pt

\theabstract

\vglue 7pt plus 0.3fil minus 3pt

{\bf AMS Classification numbers}\quad Primary:\quad \theprimaryclass

Secondary:\quad \thesecondaryclass

\vglue 5pt plus 0.3fil minus 2pt

{\bf Keywords}\quad \thekeywords

\vglue 10pt plus 0.5fil minus 5pt

{\small  Proposed: \theproposer\hfill Received: \receiveddate\nl
Seconded: \theseconders\hfill 
\ifx\reviseddate\relax                         
Accepted: \accepteddate                        
\else
Revised: \reviseddate                          
\fi}
\eject
}       

\let\maketitlepage\maketitlep
\let\maketitle\maketitlepage

\font\phead=cmsl9 scaled 950
\font=cmsl9 scaled 1050
\font\pnum=cmbx10 scaled 913
\font\lnum=cmbx10 
\font\pfoot=cmsl9 scaled 950
\font=cmsl9 scaled 1050
\ifplaintex
\headline{\vbox to 0pt{\vskip -4.5mm\line{\small\phead\ifnum
\count0=\startpage ISSN 1364-0380 (on line)
1465-3060 (printed) \hfill {\pnum\folio}\else\ifodd\count0\def\\{ }%
\ifx\theshorttitle\relax\thetitle\else\theshorttitle\fi\hfill{\pnum\folio}
\else\def\\{ and }{\pnum\folio}\hfill\ifx\theshortauthors\relax\theauthors
\else\theshortauthors\fi\fi\fi}\vss}}
\footline{\vbox to 0pt{\vglue 0mm\line{\small\pfoot\ifnum\count0=\startpage
\copyright\ \gtp\hfill\else
\gt, Volume \thevolumenumber\ (\thevolumeyear)\hfill\fi}\vss
}}
\else
\makeatletter
\def\@oddhead{{\small\lhead\ifnum\count0=\startpage ISSN 1364-0380 (on line)
1465-3060 (printed) \hfill {\lnum\number\count0}\else\ifodd\count0
\def\\{ }\ifx\theshorttitle\relax \thetitle \else\theshorttitle\fi\hfill
{\lnum\number\count0}\else\def\\{ and }{\lnum\number\count0}
\hfill\ifx\theshortauthors\relax 
\theauthors\else\theshortauthors\fi\fi\fi}}\def\@evenhead{\@oddhead}
\def\@oddfoot{\small\lfoot\ifnum\count0=\startpage\copyright\ \gtp\hfill\else
\gt, Volume \thevolumenumber\ (\thevolumeyear)\hfill\fi}
\def\@evenfoot{\@oddfoot}
\makeatother
\fi

\newwrite\gtoutfile
\long\gdef\makeheadfile{  
{\def\\{, }\def\s{ }
\immediate\openout\gtoutfile head.xxx
\immediate\write\gtoutfile{To: math@arxiv.org}
\immediate\write\gtoutfile{Subject: put or rep NNNNN:pppp}
\immediate\write\gtoutfile{--text follows this line--}
\immediate\write\gtoutfile{Proxy-for: \ifx\theasciiauthors\relax
\theauthors\else\theasciiauthors\fi\s<\ifx\theasciiemail\relax\theemail\else\theasciiemail\fi>}
\immediate\write\gtoutfile{\noexpand\\}
\immediate\write\gtoutfile{Authors: \ifx\theasciiauthors\relax
\theauthors\else\theasciiauthors\fi}
\immediate\write\gtoutfile{Title: \ifx\theasciititle\relax
\thetitle\else\theasciititle\fi}
\immediate\write\gtoutfile{Subj-class: GT or SG or MG etc}
\immediate\write\gtoutfile{MSC-class: \theprimaryclass\ifx\thesecondaryclass\relax\else, \thesecondaryclass\fi}
\immediate\write\gtoutfile{Journal-ref: Geom. Topol. \thevolumenumber
(\thevolumeyear) \startpage-\finishpage}
\immediate\write\gtoutfile{Comments: Published by Geometry and Topology at}
\immediate\write\gtoutfile{\s\s http://www.maths.warwick.ac.uk/gt/GTVol\thevolumenumber/paper\thepapernumber.abs.html}
\immediate\write\gtoutfile{\noexpand\\}
\immediate\write\gtoutfile{}
\ifx\theasciiabstract\relax
\immediate\write\gtoutfile{\theabstract}\else
\immediate\write\gtoutfile{\theasciiabstract}\fi
\immediate\write\gtoutfile{}
\immediate\write\gtoutfile{\noexpand\\}
\immediate\write\gtoutfile{}
\immediate\closeout\gtoutfile}}  

\def\maketitlepage{\maketitlep\makeheadfile}
\let\maketitle\maketitlepage

\def\ifplaintex{\expandafter\ifx\csname documentclass\endcsname\relax}

\ifplaintex 
\else
\fi

\expandafter\ifx\csname beginpicture\endcsname\relax
\expandafter\ifx\csname documentclass\endcsname\relax
\input pictex \else
\input prepictex \input pictex \input postpictex \fi\fi

\def\gt{{\mathsurround=0pt\it $\cal G\mskip-2mu$eometry \&\ 
$\cal T\!\!$opology}}        

\def\gtp{{\mathsurround=0pt\it $\cal G\mskip-2mu$eometry \&\ 
$\cal T\!\!$opology $\cal P\!$ublications}}  

\def\lognumber#1{\def\thelognumber{#1}}
\def\volumenumber#1{\def\thevolumenumber{#1}}
\def\papernumber#1{\def\thepapernumber{#1}}
\def\volumeyear#1{\def\thevolumeyear{#1}}

\def\pagenumbers#1#2{\def\startpage{#1}\def\finishpage{#2}}
\def\published#1{\def\publishdate{#1}}
\def\proposed#1{\def\theproposer{#1}}
\def\seconded#1{\def\theseconders{#1}}
\def\received#1{\def\receiveddate{#1}}

\def\accepted#1{\def\accepteddate{#1}}

\long\def\asciiabstract#1{\long\def\theasciiabstract{#1}}

\let\\\par\let\thelognumber\relax
\let\thevolumenumber\relax\let\thepapernumber\relax
\let\thevolumeyear\relax\let\thesamplenumber\relax\let\startpage\relax
\let\finishpage\relax\let\publishdate\relax\let\receiveddate\relax
\let\reviseddate\relax\let\accepteddate\relax\let\theasciititle\relax
\let\theasciiauthors\relax
\let\theasciiabstract\relax
\let\theasciiemail\relax\let\theshortauthors\relax\let\theshorttitle\relax

\long\def\maketitlep{   

\count0=\startpage

\gt\hfill      
\beginpicture
\setcoordinatesystem units <0.33truein, 0.33truein> point at 2.2 0.9
\setplotsymbol ({$\cal G$})
\plotsymbolspacing=9truept
\circulararc 315 degrees from 0 1 center at 0 0
\setplotsymbol ({$\cal T$})
\circulararc 315 degrees from 1 -1 center at 1 0
\endpicture
\break
{\small\ifx\thesamplenumber\relax 
Volume \else Sample
\fi\thevolumenumber\ (\thevolumeyear)
\startpage--\finishpage\nl
Published: \publishdate}
\vglue 0.5truein plus 0.4fil minus 0.1truein

{\parskip=0pt\leftskip 0pt plus 1fil\def\\{\par\smallskip}{\ifplaintex\large
\else\Large\fi\bf\thetitle}\par\medskip}   

\vglue 0pt plus 0.1fil 

{\parskip=0pt\leftskip 0pt plus 1fil\def\\{\par}{\sc\theauthors}
\par\medskip}

\vglue 0pt plus 0.1fil 

{\small\parskip=0pt\let\newline\\
{\leftskip 0pt plus 1fil\def\\{\par}{\sl\theaddress}\par}
\expandafter\ifx\theemail\relax    
\relax\else\vglue 5pt plus 0.02fil minus 2pt\def\\{\stdspace{\rm 
and}\stdspace} 
\cl{Email:\stdspace\tt\theemail}\fi
\ifx\theurl\relax                  
\relax\else\vglue 5pt plus 0.02fil minus 2pt\def\\{\stdspace{\rm 
and}\stdspace}
\cl{URL:\stdspace\tt\theurl}\fi\par}

\vglue 7pt plus 0.3fil minus 3pt

{\bf Abstract}
\vglue 5pt plus 0.1fil minus 2pt

\theabstract

\vglue 7pt plus 0.3fil minus 3pt

{\bf AMS Classification numbers}\quad Primary:\quad \theprimaryclass

Secondary:\quad \thesecondaryclass

\vglue 5pt plus 0.3fil minus 2pt

{\bf Keywords}\quad \thekeywords

\vglue 10pt plus 0.5fil minus 5pt

{\small  Proposed: \theproposer\hfill Received: \receiveddate\nl
Seconded: \theseconders\hfill 
\ifx\reviseddate\relax                         
Accepted: \accepteddate                        
\else
Revised: \reviseddate                          
\fi}
\eject
}       

\let\maketitlepage\maketitlep
\let\maketitle\maketitlepage

\font\phead=cmsl9 scaled 950
\font=cmsl9 scaled 1050
\font\pnum=cmbx10 scaled 913
\font\lnum=cmbx10 
\font\pfoot=cmsl9 scaled 950
\font=cmsl9 scaled 1050
\ifplaintex
\headline{\vbox to 0pt{\vskip -4.5mm\line{\small\phead\ifnum
\count0=\startpage ISSN 1364-0380 (on line)
1465-3060 (printed) \hfill {\pnum\folio}\else\ifodd\count0\def\\{ }%
\ifx\theshorttitle\relax\thetitle\else\theshorttitle\fi\hfill{\pnum\folio}
\else\def\\{ and }{\pnum\folio}\hfill\ifx\theshortauthors\relax\theauthors
\else\theshortauthors\fi\fi\fi}\vss}}
\footline{\vbox to 0pt{\vglue 0mm\line{\small\pfoot\ifnum\count0=\startpage
\copyright\ \gtp\hfill\else
\gt, Volume \thevolumenumber\ (\thevolumeyear)\hfill\fi}\vss
}}
\else
\makeatletter
\def\@oddhead{{\small\lhead\ifnum\count0=\startpage ISSN 1364-0380 (on line)
1465-3060 (printed) \hfill {\lnum\number\count0}\else\ifodd\count0
\def\\{ }\ifx\theshorttitle\relax \thetitle \else\theshorttitle\fi\hfill
{\lnum\number\count0}\else\def\\{ and }{\lnum\number\count0}
\hfill\ifx\theshortauthors\relax 
\theauthors\else\theshortauthors\fi\fi\fi}}\def\@evenhead{\@oddhead}
\def\@oddfoot{\small\lfoot\ifnum\count0=\startpage\copyright\ \gtp\hfill\else
\gt, Volume \thevolumenumber\ (\thevolumeyear)\hfill\fi}
\def\@evenfoot{\@oddfoot}
\makeatother
\fi

\newwrite\gtoutfile
\long\gdef\makeheadfile{  
{\def\\{, }\def\s{ }
\immediate\openout\gtoutfile head.xxx
\immediate\write\gtoutfile{To: math@arxiv.org}
\immediate\write\gtoutfile{Subject: put or rep NNNNN:pppp}
\immediate\write\gtoutfile{--text follows this line--}
\immediate\write\gtoutfile{Proxy-for: \ifx\theasciiauthors\relax
\theauthors\else\theasciiauthors\fi\s<\ifx\theasciiemail\relax\theemail\else\theasciiemail\fi>}
\immediate\write\gtoutfile{\noexpand\\}
\immediate\write\gtoutfile{Authors: \ifx\theasciiauthors\relax
\theauthors\else\theasciiauthors\fi}
\immediate\write\gtoutfile{Title: \ifx\theasciititle\relax
\thetitle\else\theasciititle\fi}
\immediate\write\gtoutfile{Subj-class: GT or SG or MG etc}
\immediate\write\gtoutfile{MSC-class: \theprimaryclass\ifx\thesecondaryclass\relax\else, \thesecondaryclass\fi}
\immediate\write\gtoutfile{Journal-ref: Geom. Topol. \thevolumenumber
(\thevolumeyear) \startpage-\finishpage}
\immediate\write\gtoutfile{Comments: Published by Geometry and Topology at}
\immediate\write\gtoutfile{\s\s http://www.maths.warwick.ac.uk/gt/GTVol\thevolumenumber/paper\thepapernumber.abs.html}
\immediate\write\gtoutfile{\noexpand\\}
\immediate\write\gtoutfile{}
\ifx\theasciiabstract\relax
\immediate\write\gtoutfile{\theabstract}\else
\immediate\write\gtoutfile{\theasciiabstract}\fi
\immediate\write\gtoutfile{}
\immediate\write\gtoutfile{\noexpand\\}
\immediate\write\gtoutfile{}
\immediate\closeout\gtoutfile}}  

\def\maketitlepage{\maketitlep\makeheadfile}
\let\maketitle\maketitlepage

\def\ifplaintex{\expandafter\ifx\csname documentclass\endcsname\relax}

\ifplaintex 
\else
\fi

\expandafter\ifx\csname beginpicture\endcsname\relax
\expandafter\ifx\csname documentclass\endcsname\relax
\input pictex \else
\input prepictex \input pictex \input postpictex \fi\fi

\def\gt{{\mathsurround=0pt\it $\cal G\mskip-2mu$eometry \&\ 
$\cal T\!\!$opology}}        

\def\gtp{{\mathsurround=0pt\it $\cal G\mskip-2mu$eometry \&\ 
$\cal T\!\!$opology $\cal P\!$ublications}}  

\def\lognumber#1{\def\thelognumber{#1}}
\def\volumenumber#1{\def\thevolumenumber{#1}}
\def\papernumber#1{\def\thepapernumber{#1}}
\def\volumeyear#1{\def\thevolumeyear{#1}}

\def\pagenumbers#1#2{\def\startpage{#1}\def\finishpage{#2}}
\def\published#1{\def\publishdate{#1}}
\def\proposed#1{\def\theproposer{#1}}
\def\seconded#1{\def\theseconders{#1}}
\def\received#1{\def\receiveddate{#1}}

\def\accepted#1{\def\accepteddate{#1}}

\long\def\asciiabstract#1{\long\def\theasciiabstract{#1}}

\let\\\par\let\thelognumber\relax
\let\thevolumenumber\relax\let\thepapernumber\relax
\let\thevolumeyear\relax\let\thesamplenumber\relax\let\startpage\relax
\let\finishpage\relax\let\publishdate\relax\let\receiveddate\relax
\let\reviseddate\relax\let\accepteddate\relax\let\theasciititle\relax
\let\theasciiauthors\relax
\let\theasciiabstract\relax
\let\theasciiemail\relax\let\theshortauthors\relax\let\theshorttitle\relax

\long\def\maketitlep{   

\count0=\startpage

\gt\hfill      
\beginpicture
\setcoordinatesystem units <0.33truein, 0.33truein> point at 2.2 0.9
\setplotsymbol ({$\cal G$})
\plotsymbolspacing=9truept
\circulararc 315 degrees from 0 1 center at 0 0
\setplotsymbol ({$\cal T$})
\circulararc 315 degrees from 1 -1 center at 1 0
\endpicture
\break
{\small\ifx\thesamplenumber\relax 
Volume \else Sample
\fi\thevolumenumber\ (\thevolumeyear)
\startpage--\finishpage\nl
Published: \publishdate}
\vglue 0.5truein plus 0.4fil minus 0.1truein

{\parskip=0pt\leftskip 0pt plus 1fil\def\\{\par\smallskip}{\ifplaintex\large
\else\Large\fi\bf\thetitle}\par\medskip}   

\vglue 0pt plus 0.1fil 

{\parskip=0pt\leftskip 0pt plus 1fil\def\\{\par}{\sc\theauthors}
\par\medskip}

\vglue 0pt plus 0.1fil 

{\small\parskip=0pt\let\newline\\
{\leftskip 0pt plus 1fil\def\\{\par}{\sl\theaddress}\par}
\expandafter\ifx\theemail\relax    
\relax\else\vglue 5pt plus 0.02fil minus 2pt\def\\{\stdspace{\rm 
and}\stdspace} 
\cl{Email:\stdspace\tt\theemail}\fi
\ifx\theurl\relax                  
\relax\else\vglue 5pt plus 0.02fil minus 2pt\def\\{\stdspace{\rm 
and}\stdspace}
\cl{URL:\stdspace\tt\theurl}\fi\par}

\vglue 7pt plus 0.3fil minus 3pt

{\bf Abstract}
\vglue 5pt plus 0.1fil minus 2pt

\theabstract

\vglue 7pt plus 0.3fil minus 3pt

{\bf AMS Classification numbers}\quad Primary:\quad \theprimaryclass

Secondary:\quad \thesecondaryclass

\vglue 5pt plus 0.3fil minus 2pt

{\bf Keywords}\quad \thekeywords

\vglue 10pt plus 0.5fil minus 5pt

{\small  Proposed: \theproposer\hfill Received: \receiveddate\nl
Seconded: \theseconders\hfill 
\ifx\reviseddate\relax                         
Accepted: \accepteddate                        
\else
Revised: \reviseddate                          
\fi}
\eject
}       

\let\maketitlepage\maketitlep
\let\maketitle\maketitlepage

\font\phead=cmsl9 scaled 950
\font=cmsl9 scaled 1050
\font\pnum=cmbx10 scaled 913
\font\lnum=cmbx10 
\font\pfoot=cmsl9 scaled 950
\font=cmsl9 scaled 1050
\ifplaintex
\headline{\vbox to 0pt{\vskip -4.5mm\line{\small\phead\ifnum
\count0=\startpage ISSN 1364-0380 (on line)
1465-3060 (printed) \hfill {\pnum\folio}\else\ifodd\count0\def\\{ }%
\ifx\theshorttitle\relax\thetitle\else\theshorttitle\fi\hfill{\pnum\folio}
\else\def\\{ and }{\pnum\folio}\hfill\ifx\theshortauthors\relax\theauthors
\else\theshortauthors\fi\fi\fi}\vss}}
\footline{\vbox to 0pt{\vglue 0mm\line{\small\pfoot\ifnum\count0=\startpage
\copyright\ \gtp\hfill\else
\gt, Volume \thevolumenumber\ (\thevolumeyear)\hfill\fi}\vss
}}
\else
\makeatletter
\def\@oddhead{{\small\lhead\ifnum\count0=\startpage ISSN 1364-0380 (on line)
1465-3060 (printed) \hfill {\lnum\number\count0}\else\ifodd\count0
\def\\{ }\ifx\theshorttitle\relax \thetitle \else\theshorttitle\fi\hfill
{\lnum\number\count0}\else\def\\{ and }{\lnum\number\count0}
\hfill\ifx\theshortauthors\relax 
\theauthors\else\theshortauthors\fi\fi\fi}}\def\@evenhead{\@oddhead}
\def\@oddfoot{\small\lfoot\ifnum\count0=\startpage\copyright\ \gtp\hfill\else
\gt, Volume \thevolumenumber\ (\thevolumeyear)\hfill\fi}
\def\@evenfoot{\@oddfoot}
\makeatother
\fi

\newwrite\gtoutfile
\long\gdef\makeheadfile{  
{\def\\{, }\def\s{ }
\immediate\openout\gtoutfile head.xxx
\immediate\write\gtoutfile{To: math@arxiv.org}
\immediate\write\gtoutfile{Subject: put or rep NNNNN:pppp}
\immediate\write\gtoutfile{--text follows this line--}
\immediate\write\gtoutfile{Proxy-for: \ifx\theasciiauthors\relax
\theauthors\else\theasciiauthors\fi\s<\ifx\theasciiemail\relax\theemail\else\theasciiemail\fi>}
\immediate\write\gtoutfile{\noexpand\\}
\immediate\write\gtoutfile{Authors: \ifx\theasciiauthors\relax
\theauthors\else\theasciiauthors\fi}
\immediate\write\gtoutfile{Title: \ifx\theasciititle\relax
\thetitle\else\theasciititle\fi}
\immediate\write\gtoutfile{Subj-class: GT or SG or MG etc}
\immediate\write\gtoutfile{MSC-class: \theprimaryclass\ifx\thesecondaryclass\relax\else, \thesecondaryclass\fi}
\immediate\write\gtoutfile{Journal-ref: Geom. Topol. \thevolumenumber
(\thevolumeyear) \startpage-\finishpage}
\immediate\write\gtoutfile{Comments: Published by Geometry and Topology at}
\immediate\write\gtoutfile{\s\s http://www.maths.warwick.ac.uk/gt/GTVol\thevolumenumber/paper\thepapernumber.abs.html}
\immediate\write\gtoutfile{\noexpand\\}
\immediate\write\gtoutfile{}
\ifx\theasciiabstract\relax
\immediate\write\gtoutfile{\theabstract}\else
\immediate\write\gtoutfile{\theasciiabstract}\fi
\immediate\write\gtoutfile{}
\immediate\write\gtoutfile{\noexpand\\}
\immediate\write\gtoutfile{}
\immediate\closeout\gtoutfile}}  

\def\maketitlepage{\maketitlep\makeheadfile}
\let\maketitle\maketitlepage

\def\ifplaintex{\expandafter\ifx\csname documentclass\endcsname\relax}

\ifplaintex 
\else
\fi

\expandafter\ifx\csname beginpicture\endcsname\relax
\expandafter\ifx\csname documentclass\endcsname\relax
\input pictex \else
\input prepictex \input pictex \input postpictex \fi\fi

\def\gt{{\mathsurround=0pt\it $\cal G\mskip-2mu$eometry \&\ 
$\cal T\!\!$opology}}        

\def\gtp{{\mathsurround=0pt\it $\cal G\mskip-2mu$eometry \&\ 
$\cal T\!\!$opology $\cal P\!$ublications}}  

\def\lognumber#1{\def\thelognumber{#1}}
\def\volumenumber#1{\def\thevolumenumber{#1}}
\def\papernumber#1{\def\thepapernumber{#1}}
\def\volumeyear#1{\def\thevolumeyear{#1}}

\def\pagenumbers#1#2{\def\startpage{#1}\def\finishpage{#2}}
\def\published#1{\def\publishdate{#1}}
\def\proposed#1{\def\theproposer{#1}}
\def\seconded#1{\def\theseconders{#1}}
\def\received#1{\def\receiveddate{#1}}

\def\accepted#1{\def\accepteddate{#1}}

\long\def\asciiabstract#1{\long\def\theasciiabstract{#1}}

\let\\\par\let\thelognumber\relax
\let\thevolumenumber\relax\let\thepapernumber\relax
\let\thevolumeyear\relax\let\thesamplenumber\relax\let\startpage\relax
\let\finishpage\relax\let\publishdate\relax\let\receiveddate\relax
\let\reviseddate\relax\let\accepteddate\relax\let\theasciititle\relax
\let\theasciiauthors\relax
\let\theasciiabstract\relax
\let\theasciiemail\relax\let\theshortauthors\relax\let\theshorttitle\relax

\long\def\maketitlep{   

\count0=\startpage

\gt\hfill      
\beginpicture
\setcoordinatesystem units <0.33truein, 0.33truein> point at 2.2 0.9
\setplotsymbol ({$\cal G$})
\plotsymbolspacing=9truept
\circulararc 315 degrees from 0 1 center at 0 0
\setplotsymbol ({$\cal T$})
\circulararc 315 degrees from 1 -1 center at 1 0
\endpicture
\break
{\small\ifx\thesamplenumber\relax 
Volume \else Sample
\fi\thevolumenumber\ (\thevolumeyear)
\startpage--\finishpage\nl
Published: \publishdate}
\vglue 0.5truein plus 0.4fil minus 0.1truein

{\parskip=0pt\leftskip 0pt plus 1fil\def\\{\par\smallskip}{\ifplaintex\large
\else\Large\fi\bf\thetitle}\par\medskip}   

\vglue 0pt plus 0.1fil 

{\parskip=0pt\leftskip 0pt plus 1fil\def\\{\par}{\sc\theauthors}
\par\medskip}

\vglue 0pt plus 0.1fil 

{\small\parskip=0pt\let\newline\\
{\leftskip 0pt plus 1fil\def\\{\par}{\sl\theaddress}\par}
\expandafter\ifx\theemail\relax    
\relax\else\vglue 5pt plus 0.02fil minus 2pt\def\\{\stdspace{\rm 
and}\stdspace} 
\cl{Email:\stdspace\tt\theemail}\fi
\ifx\theurl\relax                  
\relax\else\vglue 5pt plus 0.02fil minus 2pt\def\\{\stdspace{\rm 
and}\stdspace}
\cl{URL:\stdspace\tt\theurl}\fi\par}

\vglue 7pt plus 0.3fil minus 3pt

{\bf Abstract}
\vglue 5pt plus 0.1fil minus 2pt

\theabstract

\vglue 7pt plus 0.3fil minus 3pt

{\bf AMS Classification numbers}\quad Primary:\quad \theprimaryclass

Secondary:\quad \thesecondaryclass

\vglue 5pt plus 0.3fil minus 2pt

{\bf Keywords}\quad \thekeywords

\vglue 10pt plus 0.5fil minus 5pt

{\small  Proposed: \theproposer\hfill Received: \receiveddate\nl
Seconded: \theseconders\hfill 
\ifx\reviseddate\relax                         
Accepted: \accepteddate                        
\else
Revised: \reviseddate                          
\fi}
\eject
}       

\let\maketitlepage\maketitlep
\let\maketitle\maketitlepage

\font\phead=cmsl9 scaled 950
\font=cmsl9 scaled 1050
\font\pnum=cmbx10 scaled 913
\font\lnum=cmbx10 
\font\pfoot=cmsl9 scaled 950
\font=cmsl9 scaled 1050
\ifplaintex
\headline{\vbox to 0pt{\vskip -4.5mm\line{\small\phead\ifnum
\count0=\startpage ISSN 1364-0380 (on line)
1465-3060 (printed) \hfill {\pnum\folio}\else\ifodd\count0\def\\{ }%
\ifx\theshorttitle\relax\thetitle\else\theshorttitle\fi\hfill{\pnum\folio}
\else\def\\{ and }{\pnum\folio}\hfill\ifx\theshortauthors\relax\theauthors
\else\theshortauthors\fi\fi\fi}\vss}}
\footline{\vbox to 0pt{\vglue 0mm\line{\small\pfoot\ifnum\count0=\startpage
\copyright\ \gtp\hfill\else
\gt, Volume \thevolumenumber\ (\thevolumeyear)\hfill\fi}\vss
}}
\else
\makeatletter
\def\@oddhead{{\small\lhead\ifnum\count0=\startpage ISSN 1364-0380 (on line)
1465-3060 (printed) \hfill {\lnum\number\count0}\else\ifodd\count0
\def\\{ }\ifx\theshorttitle\relax \thetitle \else\theshorttitle\fi\hfill
{\lnum\number\count0}\else\def\\{ and }{\lnum\number\count0}
\hfill\ifx\theshortauthors\relax 
\theauthors\else\theshortauthors\fi\fi\fi}}\def\@evenhead{\@oddhead}
\def\@oddfoot{\small\lfoot\ifnum\count0=\startpage\copyright\ \gtp\hfill\else
\gt, Volume \thevolumenumber\ (\thevolumeyear)\hfill\fi}
\def\@evenfoot{\@oddfoot}
\makeatother
\fi

\newwrite\gtoutfile
\long\gdef\makeheadfile{  
{\def\\{, }\def\s{ }
\immediate\openout\gtoutfile head.xxx
\immediate\write\gtoutfile{To: math@arxiv.org}
\immediate\write\gtoutfile{Subject: put or rep NNNNN:pppp}
\immediate\write\gtoutfile{--text follows this line--}
\immediate\write\gtoutfile{Proxy-for: \ifx\theasciiauthors\relax
\theauthors\else\theasciiauthors\fi\s<\ifx\theasciiemail\relax\theemail\else\theasciiemail\fi>}
\immediate\write\gtoutfile{\noexpand\\}
\immediate\write\gtoutfile{Authors: \ifx\theasciiauthors\relax
\theauthors\else\theasciiauthors\fi}
\immediate\write\gtoutfile{Title: \ifx\theasciititle\relax
\thetitle\else\theasciititle\fi}
\immediate\write\gtoutfile{Subj-class: GT or SG or MG etc}
\immediate\write\gtoutfile{MSC-class: \theprimaryclass\ifx\thesecondaryclass\relax\else, \thesecondaryclass\fi}
\immediate\write\gtoutfile{Journal-ref: Geom. Topol. \thevolumenumber
(\thevolumeyear) \startpage-\finishpage}
\immediate\write\gtoutfile{Comments: Published by Geometry and Topology at}
\immediate\write\gtoutfile{\s\s http://www.maths.warwick.ac.uk/gt/GTVol\thevolumenumber/paper\thepapernumber.abs.html}
\immediate\write\gtoutfile{\noexpand\\}
\immediate\write\gtoutfile{}
\ifx\theasciiabstract\relax
\immediate\write\gtoutfile{\theabstract}\else
\immediate\write\gtoutfile{\theasciiabstract}\fi
\immediate\write\gtoutfile{}
\immediate\write\gtoutfile{\noexpand\\}
\immediate\write\gtoutfile{}
\immediate\closeout\gtoutfile}}  

\def\maketitlepage{\maketitlep\makeheadfile}
\let\maketitle\maketitlepage

\def\ifplaintex{\expandafter\ifx\csname documentclass\endcsname\relax}

\ifplaintex 
\else
\fi

\expandafter\ifx\csname beginpicture\endcsname\relax
\expandafter\ifx\csname documentclass\endcsname\relax
\input pictex \else
\input prepictex \input pictex \input postpictex \fi\fi

\def\gt{{\mathsurround=0pt\it $\cal G\mskip-2mu$eometry \&\ 
$\cal T\!\!$opology}}        

\def\gtp{{\mathsurround=0pt\it $\cal G\mskip-2mu$eometry \&\ 
$\cal T\!\!$opology $\cal P\!$ublications}}  

\def\lognumber#1{\def\thelognumber{#1}}
\def\volumenumber#1{\def\thevolumenumber{#1}}
\def\papernumber#1{\def\thepapernumber{#1}}
\def\volumeyear#1{\def\thevolumeyear{#1}}

\def\pagenumbers#1#2{\def\startpage{#1}\def\finishpage{#2}}
\def\published#1{\def\publishdate{#1}}
\def\proposed#1{\def\theproposer{#1}}
\def\seconded#1{\def\theseconders{#1}}
\def\received#1{\def\receiveddate{#1}}

\def\accepted#1{\def\accepteddate{#1}}

\long\def\asciiabstract#1{\long\def\theasciiabstract{#1}}

\let\\\par\let\thelognumber\relax
\let\thevolumenumber\relax\let\thepapernumber\relax
\let\thevolumeyear\relax\let\thesamplenumber\relax\let\startpage\relax
\let\finishpage\relax\let\publishdate\relax\let\receiveddate\relax
\let\reviseddate\relax\let\accepteddate\relax\let\theasciititle\relax
\let\theasciiauthors\relax
\let\theasciiabstract\relax
\let\theasciiemail\relax\let\theshortauthors\relax\let\theshorttitle\relax

\long\def\maketitlep{   

\count0=\startpage

\gt\hfill      
\beginpicture
\setcoordinatesystem units <0.33truein, 0.33truein> point at 2.2 0.9
\setplotsymbol ({$\cal G$})
\plotsymbolspacing=9truept
\circulararc 315 degrees from 0 1 center at 0 0
\setplotsymbol ({$\cal T$})
\circulararc 315 degrees from 1 -1 center at 1 0
\endpicture
\break
{\small\ifx\thesamplenumber\relax 
Volume \else Sample
\fi\thevolumenumber\ (\thevolumeyear)
\startpage--\finishpage\nl
Published: \publishdate}
\vglue 0.5truein plus 0.4fil minus 0.1truein

{\parskip=0pt\leftskip 0pt plus 1fil\def\\{\par\smallskip}{\ifplaintex\large
\else\Large\fi\bf\thetitle}\par\medskip}   

\vglue 0pt plus 0.1fil 

{\parskip=0pt\leftskip 0pt plus 1fil\def\\{\par}{\sc\theauthors}
\par\medskip}

\vglue 0pt plus 0.1fil 

{\small\parskip=0pt\let\newline\\
{\leftskip 0pt plus 1fil\def\\{\par}{\sl\theaddress}\par}
\expandafter\ifx\theemail\relax    
\relax\else\vglue 5pt plus 0.02fil minus 2pt\def\\{\stdspace{\rm 
and}\stdspace} 
\cl{Email:\stdspace\tt\theemail}\fi
\ifx\theurl\relax                  
\relax\else\vglue 5pt plus 0.02fil minus 2pt\def\\{\stdspace{\rm 
and}\stdspace}
\cl{URL:\stdspace\tt\theurl}\fi\par}

\vglue 7pt plus 0.3fil minus 3pt

{\bf Abstract}
\vglue 5pt plus 0.1fil minus 2pt

\theabstract

\vglue 7pt plus 0.3fil minus 3pt

{\bf AMS Classification numbers}\quad Primary:\quad \theprimaryclass

Secondary:\quad \thesecondaryclass

\vglue 5pt plus 0.3fil minus 2pt

{\bf Keywords}\quad \thekeywords

\vglue 10pt plus 0.5fil minus 5pt

{\small  Proposed: \theproposer\hfill Received: \receiveddate\nl
Seconded: \theseconders\hfill 
\ifx\reviseddate\relax                         
Accepted: \accepteddate                        
\else
Revised: \reviseddate                          
\fi}
\eject
}       

\let\maketitlepage\maketitlep
\let\maketitle\maketitlepage

\font\phead=cmsl9 scaled 950
\font=cmsl9 scaled 1050
\font\pnum=cmbx10 scaled 913
\font\lnum=cmbx10 
\font\pfoot=cmsl9 scaled 950
\font=cmsl9 scaled 1050
\ifplaintex
\headline{\vbox to 0pt{\vskip -4.5mm\line{\small\phead\ifnum
\count0=\startpage ISSN 1364-0380 (on line)
1465-3060 (printed) \hfill {\pnum\folio}\else\ifodd\count0\def\\{ }%
\ifx\theshorttitle\relax\thetitle\else\theshorttitle\fi\hfill{\pnum\folio}
\else\def\\{ and }{\pnum\folio}\hfill\ifx\theshortauthors\relax\theauthors
\else\theshortauthors\fi\fi\fi}\vss}}
\footline{\vbox to 0pt{\vglue 0mm\line{\small\pfoot\ifnum\count0=\startpage
\copyright\ \gtp\hfill\else
\gt, Volume \thevolumenumber\ (\thevolumeyear)\hfill\fi}\vss
}}
\else
\makeatletter
\def\@oddhead{{\small\lhead\ifnum\count0=\startpage ISSN 1364-0380 (on line)
1465-3060 (printed) \hfill {\lnum\number\count0}\else\ifodd\count0
\def\\{ }\ifx\theshorttitle\relax \thetitle \else\theshorttitle\fi\hfill
{\lnum\number\count0}\else\def\\{ and }{\lnum\number\count0}
\hfill\ifx\theshortauthors\relax 
\theauthors\else\theshortauthors\fi\fi\fi}}\def\@evenhead{\@oddhead}
\def\@oddfoot{\small\lfoot\ifnum\count0=\startpage\copyright\ \gtp\hfill\else
\gt, Volume \thevolumenumber\ (\thevolumeyear)\hfill\fi}
\def\@evenfoot{\@oddfoot}
\makeatother
\fi

\newwrite\gtoutfile
\long\gdef\makeheadfile{  
{\def\\{, }\def\s{ }
\immediate\openout\gtoutfile head.xxx
\immediate\write\gtoutfile{To: math@arxiv.org}
\immediate\write\gtoutfile{Subject: put or rep NNNNN:pppp}
\immediate\write\gtoutfile{--text follows this line--}
\immediate\write\gtoutfile{Proxy-for: \ifx\theasciiauthors\relax
\theauthors\else\theasciiauthors\fi\s<\ifx\theasciiemail\relax\theemail\else\theasciiemail\fi>}
\immediate\write\gtoutfile{\noexpand\\}
\immediate\write\gtoutfile{Authors: \ifx\theasciiauthors\relax
\theauthors\else\theasciiauthors\fi}
\immediate\write\gtoutfile{Title: \ifx\theasciititle\relax
\thetitle\else\theasciititle\fi}
\immediate\write\gtoutfile{Subj-class: GT or SG or MG etc}
\immediate\write\gtoutfile{MSC-class: \theprimaryclass\ifx\thesecondaryclass\relax\else, \thesecondaryclass\fi}
\immediate\write\gtoutfile{Journal-ref: Geom. Topol. \thevolumenumber
(\thevolumeyear) \startpage-\finishpage}
\immediate\write\gtoutfile{Comments: Published by Geometry and Topology at}
\immediate\write\gtoutfile{\s\s http://www.maths.warwick.ac.uk/gt/GTVol\thevolumenumber/paper\thepapernumber.abs.html}
\immediate\write\gtoutfile{\noexpand\\}
\immediate\write\gtoutfile{}
\ifx\theasciiabstract\relax
\immediate\write\gtoutfile{\theabstract}\else
\immediate\write\gtoutfile{\theasciiabstract}\fi
\immediate\write\gtoutfile{}
\immediate\write\gtoutfile{\noexpand\\}
\immediate\write\gtoutfile{}
\immediate\closeout\gtoutfile}}  

\def\maketitlepage{\maketitlep\makeheadfile}
\let\maketitle\maketitlepage

\def\ifplaintex{\expandafter\ifx\csname documentclass\endcsname\relax}

\ifplaintex 
\else
\fi

\expandafter\ifx\csname beginpicture\endcsname\relax
\expandafter\ifx\csname documentclass\endcsname\relax
\input pictex \else
\input prepictex \input pictex \input postpictex \fi\fi

\def\gt{{\mathsurround=0pt\it $\cal G\mskip-2mu$eometry \&\ 
$\cal T\!\!$opology}}        

\def\gtp{{\mathsurround=0pt\it $\cal G\mskip-2mu$eometry \&\ 
$\cal T\!\!$opology $\cal P\!$ublications}}  

\def\lognumber#1{\def\thelognumber{#1}}
\def\volumenumber#1{\def\thevolumenumber{#1}}
\def\papernumber#1{\def\thepapernumber{#1}}
\def\volumeyear#1{\def\thevolumeyear{#1}}

\def\pagenumbers#1#2{\def\startpage{#1}\def\finishpage{#2}}
\def\published#1{\def\publishdate{#1}}
\def\proposed#1{\def\theproposer{#1}}
\def\seconded#1{\def\theseconders{#1}}
\def\received#1{\def\receiveddate{#1}}

\def\accepted#1{\def\accepteddate{#1}}

\long\def\asciiabstract#1{\long\def\theasciiabstract{#1}}

\let\\\par\let\thelognumber\relax
\let\thevolumenumber\relax\let\thepapernumber\relax
\let\thevolumeyear\relax\let\thesamplenumber\relax\let\startpage\relax
\let\finishpage\relax\let\publishdate\relax\let\receiveddate\relax
\let\reviseddate\relax\let\accepteddate\relax\let\theasciititle\relax
\let\theasciiauthors\relax
\let\theasciiabstract\relax
\let\theasciiemail\relax\let\theshortauthors\relax\let\theshorttitle\relax

\long\def\maketitlep{   

\count0=\startpage

\gt\hfill      
\beginpicture
\setcoordinatesystem units <0.33truein, 0.33truein> point at 2.2 0.9
\setplotsymbol ({$\cal G$})
\plotsymbolspacing=9truept
\circulararc 315 degrees from 0 1 center at 0 0
\setplotsymbol ({$\cal T$})
\circulararc 315 degrees from 1 -1 center at 1 0
\endpicture
\break
{\small\ifx\thesamplenumber\relax 
Volume \else Sample
\fi\thevolumenumber\ (\thevolumeyear)
\startpage--\finishpage\nl
Published: \publishdate}
\vglue 0.5truein plus 0.4fil minus 0.1truein

{\parskip=0pt\leftskip 0pt plus 1fil\def\\{\par\smallskip}{\ifplaintex\large
\else\Large\fi\bf\thetitle}\par\medskip}   

\vglue 0pt plus 0.1fil 

{\parskip=0pt\leftskip 0pt plus 1fil\def\\{\par}{\sc\theauthors}
\par\medskip}

\vglue 0pt plus 0.1fil 

{\small\parskip=0pt\let\newline\\
{\leftskip 0pt plus 1fil\def\\{\par}{\sl\theaddress}\par}
\expandafter\ifx\theemail\relax    
\relax\else\vglue 5pt plus 0.02fil minus 2pt\def\\{\stdspace{\rm 
and}\stdspace} 
\cl{Email:\stdspace\tt\theemail}\fi
\ifx\theurl\relax                  
\relax\else\vglue 5pt plus 0.02fil minus 2pt\def\\{\stdspace{\rm 
and}\stdspace}
\cl{URL:\stdspace\tt\theurl}\fi\par}

\vglue 7pt plus 0.3fil minus 3pt

{\bf Abstract}
\vglue 5pt plus 0.1fil minus 2pt

\theabstract

\vglue 7pt plus 0.3fil minus 3pt

{\bf AMS Classification numbers}\quad Primary:\quad \theprimaryclass

Secondary:\quad \thesecondaryclass

\vglue 5pt plus 0.3fil minus 2pt

{\bf Keywords}\quad \thekeywords

\vglue 10pt plus 0.5fil minus 5pt

{\small  Proposed: \theproposer\hfill Received: \receiveddate\nl
Seconded: \theseconders\hfill 
\ifx\reviseddate\relax                         
Accepted: \accepteddate                        
\else
Revised: \reviseddate                          
\fi}
\eject
}       

\let\maketitlepage\maketitlep
\let\maketitle\maketitlepage

\font\phead=cmsl9 scaled 950
\font=cmsl9 scaled 1050
\font\pnum=cmbx10 scaled 913
\font\lnum=cmbx10 
\font\pfoot=cmsl9 scaled 950
\font=cmsl9 scaled 1050
\ifplaintex
\headline{\vbox to 0pt{\vskip -4.5mm\line{\small\phead\ifnum
\count0=\startpage ISSN 1364-0380 (on line)
1465-3060 (printed) \hfill {\pnum\folio}\else\ifodd\count0\def\\{ }%
\ifx\theshorttitle\relax\thetitle\else\theshorttitle\fi\hfill{\pnum\folio}
\else\def\\{ and }{\pnum\folio}\hfill\ifx\theshortauthors\relax\theauthors
\else\theshortauthors\fi\fi\fi}\vss}}
\footline{\vbox to 0pt{\vglue 0mm\line{\small\pfoot\ifnum\count0=\startpage
\copyright\ \gtp\hfill\else
\gt, Volume \thevolumenumber\ (\thevolumeyear)\hfill\fi}\vss
}}
\else
\makeatletter
\def\@oddhead{{\small\lhead\ifnum\count0=\startpage ISSN 1364-0380 (on line)
1465-3060 (printed) \hfill {\lnum\number\count0}\else\ifodd\count0
\def\\{ }\ifx\theshorttitle\relax \thetitle \else\theshorttitle\fi\hfill
{\lnum\number\count0}\else\def\\{ and }{\lnum\number\count0}
\hfill\ifx\theshortauthors\relax 
\theauthors\else\theshortauthors\fi\fi\fi}}\def\@evenhead{\@oddhead}
\def\@oddfoot{\small\lfoot\ifnum\count0=\startpage\copyright\ \gtp\hfill\else
\gt, Volume \thevolumenumber\ (\thevolumeyear)\hfill\fi}
\def\@evenfoot{\@oddfoot}
\makeatother
\fi

\newwrite\gtoutfile
\long\gdef\makeheadfile{  
{\def\\{, }\def\s{ }
\immediate\openout\gtoutfile head.xxx
\immediate\write\gtoutfile{To: math@arxiv.org}
\immediate\write\gtoutfile{Subject: put or rep NNNNN:pppp}
\immediate\write\gtoutfile{--text follows this line--}
\immediate\write\gtoutfile{Proxy-for: \ifx\theasciiauthors\relax
\theauthors\else\theasciiauthors\fi\s<\ifx\theasciiemail\relax\theemail\else\theasciiemail\fi>}
\immediate\write\gtoutfile{\noexpand\\}
\immediate\write\gtoutfile{Authors: \ifx\theasciiauthors\relax
\theauthors\else\theasciiauthors\fi}
\immediate\write\gtoutfile{Title: \ifx\theasciititle\relax
\thetitle\else\theasciititle\fi}
\immediate\write\gtoutfile{Subj-class: GT or SG or MG etc}
\immediate\write\gtoutfile{MSC-class: \theprimaryclass\ifx\thesecondaryclass\relax\else, \thesecondaryclass\fi}
\immediate\write\gtoutfile{Journal-ref: Geom. Topol. \thevolumenumber
(\thevolumeyear) \startpage-\finishpage}
\immediate\write\gtoutfile{Comments: Published by Geometry and Topology at}
\immediate\write\gtoutfile{\s\s http://www.maths.warwick.ac.uk/gt/GTVol\thevolumenumber/paper\thepapernumber.abs.html}
\immediate\write\gtoutfile{\noexpand\\}
\immediate\write\gtoutfile{}
\ifx\theasciiabstract\relax
\immediate\write\gtoutfile{\theabstract}\else
\immediate\write\gtoutfile{\theasciiabstract}\fi
\immediate\write\gtoutfile{}
\immediate\write\gtoutfile{\noexpand\\}
\immediate\write\gtoutfile{}
\immediate\closeout\gtoutfile}}  

\def\maketitlepage{\maketitlep\makeheadfile}
\let\maketitle\maketitlepage

\lognumber{168}
\volumenumber{5}\papernumber{20}\volumeyear{2001}
\pagenumbers{609}{650}
\received{23 February 2001}
\accepted{5 June 2001}
\proposed{Cameron Gordon}
\seconded{Joan Birman, Robion Kirby}
\published{27 June 2001}

\usepackage{amssymb,amsmath,graphicx} 
\nocolon
\newtheorem{theorem}{Theorem}[section]
\newtheorem{lemma}[theorem]{Lemma}
\newtheorem{corollary}[theorem]{Corollary}
\newtheorem{proposition}[theorem]{Proposition}

\let\oldcaption\caption
\def\caption#1{\vspace{-20pt}\oldcaption{#1}}

\theoremstyle{definition}
\newtheorem{definition}[theorem]{Definition}

\newtheorem{remark}[theorem]{Remark}

\begin{document}
\title{Heegaard splittings of exteriors of\\two bridge knots }                    
\author{Tsuyoshi Kobayashi }                  
\address{ Department of Mathematics, Nara
Women's University\\Kita-Uoya Nishimachi, Nara 630-8506, JAPAN }               \email{tsuyoshi@cc.nara-wu.ac.jp }                     
\begin{abstract} 
In this paper, we show that, 
for each non-trivial two bridge knot $K$ and for each $g \ge 3$, 
every genus $g$ Heegaard splitting of the exterior $E(K)$ of 
$K$ is reducible. 
\end{abstract}
\asciiabstract{
In this paper, we show that, 
for each non-trivial two bridge knot K and for each g > 2, 
every genus g Heegaard splitting of the exterior E(K) of 
K is reducible.}


\maketitlepage

\section{Introduction}
In this paper, we prove the following theorem. 

\begin{theorem}\label{main theorem}
Let $K$ be a non-trivial two bridge knot. 
Then, for each $g \ge 3$, 
every genus $g$ Heegaard splitting of the exterior $E(K)$ of $K$ is reducible. 
\end{theorem}

We note that since $E(K)$ is irreducible, the above theorem together 
with the classification of the Heegaard splittings of the 3--sphere $S^3$ 
(F Waldhausen \cite{W}) implies the next corollary. 

\begin{corollary}\label{corollary 1}
Let $K$ be a non-trivial two bridge knot. 
Then, for each $g \ge 3$, 
every genus $g$ Heegaard splitting of $E(K)$ is stabilized. 
\end{corollary}

By H Goda, M Scharlemann, and A Thompson \cite{GST} 
(see also K Morimoto's paper \cite{Mo}) or \cite{K}, 
it is shown that, for each non-trivial two bridge knot $K$, 
every genus two Heegaard splitting of $E(K)$ 
is isotopic to either one of six typical 
Heegaard splittings (see Figure~\ref{7.1}). 
We note that Y Hagiwara \cite{Hag} proved that 
genus three Heegaard splittings obtained by stabilizing 
the six Heegaard splittings are mutually isotopic. 
This result together with Corollary~\ref{corollary 1} implies the 
following. 

\begin{corollary}\label{corollary 2}
Let $K$ be a non-trivial two bridge knot. 
Then, for each $g \ge 3$, the genus $g$ Heegaard splittings of 
$E(K)$ are mutually isotopic, 
ie, there is exactly one isotopy class of Heegaard splittings 
of genus $g$. 
\end{corollary}

We note that this result is proved for figure eight knot by 
D Heath \cite{Hea}. 

The author would like to express his thanks to Dr Kanji Morimoto 
for careful readings of a manuscript of this paper. 

\section{Preliminaries}
Throughout this paper, 
we work in the differentiable category. 
For a submanifold $H$ of a manifold $M$, 
$N(H,M)$ denotes a regular neighborhood of $H$ in $M$. 
When $M$ is well understood, we often abbreviate $N(H,M)$ to $N(H)$. 
Let $N$ be a manifold embedded in a manifold $M$
with dim$N=$dim$M$.
Then $\text{Fr}_M N$ denotes the frontier of $N$ in $M$.
For the definitions of standard terms in 
3--dimensional topology, 
we refer to \cite{He} or \cite{Ja}. 

\subsection*{2.A\qua Heegaard splittings}

A 3--manifold $C$ is a {\it compression body} if there exists a compact, 
connected (not necessarily closed) surface $F$ 
such that $C$ is obtained 
from $F \times [0,1]$ by attaching 2--handles along mutually disjoint 
simple closed curves in $F \times \{ 1 \}$ and capping off the resulting 
2--sphere boundary components which are disjoint from $F \times \{ 0 \}$ 
by 3--handles. 
The subsurface of $\partial C$ corresponding to $F \times \{ 0 \}$ 
is denoted by $\partial_+C$. 
Then $\partial_-C$ denotes the subsurface 
$\text{cl} (\partial C -(\partial F \times [0,1] \cup \partial_+C))$ 
of $\partial C$. 
A compression body $C$ is said to be {\it trivial} if either 
$C$ is a 3--ball with $\partial_+ C = \partial C$, or 
$C$ is homeomorphic to $F \times [0,1]$ with $\partial_-C$ corresponding 
to $F \times \{ 0 \}$. 
A compression body $C$ is called a {\it handlebody} if 
$\partial_-C = \emptyset$. 
A compressing disk $D (\subset C)$ of $\partial_+ C$ is called 
a {\it meridian disk} of the compression body $C$. 

\begin{remark}\label{remark of compression body}
The following properties are known for compression bodies. 

\begin{enumerate}

\item 
Compression bodies are irreducible. 

\item 
By extending the cores of the 2--handles in the definition of 
the compression body $C$ vertically to $F \times [0,1]$, 
we obtain a union of mutually disjoint meridian disks 
${\cal D}$ of $C$ such that the manifold 
obtained from $C$ by cutting along ${\cal D}$ is homeomorphic 
to a union of $\partial_- C \times [0,1]$ and some 
(possibly empty) 3--balls. 
This gives a dual description of compression bodies. 
That is, a connected 3--manifold $C$ is a compression body 
if there exists a compact (not necessarily connected) surface 
${\cal F}$ without 2--sphere components and a union of 
(possibly empty) 3--balls ${\cal B}$ such that $C$ is obtained from 
${\cal F} \times [0,1] \cup {\cal B}$ by attaching 1--handles to 
${\cal F} \times \{ 0 \} \cup \partial {\cal B}$. 
We note that $\partial_-C$ is the surface corresponding to 
${\cal F} \times \{ 1 \}$. 

\item 
Let ${\cal D}$ be a union of mutually disjoint meridian disks of a 
compression body $C$, and $C'$ a component of the manifold obtained 
from $C$ by cutting along ${\cal D}$. 
Then, by using the above fact~2, we can show that 
$C'$ inherits a compression body structure from $C$, 
ie, 
$C'$ is a compression body such that 
$\partial_- C' = \partial_- C \cap C'$ 
and 
$\partial_+ C' = (\partial_+ C \cap C') \cup \text{Fr}_C C'$. 

\item 
Let $S$ be an incompressible surface in $C$ such that 
$\partial S \subset \partial_+ C$. 
If $S$ is not a meridian disk, 
then, by using the above fact~2, we can show that 
$S$ is $\partial$--compressible into $\partial_+ C$, 
ie, there exists a disk $\Delta$ such that 
$\Delta \cap S = \partial \Delta \cap S = a$ is an essential arc
in $S$, and 
$\Delta \cap \partial C = \text{cl}(\partial \Delta - a)$ 
with $\Delta \cap \partial C \subset \partial_+ C$. 

\end{enumerate}
\end{remark}

Let $N$ be a cobordism rel~$\partial$ between two surfaces $F_1$, $F_2$ 
(possibly $F_1 = \emptyset$ or $F_2 = \emptyset$), 
ie, 
$F_1$ and $F_2$ are mutually disjoint surfaces in $\partial N$ 
with $\partial F_1 \cong \partial F_2$ such that 
$\partial N = F_1 \cup F_2 \cup (\partial F_1 \times [0, 1])$. 

\begin{definition}\label{Heegaard splitting}
We say that $C_1 \cup_P C_2$ (or $C_1 \cup C_2)$ 
is a {\it Heegaard splitting} 
of $(N, F_1, F_2)$ (or simply, $N$) if it satisfies the following 
conditions. 

\begin{enumerate}
\item 
$C_i$ $(i=1,2)$ is a compression body in $N$ such that 
$\partial_- C_i = F_i$, 

\item 
$C_1 \cup C_2 = N$, and 

\item 
$C_1 \cap C_2 = \partial_+ C_1 = \partial_+ C_2 = P$.

\end{enumerate}

The surface $P$ is called a {\it Heegaard surface} of 
$(N, F_1, F_2)$ (or, $N$). 
In particular, if $P$ is a closed surface, then the genus of 
$P$ is called the {\it genus} of the Heegaard splitting. 
\end{definition}

\begin{definition}\label{reducible Heegaard splitting}

\ 
\begin{enumerate}

\item 
A Heegaard splitting $C_1 \cup_P C_2$ is {\it reducible} 
if there exist meridian disks $D_1$, $D_2$ of 
the compression bodies $C_1$, $C_2$ 
respectively such that 
$\partial D_1= \partial D_2$ 

\item 
A Heegaard splitting $C_1 \cup_P C_2$ is {\it weakly reducible} 
if there exist meridian disks $D_1$, $D_2$ of 
the compression bodies $C_1$, $C_2$ 
respectively such that 
$\partial D_1 \cap  \partial D_2 = \emptyset$. 
If $C_1 \cup_P C_2$ is not weakly reducible, then it is called 
{\it strongly irreducible}. 

\item 
A Heegaard splitting $C_1 \cup_P C_2$ is {\it stabilized} 
if there exists another Heegaard splitting 
$C_1' \cup_{P'} C_2'$ such that the pair 
$(N, P)$ is isotopic to a connected sum of pairs 
$(N, P') \sharp (S^3, T)$, 
where $T$ is a genus one Heegaard surface of the 3--sphere $S^3$. 

\item 
A Heegaard splitting $C_1 \cup_P C_2$ is {\it trivial} 
if either $C_1$ or $C_2$ is a trivial compression body. 
\end{enumerate}
\end{definition}

\begin{remark}

\ 
\begin{enumerate}

\item 
We note that $C_1 \cup_P C_2$ is stabilized if and only if 
there exist meridian disks $D_1$, $D_2$ of $C_1$, $C_2$ 
respectively such that 
$\partial D_1$ and $\partial D_2$ intersect transversely in 
one point. 

\item 
If $C_1 \cup_P C_2$ is stabilized and not a genus one Heegaard splitting 
of $S^3$, then $C_1 \cup_P C_2$ is reducible. 
\end{enumerate}

\end{remark}

\subsection*{2.B\qua Orbifold version of Heegaard splittings}

Throughout this subsection, let $N$ be a compact, orientable 3--manifold, 
$\gamma$ a 1--manifold properly embedded in $N$, and 
$F$, $F_1$, $F_2$, $D$, $S$, connected surfaces embedded in $N$, 
which are in general position 
with respect to $\gamma$.

\begin{definition}
We say that $D$ is a {\it $\gamma$--disk}, 
if (1) $D$ is a disk, and 
(2) either $D \cap \gamma = \emptyset$, or 
$D$ intersects $\gamma$ transversely in one point. 
\end{definition}

Let $\ell (\subset F)$ be a simple closed curve such that 
$\ell \cap \gamma = \emptyset$.

\begin{definition}
We say that $\ell$ is $\gamma$--{\it inessential} 
if $\ell$ bounds a $\gamma$--disk in $F$. 
We say that $\ell$ is $\gamma$--{\it essential} if 
it is not $\gamma$--inessential.
\end{definition}

\begin{figure}[ht!]
\begin{center}
\includegraphics[width=6cm, clip]{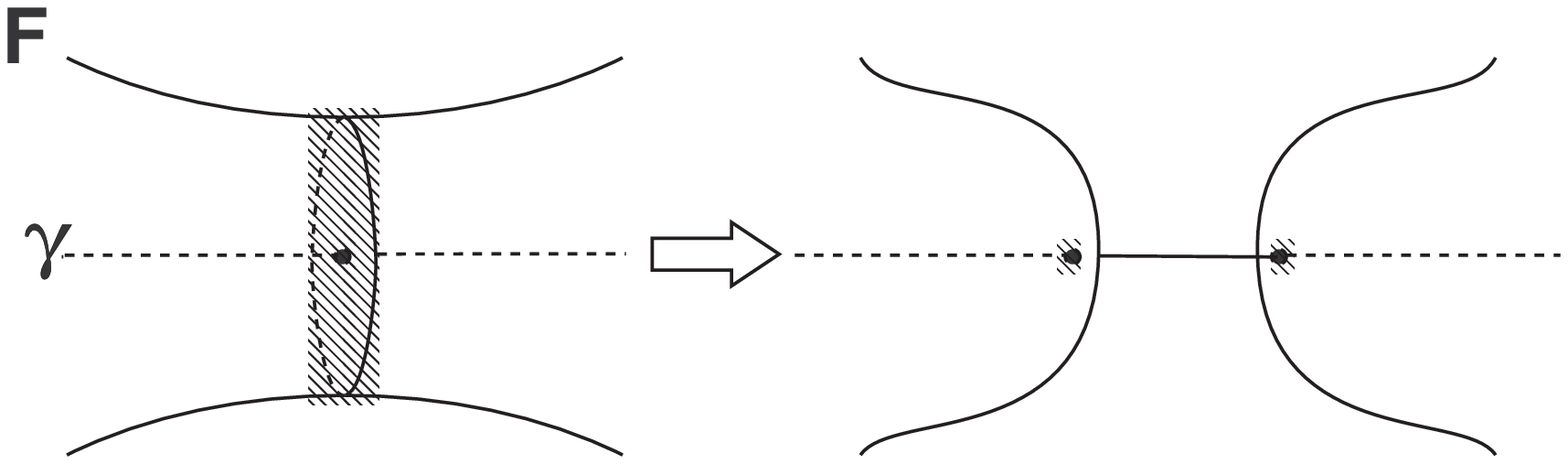}
\end{center}\vspace{10pt}
\caption{}\label{2.2.1}
\end{figure}

\begin{definition}
We say that $D$ is a $\gamma$--{\it compressing disk} 
for $F$ if 
$D$ is a $\gamma$--disk, $D \cap F = \partial D$, and 
$\partial D$ is a $\gamma$--essential simple closed curve in $F$. 
The surface $F$ is $\gamma$--{\it compressible} 
if it admits a $\gamma$--compressing disk, and 
$F$ is $\gamma$--{\it incompressible} if it is not 
$\gamma$--compressible. 
We note that if $D$ is a $\gamma$--compressing disk for $F$, then we 
can perform a $\gamma$--compression on $F$ along $D$ (Figure~\ref{2.2.1}). 
\end{definition}

\begin{definition}
Suppose that $\partial F_1 = \partial F_2$, 
or 
$\partial F_1 \cap \partial F_2 = \emptyset$. 
We say that $F_1$ and $F_2$ are $\gamma$--{\it parallel}, 
if there is a submanifold $R$ in $N$ such that 
$(R, R \cap \gamma )$ is homeomorphic to 
$(F_1 \times [0,1], {\cal P} \times [0,1])$ 
as a pair, where 
(1) ${\cal P}$ is a union of points in
$\text{Int} F_1$, and 
(2) $\partial F_1 = \partial F_2$ and 
$F_1$ ($F_2$ respectively) is the subsurface of $\partial R$ 
corresponding to the closure of the component of 
$\partial (F_1 \times [0,1]) - (\partial F_1 \times \{ 1/2 \})$ 
containing $F_1 \times \{ 0 \}$ ($F_1 \times \{ 1 \}$ respectively), or 
$\partial F_1 \cap \partial F_2 = \emptyset$ and 
$F_1$ ($F_2$ respectively) is the subsurface of $\partial R$ corresponding to 
$F_1 \times \{ 0 \}$ ($F_1 \times \{ 1 \}$ respectively). 
The submanifold $R$ is called a {\it $\gamma$--parallelism} 
between $F_1$ and $F_2$. 

We say that $F$ is $\gamma$--{\it boundary parallel} 
if there is a subsurface $F'$ in $\partial N$  such that 
$F$ and $F'$ are $\gamma$--parallel. 
\end{definition}

\begin{definition}
We say that $S$ is a {\it $\gamma$--sphere} 
if (1) $S$ is a sphere, and 
(2) either $S \cap \gamma = \emptyset$, or 
$S$ intersects $\gamma$ transversely in two points. 
We say that a 3--ball $B^3$ in $N$ is a {\it $\gamma$--ball} 
if either $B^3 \cap \gamma = \emptyset$, or 
$B^3 \cap \gamma$ is an unknotted arc properly embedded in $B^3$. 
A $\gamma$--sphere $S$ is {\it $\gamma$--compressible} if there 
exists a $\gamma$--ball $B^3$ in $N$ such that 
$\partial B^3 = S$. 
A $\gamma$--sphere $S$ is {\it $\gamma$--incompressible} if it is not 
$\gamma$--compressible. 
We say that $N$ is {\it $\gamma$--reducible} if $N$ contains a 
$\gamma$--incompressible 2--sphere. 
The manifold $N$ is {\it $\gamma$--irreducible} if it is not 
$\gamma$--reducible. 
\end{definition}

\begin{definition}
We say that $F$ is {\it $\gamma$--essential} 
if $F$ is $\gamma$--incompressible, and 
not $\gamma$--boundary parallel. 
\end{definition}

Let $a$ be an arc properly embedded in $F$ 
with $a \cap \gamma = \emptyset$.

\begin{definition}
We say that $a$ is {\it $\gamma$--inessential} 
if there is a subarc $b$ of $\partial F$ 
such that $\partial b = \partial a$, 
and $a \cup b$ bounds a disk $D$ in $F$ such that 
$D \cap \gamma = \emptyset$, 
and $a$ is {\it $\gamma$--essential} if it is not 
$\gamma$--inessential. 
\end{definition}

\begin{definition}
We say that $\Delta$ is a {\it $\gamma$--boundary compressing disk} 
for $F$ if $\Delta$ is a disk disjoint from $\gamma$, 
$\Delta \cap F = \partial \Delta \cap F = \alpha$ 
is a $\gamma$--essential arc in $F$, and 
$\Delta \cap \partial N = \partial \Delta \cap \partial N = 
\text{cl} (\partial \Delta - \alpha)$. 
The surface $F$ is {\it $\gamma$--boundary compressible} 
if it admits a $\gamma$--boundary compressing disk. 
The surface $F$ is {\it $\gamma$--boundary incompressible} 
if it is not $\gamma$--boundary compressible. 
We note that if $\Delta$ is a $\gamma$--boundary compressing disk for $F$, 
then we can perform a {\it $\gamma$--boundary compression} on $F$ along 
$\Delta$. 

\end{definition}

\begin{definition}
We say that $F_1$ and $F_2$ are {\it $\gamma$--isotopic} 
if there is an ambient isotopy 
$\phi_t$ $(0 \le t \le 1)$ of $N$ such that 
$\phi_0 = id_N$, 
$\phi_1(F_1) = F_2$, and 
$\phi_t(\gamma ) = \gamma$ for each $t$. 
\end{definition}

The next definition gives an orbifold version of compression body 
(cf (2) of Remark~\ref{remark of compression body}). 

\begin{definition}\label{o-compression body}
Suppose that $N$ is a cobordism rel~$\partial$ between 
two surfaces $G_+$, $G_-$. 
We say that $(N, \gamma)$ is an 
{\it orbifold compression body} 
(or {\it O--compression body}) 
(with $\partial _+ N = G_+$, and $\partial_-N = G_-$) 
if the following conditions are satisfied. 

\begin{enumerate}

\item 
$G_+$ is not empty, and is connected
(possibly, $G_-$ is empty). 

\item 
No component of $G_-$ is a $\gamma$--sphere. 

\item 
$\partial \gamma \subset \text{Int}(G_+ \cup G_-)$. 

\item 
There exists a union of mutually disjoint 
$\gamma$--compressing disks, say ${\cal D}$, for $G_+$ such that, 
for each component $E$ of the manifold obtained from $N$ 
by cutting along ${\cal D}$, 
either 
$E$ is a $\gamma$--ball with $E \cap G_- = \emptyset$, or 
$(E, \gamma \cap E)$ is homeomorphic to 
$( G \times [0,1], {\cal P} \times [0,1])$, 
where $G$ is a component of $G_-$ with 
$E \cap G_- = G \times \{ 0 \} = G$ 
and ${\cal P}$ is a union of mutually disjoint 
(possibly empty) points in $G$ (see Figure~\ref{2.2.2}). 

\end{enumerate}

\begin{figure}[ht!]
\begin{center}
\includegraphics[width=7cm, clip]{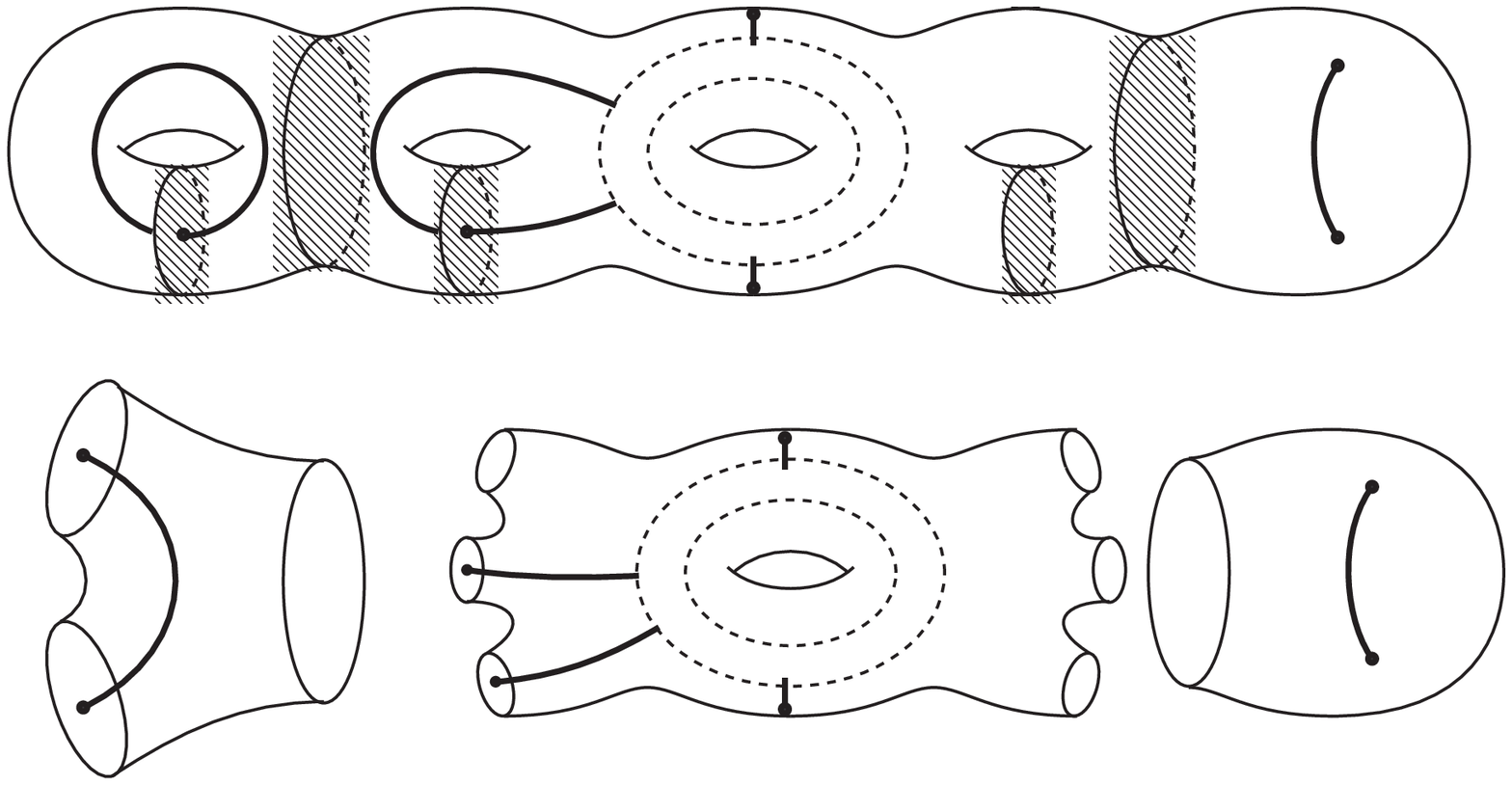}
\end{center}
\caption{}\label{2.2.2}
\end{figure}

Note that the condition~1 of Definition~\ref{o-compression body} 
implies that $N$ is connected. 
We say that an O--compression body $(N, \gamma )$ is 
{\it trivial} if either $N$ is a $\gamma$--ball with 
$\partial_+ N = \partial N$, or 
$(N, \gamma )$ is homeomorphic to 
$( G_- \times [0,1], {\cal P}' \times [0,1])$
with $G_-$ $(\subset \partial N)$ corresponding to 
$G_- \times \{ 1 \}$, and ${\cal P}'$ a union of 
mutually disjoint points in $G_-$. 
An O--compression body $(N, \gamma )$ is called an {\it O--handlebody} 
if $\partial_-N = \emptyset$. 
A $\gamma$--compressing disk of $\partial_+ N$ is called 
a {\it ($\gamma$--)meridian disk} 
of the $O$--compression body $(N, \gamma )$. 
\end{definition}

By ${\mathbb Z}_2$--equivariant loop theorem \cite[Lemma 3]{KT}, and 
${\mathbb Z}_2$--equivariant cut and paste argument as in 
\cite[Proof of 10.3]{He}, we can prove the following 
(the proof is omitted). 

\begin{proposition}\label{2-fold cover, incompressible}
Let $N$ be a compact, orientable 3--manifold, and 
$\gamma$ a 1--manifold properly embedded in $N$. 
Suppose that $N$ admits a 2--fold branched cover 
$p\co  \tilde{N} \rightarrow N$ with branch set $\gamma$. 
Let $F$ be a (possibly closed) surface properly embedded 
in $N$, which is in general position with respect to $\gamma$. 
Then $F$ is $\gamma$--incompressible 
($\gamma$--boundary incompressible respectively) 
if and only if 
$p^{-1}(F)$ is incompressible 
(boundary incompressible respectively) 
in $\tilde{N}$. 
\end{proposition}

By (2) of Remark~\ref{remark of compression body}, 
Definition~\ref{o-compression body}, 
${\mathbb Z}_2$--equivariant cut and paste argument as in 
\cite[Proof of 10.3]{He}, 
and ${\mathbb Z}_2$--Smith conjecture \cite{W}, 
we immediately have the following. 

\begin{proposition}\label{2-fold cover of compression body}
Let $N$, $\gamma$ be as in Proposition~\ref{2-fold cover, incompressible}. 
Then $(N, \gamma )$ is an $O$--compression body with 
$\partial_{\pm} N = G_{\pm}$, 
if and only if $\tilde{N}$ is a compression body with 
$\partial_{\pm} \tilde{N} = p^{-1}(G_{\pm})$. 
\end{proposition}

Since the compression bodies are irreducible 
(see (1) of Remark~\ref{remark of compression body}), 
Proposition~\ref{2-fold cover of compression body}
together with ${\mathbb Z}_2$--Smith conjecture \cite{W}
implies the following. 

\begin{corollary}\label{gamma-compression body is gamma-irreducible}
Let $(N, \gamma )$ be an O--compression body. 
Suppose that $N$ admits a 2--fold branched cover 
with branch set $\gamma$. 
Then $N$ is $\gamma$--irreducible. 
\end{corollary}

By (4) of Remark~\ref{remark of compression body}, and 
${\mathbb Z}_2$--equivariant cut and paste argument as in 
\cite[Proof of 10.3]{He}, we have the following. 

\begin{corollary}\label{incomp., boundary incomp. surface in gamma-comp.}
\label{gamma-incomp., boundary incomp. surface in gamma-comp. body}
Let $(N, \gamma )$ be an O--compression body such that 
$N$ admits a 2--fold branched cover with branch set $\gamma$. 
Let $F$ be a connected $\gamma$--incompress\-ible surface 
properly embedded in $N$, which is not a $\gamma$--meridian disk. 
Suppose that 
$\partial F \subset \partial_+ N$. 
Then there exists a $\gamma$--boundary compressing disk $\Delta$ 
for $F$ such that $\Delta \cap \partial N \subset \partial_+ N$. 

\end{corollary}

Let $M$ be a compact, orientable 3--manifold, and 
$\delta$ a 1--manifold properly embedded in $M$. 
Let $C$ be a 3--dimensional manifold embedded in $M$. 
We say that $C$ is a {\it $\delta$--compression body} 
if $(C, \delta \cap C)$ is an O--compression body. 
Suppose that $M$ is a cobordism rel~$\partial$ between two surfaces 
$G_1$, $G_2$ 
(possibly $G_1 = \emptyset$ or $G_2 = \emptyset$) 
such that 
$\partial \delta \subset \text{Int} (G_1 \cup G_2)$. 

\begin{definition}\label{orbifold version of Heegaard splitting}
We say that $C_1 \cup_P C_2$ is a {\it Heegaard splitting} 
of $(M, \delta , G_1, G_2)$ (or simply $( M, \delta )$) if 
it satisfies the following conditions. 

\begin{enumerate}
\item 
$C_i$ $(i=1,2)$ is a $\delta$--compression body such that 
$\partial_- C_i = G_i$, 

\item 
$C_1 \cup C_2 = M$, and 

\item 
$C_1 \cap C_2 = \partial_+ C_1 = \partial _+ C_2 =P$. 
\end{enumerate}

The surface $P$ is called a {\it Heegaard surface} of 
$(M, \delta , G_1, G_2)$ (or $( M, \delta )$). 
\end{definition}

\begin{definition}\label{reducible o-Heegaard splitting}

\ 
\begin{enumerate}
\item 
A Heegaard splitting $C_1 \cup_P C_2$ of $(M, \delta )$ is 
{\it $\delta$--reducible} if there exist 
$\delta$--meridian disks $D_1$, $D_2$ of 
the $\delta$--compression bodies $C_1$, $C_2$ respectively such that 
$\partial D_1 = \partial D_2$. 

\item 
A Heegaard splitting $C_1 \cup_P C_2$ 
of $( M, \delta )$ is {\it weakly $\delta$--reducible} 
if there exist $\delta$--meridian disks $D_1$, $D_2$ of 
the $\delta$--compression bodies $C_1$, $C_2$ 
respectively such that 
$\partial D_1 \cap  \partial D_2 = \emptyset$. 
If $C_1 \cup_P C_2$ is not weakly $\delta$--reducible, then it is called 
{\it strongly $\delta$--irreducible}. 

\item 
A Heegaard splitting $C_1 \cup_P C_2$ 
of $( M, \delta )$ is {\it trivial} 
if either $C_1$ or $C_2$ is a trivial $\delta$--compression body. 
\end{enumerate}
\end{definition}

\subsection*{2.C\qua Genus $g$, $n$--bridge positions}

We first recall the definition of a genus $g$, $n$--bridge position 
of H.Doll \cite{D}. 
Let $\Gamma = \gamma_1 \cup \cdots \cup \gamma_n$ be a union 
of mutually disjoint arcs $\gamma_i$ properly embedded in a 3--manifold $N$. 

\begin{definition}\label{trivial arcs}
We say that $\Gamma$ is {\it trivial} if there exist 
mutually disjoint disks 
$D_1, \dots , D_n$ in $N$ 
such that 
(1) $D_i \cap \Gamma = \partial D_i \cap \gamma_i = \gamma_i$, and 
(2) $D_i \cap \partial N = \text{cl}(\partial D_i - \gamma_i)$. 
\end{definition}


Let $K$ be a link in a closed 3--manifold $M$. 
Let $X \cup_Q Y$ be a genus $g$ Heegaard splitting of $M$. 
Then, with following \cite{D}, we say that 
$K$ is in a {\it genus $g$, $n$--bridge position} 
(with respect to the Heegaard splitting $X \cup_Q Y$) if 
$K \cap X$ ($K \cap Y$ respectively) is a union of 
$n$ arcs which is trivial in $X$ ($Y$ respectively). 

A proof of the next lemma is elementary, and we omit it. 

\begin{lemma}
Let $\Gamma$ be a union 
of mutually disjoint arcs properly embedded in a handlebody $H$. 
Then $\Gamma$ is trivial if and only if $(H, \Gamma )$ is a 
O--handlebody. 
\end{lemma}

This lemma allows us to generalize 
the definition of genus $g$, $n$--bridge positions 
as in the following form. 
Let $K$, $M$, and $X \cup_Q Y$ be as above. 

\begin{definition}\label{def. of genus g, n-bridge pos.}
We say that $K$ is in a {\it genus $g$, $n$--bridge position} 
(with respect to the Heegaard splitting $X \cup_Q Y$) if 
$X \cup_Q Y$ gives a Heegaard splitting of $(M,K)$ such that 
$\text{genus}(Q)=g$, and $K \cap Q$ consists of $2n$ points. 
\end{definition}

\begin{remark}
This definition allows genus $g$, $0$--bridge position of $K$. 
\end{remark}

In this paper, we abbreviate genus 0, $n$--bridge 
position to {\it $n$--bridge position}.

\begin{definition}\label{n-bridge knot}
A knot $K$ in the 3--sphere $S^3$ is called a {\it $n$--bridge knot}, 
if it admits a $n$--bridge position. 
\end{definition}

\section{Weakly $\gamma$--reducible Heegaard splittings} 

In \cite{Ha}, W Haken proved that the Heegaard splittings of a 
reducible 3--manifold are reducible. 
As a sequel of this, Casson--Gordon \cite{CG} proved that 
each non-trivial Heegaard splitting of a $\partial$--reducible 
3--manifold is weakly reducible. 
In this section, we prove orbifold versions of these results. 
In fact, we prove the following. 

\begin{proposition}\label{Haken's theorem}
Let $N$ be a compact orientable 3--manifold, and 
$\gamma$ a 1--manifold properly embedded in $N$ 
such that 
$N$ admits a 2--fold branched cover with branch set $\gamma$. 
Suppose that $N$ is a cobordism rel~$\partial$ between two surfaces 
$F_1$, $F_2$ 
(possibly $F_1 = \emptyset$ or $F_2 = \emptyset$) 
such that 
$\partial \gamma \subset \text{Int} (F_1 \cup F_2)$, and 
no component of $F_1 \cup F_2$ is a $\gamma$--disk. 
If $N$ is $\gamma$--reducible, 
then every Heegaard splitting of $(N, \gamma , F_1, F_2)$ is 
weakly $\gamma$--reducible. 
\end{proposition}

\begin{proposition}\label{Casson-Gordon's lemma}
Let $N$, $\gamma$, $F_1$, $F_2$ be as in Proposition~\ref{Haken's theorem}. 
If $F_1 \cup F_2$ is $\gamma$--compressible in $N$, 
then every non-trivial Heegaard splitting of 
$(N, \gamma , F_1, F_2 )$ is weakly $\gamma$--reducible. 
\end{proposition}

\begin{remark}
In the conclusion of Proposition~\ref{Haken's theorem}, 
we can have just \lq\lq weakly $\gamma$--reducible \rq\rq, 
not \lq\lq $\gamma$--reducible \rq\rq. 
For example, let $K$ be a connected sum of two trefoil knots, 
and $C_1 \cup C_2$ the Heegaard splitting of $(S^3, K)$ 
as in Figure~\ref{3.1}. 
We note that $(S^3, K)$ is $K$--reducible 
(in fact, a 2--sphere giving prime decomposition of $K$ is 
$K$--incompressible). 
Since the Heegaard splitting gives a minimal genus 
Heegaard splitting of $E(K)$, we can show that 
$C_1 \cup C_2$ is not $\gamma$--reducible. 
But $C_1 \cup C_2$ admits a pair of weakly $K$--reducing 
disks $D_1$, $D_2$ as in Figure~\ref{3.1}. 

\begin{figure}[ht!]
\begin{center}
\includegraphics[width=5cm, clip]{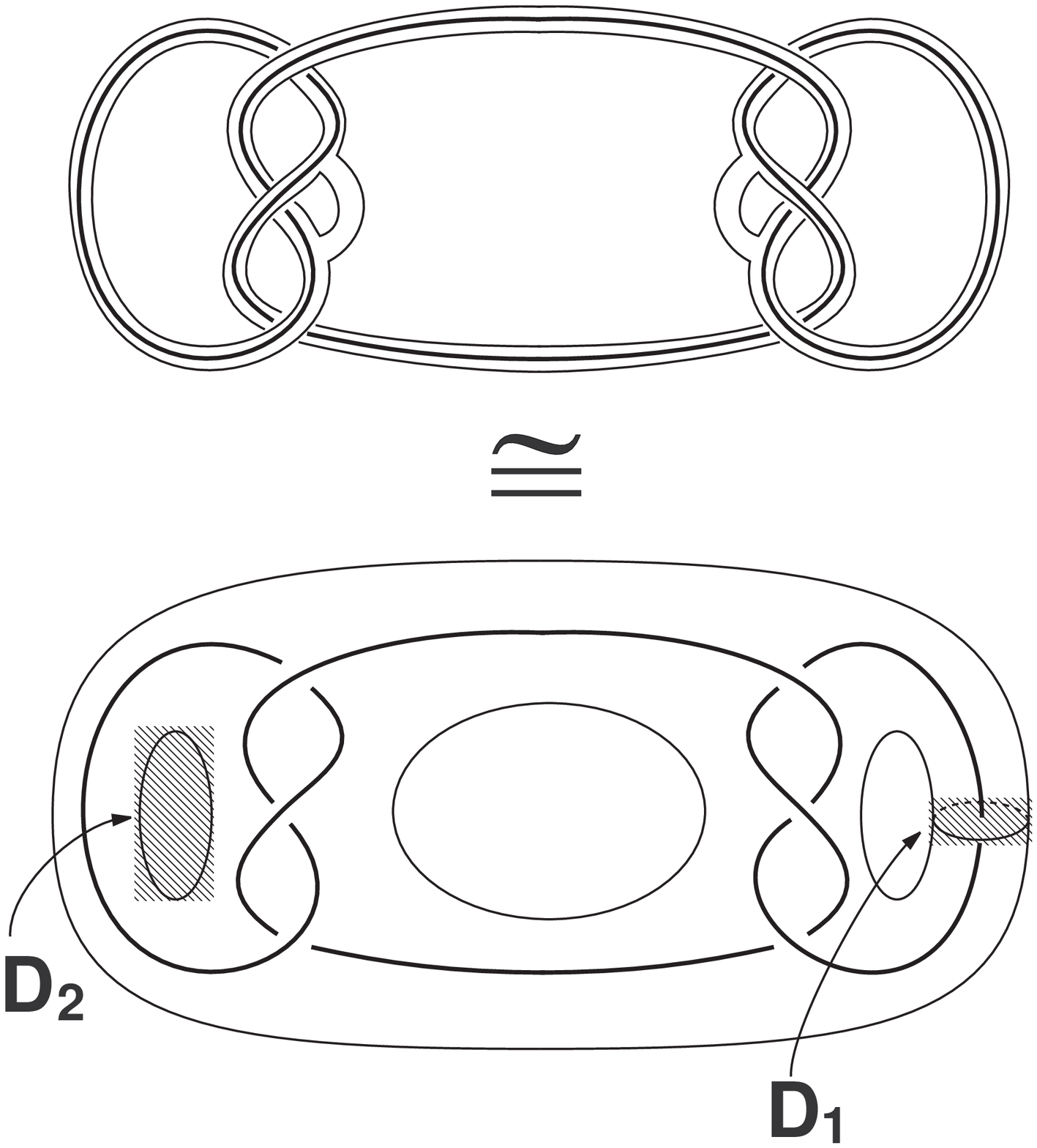}
\end{center}
\caption{}\label{3.1}
\end{figure}
\end{remark}

Then, by using Proposition~\ref{Casson-Gordon's lemma}, 
we prove an orbifold version of a lemma of Rubinstein--Scharlemann
\cite[Lemma 4.5]{RS}. 

\begin{proposition}\label{Rubinstein-Scharlemann's lemma}
Let $M$ be a closed orientable 3--manifold, and 
$K$ a link in $M$ such that $M$ admits a 2--fold branched cover 
with branch set $K$. 
Let $A \cup_P B$, $X \cup_Q Y$ be Heegaard splittings of $(M,K)$. 
Suppose that $A \subset \text{Int} X$, 
and there exists a $K$--meridian disk $D$ of $X$ such that 
$D \cap A = \emptyset$. 
Then we have one of the following. 

\begin{enumerate}

\item 
$M$ is homeomorphic to the 3--sphere, 
and either $K = \emptyset$ or $K$ is a trivial knot. 

\item 
$X \cup_Q Y$ is weakly $K$--reducible. 
\end{enumerate}
\end{proposition}

\subsection*{3.A\qua Heegaard splittings of 
$(\hat{N}, \hat{F_1}, \hat{F_2})$} 

For the proofs of Propositions~\ref{Haken's theorem}, and 
\ref{Casson-Gordon's lemma}, we show that we can derive  Heegaard splittings 
of $\text{cl} (N - N(\gamma ))$ from Heegaard splittings of 
$(N, \gamma )$. 

\begin{lemma}\label{o-compression body, compression body}
Let $(C, \beta )$ be a O--compression body such that 
$C$ admits a 2--fold branched cover 
$q\co  \tilde{C} \rightarrow C$ with branch set $\beta$. 
Let $\hat{C} = \text{cl}(C-N(\beta ))$, 
$S = \text{cl}(\partial_+ C - N( \beta ))$. 
Then $\hat{C}$ is a compression body with 
$\partial_+ \hat{C} = S$. 
\end{lemma}

\begin{proof}
Let ${\cal D}$ be the union of mutually disjoint 
$\beta$--compressing disks for $\partial_+ C$ as in 
Definition~\ref{o-compression body}. 
Let ${\cal D}_0$ (${\cal D}_1$ respectively) be the union of 
the components of ${\cal D}$ which are disjoint from $\beta$
(which intersect $\beta$ respectively). 
Let $E$ be a component of the manifold obtained from $C$ 
by cutting along ${\cal D}_0$, 
and $\hat{E} = \text{cl}(E - N( \beta ))$. 
Let ${\cal D}_{1,E}$ be the union of the components of ${\cal D}_1$ 
that are contained in $E$. 
Let $E'$ be the manifold obtained from $E$ by cutting along 
${\cal D}_{1,E}$, and $\hat{E'} = \text{cl}(E' - N( \beta ))$. 
Then we have the following cases. 

\medskip
\noindent
{\bf Case 1}\qua  $E \cap \beta = \emptyset$. 

\medskip

In this case, ${\cal D}_{1,E} = \emptyset$, 
and we have 
$E = \hat{E} = E' = \hat{E'}$. 
By the definition of $\beta$--compression body 
(Definition~\ref{o-compression body}), we see that 
$\hat{E} (= E)$ is a trivial compression body 
such that $\hat{E} \cap {\cal D}_0 \subset \partial_+ \hat{E}$. 

\medskip
\noindent
{\bf Case 2}\qua  $E \cap \beta \ne \emptyset$, and 
$E \cap \partial_- C = \emptyset$. 

\medskip
By the definition of $\beta$--compression body, 
we see that each component of $E'$ is a $\beta$--ball 
intersecting $\beta$ with $E' \cap \partial_- C = \emptyset$. 
Hence each component of $\hat{E}'$ is a solid torus, say $T$, 
such that $T \cap N( \beta )$ is an annulus which is a neighborhood 
of a longitude of $T$. 
This implies that 
each component of $\hat{E'}$ is a trivial compression body 
such that the union of the $\partial_+$ boundaries is 
$\text{cl}(\partial \hat{E}' - N( \beta ))$. 
Since $\hat{E}$ is recovered from $\hat{E'}$ by identifying 
pairs of annuli corresponding to 
$\text{cl}({\cal D}_{1,E} - N( \beta ))$, 
we see that the triviality can be pulled back 
to show that $\hat{E}$ is a trivial compression body 
with $\partial_+ \hat{E} = 
\text{cl}(\partial \hat{E} - N( \beta )) = 
\text{cl}(\partial_+ E - N( \beta ))$, 
where $\hat{E} \cap {\cal D}_0 \subset \partial_+ \hat{E}$. 
In fact, we see that either $E$ is a $\beta$--ball or 
$E$ is a solid torus with $\beta \cap E$ a core circle. 

\medskip
\noindent
{\bf Case 3}\qua  $E \cap \beta \ne \emptyset$, and 
$E \cap \partial_- C \ne \emptyset$. 

\medskip
By the definition of $\beta$--compression body, 
for each component $E^*$ of $E'$, we have either 
$E^*$ is a $\beta$--ball 
intersecting $\beta$ with $E^* \cap \partial_- C = \emptyset$, 
or $(E^*, E^* \cap \beta )$ is a trivial $\beta$--compression body 
such that the $\partial_-$ boundary is a component of $\partial_- C$. 
In either case, 
$\hat{E}^* = \text{cl} (E^*-N( \beta ))$ 
is a trivial compression body such that 
$\partial_+ \hat{E}^* = \text{cl} (\partial_+ E^*-N( \beta ))$. 
Hence $\hat{E'}$ is a union of trivial compression bodies 
such that the union of the $\partial_+$ boundaries is 
$\text{cl}(\partial_+ E' - N( \beta ))$. 
Since $\hat{E}$ is recovered from $\hat{E'}$ by identifying 
pairs of annuli corresponding to 
$\text{cl}({\cal D}_{1,E} - N( \beta ))$, 
we see that the triviality can be pulled back 
to show that $\hat{E}$ is a trivial compression body 
with $\partial_+ \hat{E} = \text{cl}(\partial_+ E - N( \beta ))$, 
where $\hat{E} \cap {\cal D}_0 \subset \partial_+ \hat{E}$. 

\medskip
By the conclusions of Cases~1, 2 and 3, we see that 
$\hat{C}$ is recovered from a union of trivial compression bodies 
by identifying the pairs of disks in $\partial_+$ boundaries, 
which are corresponding to ${\cal D}_0$, 
and this implies that $\hat{C}$ is a compression body 
(see (2) of Remark~\ref{remark of compression body}). 
Moreover, since the $\partial_+$ boundary of each trivial compression body
$\hat{E}$ is $\text{cl}(\partial_+ E- N( \beta ))$, 
we see that 
$\partial_+\hat{C} = \text{cl}(\partial_+ C - N( \beta ))$. 
\end{proof}

Let 
$C_1 \cup_P C_2$ be a Heegaard splitting of 
$(N, \gamma , F_1, F_2)$. 
Then let 
$\hat{N} = \text{cl} (N-N( \gamma ))$, 
$\hat{P} = \text{cl} (P-N( \gamma ))$, 
$\hat{C_i} = \text{cl} (C_i-N( \gamma ))$, and 
$\hat{F_i} = \text{cl} (\partial \hat{C_i}-N(\hat{P}, \partial \hat{C_i}))$ 
$(i=1,2)$. 
By Lemma~\ref{o-compression body, compression body}, 
we see that $\hat{C_1} \cup_{\hat{P}} \hat{C_2}$ is a 
Heegaard splitting of $(\hat{N}, \hat{F_1}, \hat{F_2})$. 
By the definitions of strongly irreducible Heegaard splittings, and 
strongly $\gamma$--irreducible Heegaard splittings, 
we immediately have the following. 

\begin{lemma}\label{strongly gamma-irreducible, strongly irreducible}
If $C_1 \cup_P C_2$ is strongly $\gamma$--irreducible, 
then $\hat{C_1} \cup_{\hat{P}} \hat{C_2}$ is strongly irreducible. 
\end{lemma}

\subsection*{3.B\qua Proof of Proposition~\ref{Haken's theorem}}

Let $N$, $\gamma$ be as in Proposition~\ref{Haken's theorem}, 
and 
$C_1 \cup_P C_2$ a Heegaard splitting of $(N, \gamma )$. 
Let $\hat{N} = \text{cl}(N-N( \gamma ))$, and 
$\hat{C_1} \cup_{\hat{P}} \hat{C_2}$ a Heegaard splitting of 
$(\hat{N}, \hat{F}_1, \hat{F}_2)$ 
obtained from $C_1 \cup_P C_2$ as in Section~3.A. 
Since $(N, \gamma )$ is $\gamma$--reducible, 
there exists a $\gamma$--incompressible $\gamma$--sphere $S$ in $N$. 
Then we have the following two cases. 

\medskip
\noindent
{\bf Case 1}\qua $S \cap \gamma = \emptyset$. 

\medskip
In this case, we may regard that $S$ is a 2--sphere in $\hat{N}$. 
It is clear that $S$ is an incompressible 2--sphere in $\hat{N}$. 
Hence, by \cite[Lemma 1.1]{CG}, 
we see that there exists an incompressible 2--sphere $S'$ in $\hat{N}$ 
such that $S'$ intersects $\hat{P}$ in a circle. 
Since $\hat{N} \subset N$, 
we may regard $S'$ is a 2--sphere in $N$. 
It is clear that $S' \cap P$ is a $\gamma$--essential simple closed 
curve in $P$, hence, $S' \cap C_i$ $(i=1,2)$ is a $\gamma$--meridian 
disk of $C_i$. 
This shows that $C_1 \cup_P C_2$ is $\gamma$--reducible. 

\medskip
\noindent
{\bf Case 2}\qua $S \cap \gamma \ne \emptyset$ 
(ie,  $S \cap \gamma$ consists of two points). 

\medskip
We may suppose that $(S \cap \gamma ) \cap P = \emptyset$. 
Let $\hat{S} = \text{cl}(S-N( \gamma ))$. 
Then $\hat{S}$ is an annulus properly 
embedded in $\hat{N}$ such that 
$\partial \hat{S} \subset \text{Fr}_N N( \gamma )$, and 
$\partial \hat{S} \cap \hat{P} = \emptyset$. 

\medskip
\noindent
{\bf Claim 1}\qua $\hat{S}$ is incompressible in $\hat{N}$. 

\begin{proof}
If there is a compressing disk $D$ for $\hat{S}$, 
then by compressing $S$ along $D$, we obtain two 2--spheres, 
each of which intersects $\gamma$ in one point. 
This contradicts the existence of a 2--fold branched cover of $N$ 
with branch set $\gamma$. 
\end{proof}

\medskip
\noindent
{\bf Claim 2}\qua $\hat{S}$ is not $\partial$--parallel in $\hat{N}$. 

\begin{proof}
Suppose that $\hat{S}$ is parallel to an annulus ${\cal A}$ 
in $\partial \hat{N}$. 
Let $s = \text{cl}(\partial N-(F_1 \cup F_2))$. 
Note that $s$ is a (possibly empty) union of annulus. 
Let $F_i' = \text{cl}(F_i - N( \gamma ))$. 
Then $\partial \hat{N} = 
s \cup F_1' \cup F_2' \cup \text{Fr}_N N( \gamma )$. 
Since $S$ is $\gamma$--incompressible, 
we see that 
$(F_1' \cup F_2') \cap {\cal A} \ne \emptyset$. 
Since no component of $F_1 \cup F_2$ is a $\gamma$--disk, 
each component of 
$(F_1' \cup F_2') \cap {\cal A}$ is an annulus. 
Let $A^*$ be a component of $\text{Fr}_N N( \gamma )$ such that 
$A^*$ contains a component of $\partial \hat{S}$. 
Let $F^*$ be the component of $(F_1' \cup F_2') \cap {\cal A}$ 
such that $F^* \cap A^* \ne \emptyset$. 
Note that $F^* \cap A^*$ is a component of $\partial A^*$ and 
is also a component of $\partial F^*$. 
Let $A'$ be the component of 
$\text{cl} (\partial \hat{N} - (F_1' \cup F_2'))$ such that 
$A' \cap F^*$ is the component of $\partial F^*$ other than $F^* \cap A^*$. 
Then $A'$ is an annulus which is 
either a component of $\text{Fr}_N N( \gamma )$, 
or a component of $s$. 
If $A'$ is a component of $\text{Fr}_N N( \gamma )$, then the component of 
$F_1 \cup F_2$ corresponding to $F^*$ is a $\gamma$--sphere, 
hence, a component of $\partial_- C_1$ or $\partial_- C_2$ 
is a $\gamma$--sphere, a contradiction. 
If $A'$ is a component of $s$, then the component of 
$F_1 \cup F_2$ corresponding to $F^*$ is a $\gamma$--disk, 
contradicting the assumption of Proposition~\ref{Haken's theorem}. 
\end{proof}

By Claims~1 and 2, $\hat{S}$ is $\gamma$--essential in $\hat{N}$. 
Suppose, for a contradiction, that 
$C_1 \cup_P C_2$ is strongly $\gamma$--irreducible. 
By Lemma~\ref{strongly gamma-irreducible, strongly irreducible}, 
$\hat{C_1} \cup_{\hat{P}} \hat{C_2}$ is strongly irreducible. 
Then, by \cite[Lemma 6]{Schu} or \cite[Lemma 2.3]{Mo'}, 
$\hat{S}$ is ambient isotopic rel~$\partial$ to a surface $\hat{S}'$ 
such that $\hat{S}' \cap \hat{P}$ consists of essential simple closed 
curves in $\hat{S}'$. 
We regard $\hat{S} = \hat{S}'$. 
This means that each component of $S \cap P$ is a simple closed 
curve which separates the points $S \cap \gamma$. 
We suppose that $\vert S \cap P \vert$ is minimal 
among the $\gamma$--incompressible $\gamma$--spheres with this property. 
Let $n = \vert S \cap P \vert$. 
Suppose that $n=1$, ie,  $S \cap P$ consists of a 
simple closed curve, say $\ell_1$. 
Then $\ell_1$ separates $S$ into two $\gamma$--disks, 
which are $\gamma$--meridian disks in $C_1$ and $C_2$ respectively. 
This shows that $C_1 \cup_P C_2$ is $\gamma$--reducible, a contradiction. 
Suppose that $n \ge 2$. 
Let $D_1$ be the closure of a component of $S-P$ such that 
$D_1 \cap \gamma \ne \emptyset$. 
Note that $D_1$ is a $\gamma$--disk. 
Without loss of generality, we may suppose that 
$D_1 \subset C_1$. 
By the minimality of $\vert S \cap P \vert$, 
we see that $D_1$ is a $\gamma$--meridian disk of $C_1$. 
Let $A_2$ be the closure of the component of $S-P$ 
such that $A_2 \cap D_1 \ne \emptyset$. 

\medskip
\noindent
{\bf Claim 3}\qua $A_2$ is $\gamma$--incompressible in $C_2$. 

\begin{proof}
Suppose that there is a $\gamma$--compressing disk $D$ for $A_2$ in $C_2$. 
If $D \cap \gamma = \emptyset$, then we have a contradiction as in 
the proof of Claim~1. 
Suppose that $D \cap \gamma \ne \emptyset$. 
Let $D_2$ be the disk obtained from $A_2$ by $\gamma$--compressing 
$A_2$ along $D$ such that $\partial D_2 = \partial D_1$. 
Since $\partial D_1$ is $\gamma$--essential in $P$, 
this shows that $D_2$ is a $\gamma$--meridian disk of $C_2$. 
Hence $C_1 \cup C_2$ is $\gamma$--reducible, a contradiction. 
\end{proof}

Note that $\partial A_2 \subset \partial_+ C_2$. 
There is a $\gamma$--boundary compressing disk $\Delta$ 
for $A_2$ in $C_2$ such that 
$\Delta \cap \partial C_2 \subset \partial_+ C_2$ 
(Corollary~\ref{gamma-incomp., boundary incomp. surface in gamma-comp. body}). 
By the minimality of $\vert S \cap P \vert$, 
we see that $A_2$ is not $\gamma$--parallel to a surface 
in $\partial_+ C_2$. 
Hence, by $\gamma$--boundary compressing $A_2$ along $\Delta$, 
and applying a tiny isotopy, 
we obtain a $\gamma$--meridian disk $D_2$ in $C_2$ such that 
$D_1 \cap D_2 = \emptyset$. 
Hence $C_1 \cup_P C_2$ is weakly $\gamma$--reducible, 
a contradiction. 

This completes the proof of Proposition~\ref{Haken's theorem}. 

\subsection*{3.C\qua Proof of Proposition~\ref{Casson-Gordon's lemma}}

Let $N$, $\gamma$ be as in Proposition~\ref{Casson-Gordon's lemma}
and 
$C_1 \cup_P C_2$ a Heegaard splitting of $(N, \gamma )$. 
Let $\hat{N} = \text{cl}(N-N( \gamma ))$, and 
$\hat{C_1} \cup_{\hat{P}} \hat{C_2}$ the Heegaard splitting of 
$(\hat{N}, \hat{F}_1, \hat{F}_2)$ 
obtained from $C_1 \cup_P C_2$ as in Section~3.A. 
Let $D$ be a $\gamma$--compressing disk for $F_1 \cup F_2$. 

\medskip
\noindent
{\bf Case 1}\qua $D \cap \gamma = \emptyset$. 

\medskip
In this case, we may regard that $D$ is a disk in $\hat{N}$. 
It is clear that $D$ is a compressing disk of 
$\hat{F}_1 \cup \hat{F}_2$. 
Hence, by \cite[Lemma 1.1]{CG}, 
we see that $\hat{C_1} \cup_{\hat{P}} \hat{C_2}$ is 
weakly reducible. 
This implies that $C_1 \cup_P C_2$ is weakly $\gamma$--reducible. 

\medskip
\noindent
{\bf Case 2}\qua $D \cap \gamma \ne \emptyset$ 
(ie, $D \cap \gamma$ consists of a point). 

\medskip
Let $\hat{D} = \text{cl}(D-N( \gamma ))$. 

\medskip
\noindent
{\bf Claim}\qua $\hat{D}$ is an essential annulus in $\hat{N}$. 

\begin{proof}
By using the argument as in Claim~1 of Case~2 of Section~3.B, 
we can show that $\hat{D}$ is incompressible in $\hat{N}$. 
Suppose that $\hat{D}$ is parallel to an annulus ${\cal A}$ 
in $\partial \hat{N}$. 
Let $s$, $F_i'$ $(i=1,2)$ be as in Claim~2 of Case~2 of Section~3.B. 
Let $A^*$ be the component of $\text{Fr}_N N( \gamma )$ such that 
$\partial D \subset A^*$, and 
$F^*$ the component of $F_1' \cup F_2'$ such that 
$F^* \supset \partial D$. 
By using the argument of the proof of Claim~2 of Case~2 of Section~3.B, 
we see that ${\cal A}$ is disjoint from 
$s \cup (\text{Fr}_N N( \gamma ) - A^*)$, hence 
$\text{cl}({\cal A} - A^*) \subset F^*$. 
Hence $F^* \cap {\cal A}$ is an annulus, and this shows that 
$\partial D$ bounds a $\gamma$--disk in $F_1 \cup F_2$, 
a contradiction. 
\end{proof}

Suppose, for a contradiction, that 
$C_1 \cup_P C_2$ is strongly $\gamma$--irreducible. 
By Lemma~\ref{strongly gamma-irreducible, strongly irreducible}, 
$\hat{C_1} \cup_{\hat{P}} \hat{C_2}$ is strongly irreducible. 
Then, by \cite[Lemma 6]{Schu} or \cite[Lemma 2.3]{Mo'}, 
$\hat{D}$ is ambient isotopic rel~$\partial$ to a surface $\hat{D}'$ 
such that $\hat{D}' \cap \hat{P}$ consists of essential simple closed 
curves in $\hat{D}'$. 
We regard $\hat{D} = \hat{D}'$. 
This means that each component of $D \cap P$ is a simple closed curve 
bounding a disk in $D$, which contains the point $D \cap \gamma$. 
We suppose that $\vert D \cap P \vert$ is minimal 
among the $\gamma$--compressing disks for $F_1 \cup F_2$ with this property. 
Let $n = \vert D \cap P \vert$. 

Suppose that $n=1$, ie, $D \cap P$ consists of a 
simple closed curve, say $\ell_1$. 
Then the closures of the components of $D - \ell_1$ consists of 
a disk, say $D_1$, and an annulus, say $A_2$. 
Without loss of generality, we may suppose that 
$D_1 \subset C_1$, and $A_2 \subset C_2$. 
Note that a component of $\partial A_2$ is contained in 
$\partial_- C_2$, and the other in $\partial_+ C_2$. 
Since $C_2$ is not trivial, there exists a $\gamma$--meridian 
disk $D_2$ in $C_2$. 
It is elementary to show, by applying cut and paste arguments on
$D_2$ and $A_2$, that there is a $\gamma$--meridian disk 
$D_2'$ in $C_2$ such that 
$D_2' \cap A_2 = \emptyset$. 
Hence $D_1 \cap D_2' = \emptyset$, and 
this shows that $C_1 \cup_P C_2$ is weakly $\gamma$--reducible, 
a contradiction. 

Suppose that $n \ge 2$. 
Let $D_1$ be the closure of the component of $D-P$ such that 
$D_1 \cap \gamma \ne \emptyset$. 
Note that $D_1$ is a $\gamma$--disk. 
Without loss of generality, we may suppose that 
$D_1 \subset C_1$. 
By the minimality of $\vert D \cap P \vert$, 
we see that $D_1$ is a $\gamma$--meridian disk of $C_1$. 
Let $A_2$ be the closure of the component of $D-P$ 
such that $A_2 \cap D_1 \ne \emptyset$. 
Then, 
by using the arguments as in the proof of Claim~3 of Case~2 of 
Section~3.B, we can show that $A_2$ is $\gamma$--incompressible in $C_2$. 
Note that $\partial A_2 \subset \partial_+ C_2$. 
There is a $\gamma$--boundary compressing disk $\Delta$ 
for $A_2$ in $C_2$ such that 
$\Delta \cap \partial C_2 \subset \partial_+ C_2$ 
(Corollary~\ref{gamma-incomp., boundary incomp. surface in gamma-comp. body}). 
By the minimality of $\vert S \cap P \vert$, 
we see that $A_2$ is not $\gamma$--parallel to a surface 
in $\partial_+ C_2$. 
Hence, by $\gamma$--boundary compressing $A_2$ along $\Delta$, 
and applying a tiny isotopy, 
we obtain a $\gamma$--meridian disk $D_2$ of $C_2$ such that 
$D_1 \cap D_2 = \emptyset$. 
Hence $C_1 \cup_P C_2$ is weakly $\gamma$--reducible, 
a contradiction. 

This completes the proof of Proposition~\ref{Casson-Gordon's lemma}. 

\subsection*{3.D\qua Proof of Proposition~\ref{Rubinstein-Scharlemann's lemma}}

Let ${\cal D}$ be a union of mutually disjoint, 
non $K$--parallel, $K$--meridian disks for $X$ such that 
${\cal D} \cap A = \emptyset$. 
We suppose that ${\cal D}$ is maximal among the 
unions of $K$--meridian disks 
with the above properties. 
Let $Z'= N(\partial X,X) \cup N({\cal D},X)$. 
Then we have the following two cases. 

\medskip
\noindent
{\bf Case 1}\qua  A component of $\partial Z' - \partial X$ bounds a 
$K$--ball, say $B_K$, such that $B_K \supset A$. 

\medskip
In this case, since $\partial B_K \subset B$, and $B$ is $K$--irreducible, 
$\partial B_K$ bounds a $K$--ball $B_K'$ in $B$ 
(Corollary~\ref{gamma-compression body is gamma-irreducible}). 
Hence $M= B_K \cup B_K'$ is a 3--sphere. 
In particular, if $K \ne \emptyset$, then 
$K \cap B_K$ ($K \cap B_K'$ respectively) is a trivial arc 
properly embedded in $B_K$ ($B_K'$ respectively). 
Hence $K$ is a trivial knot. 
This shows that we have conclusion~1. 

\medskip
\noindent
{\bf Case 2}\qua  No component of $\partial Z' - \partial X$ bounds a $K$-ball 
which contains $A$. 

\medskip
Since $X$ is $K$--irreducible, 
each of the $K$--sphere components of $\partial Z' - \partial X$ 
(if exists) bounds $K$--balls in $X$. 
By the construction of $Z'$, 
it is easy to see that the $K$--balls are mutually disjoint. 
Let $Z = Z' \cup \text{(the $K$--balls)}$. 
By (3) of Remark~\ref{remark of compression body} and 
Proposition~\ref{2-fold cover of compression body}, 
we see that 
$Z$ is a $K$--compression body with 
$\partial_+ Z = \partial X$, and 
by the maximality of ${\cal D}$, we see that
$\partial_- Z$ consists of one component, say $F$, such that 
$F$ bounds a $K$--handlebody which contains $A$. 
Let $N = Y \cup Z$. 
Note that $Y \cup_Q Z$ is a Heegaard splitting of $(N, K \cap N)$. 
Since $\partial N = F$ is a closed surface 
contained in $B$, it is $K$--compressible in $B$
(Proposition~\ref{2-fold cover, incompressible}). 
By the maximality of ${\cal D}$, 
we see that the compressing disk lies in $N$. 
Hence, by Proposition~\ref{Casson-Gordon's lemma}, we see that 
$Y \cup_Q Z$ is weakly $K$--reducible. 
This obviously implies that $X \cup_Q Y$ is weakly $K$--reducible, 
and we have conclusion~2. 

This completes the proof of Proposition~\ref{Rubinstein-Scharlemann's lemma}. 

\section{The Casson--Gordon theorem}

A Casson and C\,McA Gordon proved that if a Heegaard 
splitting of a closed 3--manifold $M$ is weakly reducible, then 
either the splitting is reducible, or 
$M$ contains an incompressible surface \cite[Theorem 3.1]{CG}. 
In this section, 
we generalize this result for compact $M$. 
The author thinks that this generalization is well known 
(eg, \cite{Schu'}). 
However, the formulation given here will be useful for 
the proof of Theorem~\ref{main theorem} (Section~7.C). 

Let $M$ be a compact, orientable 3--manifold, and 
$C_1 \cup_P C_2$ a Heegaard splitting of $M$ such that 
$P$ is a closed surface, 
ie, $\partial_- C_1 \cup \partial_- C_2 = \partial M$. 
Let 
$\Delta = \Delta_1 \cup \Delta_2$ 
be a weakly reducing collection of disks for $P$, 
ie, 
$\Delta_i$ $(i=1,2)$ is a union of mutually disjoint, 
non-empty meridian disks of $C_i$ such that 
$\Delta_1 \cap \Delta_2 = \emptyset$. 
Then $P( \Delta )$ denotes the surface obtained from $P$ 
by compressing $P$ along $\Delta$. 
Let 
$\hat{P}( \Delta ) = P( \Delta ) - ($the components of 
$P( \Delta )$ 
which are contained in $C_1$ or $C_2)$. 

\begin{lemma}\label{2-sphere implies reducible}
If there is a 2--sphere component in 
$\hat{P}( \Delta )$, then $C_1 \cup_P C_2$ is reducible. 
\end{lemma}

\begin{proof}
Suppose that there is a 2--sphere component $S$ of 
$\hat{P}( \Delta )$. 
We note that $S \cap C_i$ $(i=1,2)$ is a union of 
non-empty meridian disks of $C_i$. 
Let 
$\hat{S} = \text{cl}(S-(C_1 \cup C_2))$. 
Note that $\hat{S}$ is a planar surface in $P$. 
Let ${\cal A}_1 \cup {\cal A}_2$ be a union of mutually 
disjoint arcs properly embedded in $\hat{S}$ 
such that 
$\partial {\cal A}_i \subset \partial (S \cap C_i)$, 
and that 
$\text{cl}(\hat{S}-N({\cal A}_1 \cup {\cal A}_2, \hat{S}))$ 
is an annulus, say $A'$. 
Let $S'$ be a 2--sphere obtained from $S$ by pushing ${\cal A}_1$ 
into $C_1$, and ${\cal A}_2$ into $C_2$ such that $S' \cap P =A'$. 
It is clear that $S' \cap C_i$ $(i=1,2)$ consists of a disk, say $D_i$, 
obtained from $S \cap C_i$ by banding along ${\cal A}_i$. 

\medskip
\noindent
{\bf Claim}\qua $D_i$ is a meridian disk of the compression body $C_i$ $(i=1,2)$. 

\begin{proof}
Suppose, for a contradiction, that either $D_1$ or $D_2$, 
say $D_1$, is not a meridian disk, 
ie, there is a disk $D$ in $P$ such that 
$\partial D = \partial D_1$. 
Note that we have either 
$N({\cal A}_1, \hat{S}) \subset D$, or 
$N({\cal A}_1, \hat{S}) \subset \text{cl}(P-D)$. 
If 
$N({\cal A}_1, \hat{S}) \subset D$, 
then $\partial (S \cap C_1)$ is recovered from 
$\partial D$ by banding along arcs properly embedded in $D$. 
This shows that $\partial (S \cap C_1) \subset D$, 
and this implies that each component of $S \cap C_1$ is not a 
meridian disk, a contradiction. 
On the other hand, if 
$N({\cal A}_1, \hat{S}) \subset \text{cl}(P-D)$, 
then $\text{cl}(\hat{S} - N({\cal A}_1, \hat{S})) \subset D$. 
This shows that $\partial (S \cap C_2) \subset D$, 
and this implies that each component of $S \cap C_2$ 
is not a meridian disk, a contradiction. 
\end{proof}

Since $\partial D_1$ and $\partial D_2$ are parallel in $P$, 
we see by Claim that $C_1 \cup_P C_2$ is reducible. 
\end{proof}

Now we define a complexity $c(F)$ of a closed surface $F$ as follows. 
$$c(F) = \sum (\chi (F_i) -1),$$ 
where the sum is taken for all components of $F$. 
Then we suppose that $c (\hat{P}( \Delta ))$ is maximal 
among all weakly reducing collections of disks for $P$. 
By Lemma~\ref{2-sphere implies reducible}, 
we see that if the complexity of a 
component of $\hat{P}( \Delta )$ is positive, then 
$C_1 \cup_P C_2$ is reducible. 
Suppose that the complexities of the components of 
$\hat{P}( \Delta )$ are strictly negative, 
ie, 
each component of $\hat{P}( \Delta )$ is not a 2--sphere. 
Then, by the argument of the proof of \cite[Theorem 3.1]{CG}, 
we see that $\hat{P}( \Delta )$ 
is incompressible in $M$. 
Hence we have the next proposition. 

\begin{proposition}\label{CG Theorem 3.1}
Let $M$ be a compact, orientable 3--manifold, 
and $C_1 \cup_P C_2$ a Heegaard splitting of $M$ with 
$\partial_- C_1 \cup \partial_- C_2 = \partial M$. 
Suppose that $C_1 \cup_P C_2$ is weakly reducible. 
Then either 

\begin{enumerate}
\item 
$C_1 \cup_P C_2$ is reducible, or 

\item
there exists a weakly reducing collection of disks 
$\Delta$ for $P$ such that each component of 
$\hat{P}( \Delta )$ is an incompressible surface in $M$, 
which is not a 2--sphere. 
\end{enumerate}
\end{proposition}

Note that, in \cite{CG}, $M$ is assumed to be closed. 
However, it is easy to see that 
the arguments there work for Heegaard splittings 
$C_1 \cup C_2$ such that 
$\partial_-C_1 \cup \partial_-C_2 = \partial M$. 

The following is a slight extension of \cite[Lemme 1.4]{BO}. 
Let $M$, $C_1 \cup_P C_2$, $\Delta$ be as above. 
Suppose that we have conclusion~2 of Proposition~\ref{CG Theorem 3.1}. 
Let $M_1, \dots , M_n$ be the closures of the components of 
$M - \hat{P}( \Delta )$. 
Let $M_{j,i} = M_j \cap C_i$ 
$(j=1, \dots , n, i = 1,2 )$. 

\begin{lemma}\label{HS from weakly reducing collection of disks}
For each $j$, we have either one of the following. 

\begin{enumerate}

\item 
$M_{j,2} \cap P \subset \text{Int} (M_{j,1} \cap P)$, 
and 
$M_{j,1}$ is connected. 

\item 
$M_{j,1} \cap P \subset \text{Int} (M_{j,2} \cap P)$, 
and 
$M_{j,2}$ is connected. 

\end{enumerate}
\end{lemma}

\begin{proof}
Recall that $\Delta_i$ is the union of the components of $\Delta$ 
that are contained in $C_i$ $(i = 1, 2)$. 
We see, from the definition of $\hat{P}( \Delta )$, 
that each $M_j$ is obtained as in the following manner. 

\begin{quote}
(*) Take a component $N$ of 
$\text{cl}(C_i - N( \Delta_i, C_i))$ ($i=1$ or $2$, say $1$) 
such that there exists a component $D_2$ of $\Delta_2$ such that 
$\partial D_2 \subset N$. 
Let 
$N' = N \cup ($the components of $N( \Delta_2, C_2)$ intersecting $N)$. 
Then 
$M_j = N' \cup($the union of components $N_2$ of 
$\text{cl}(C_2 - N( \Delta_2, C_2))$ 
such that $(N_2 \cap P) \subset (N \cap P))$. 
\end{quote}

It is clear that this construction process gives conclusion~1. 
If $N$ is a component of 
$\text{cl}(C_2 - N( \Delta_2, C_2))$, 
then we have conclusion~2. 
\end{proof}

We note that each component of 
$\text{Fr}_{C_i} (M_{j,i})$ is a meridian disk of $C_i$, 
which is parallel to a component of $\Delta$. 
Recall that $\hat{P}( \Delta )$ is obtained from 
$P( \Delta )$ by discarding the components each of which 
is contained in $C_1$ or $C_2$. 
These imply that each component $E$ of $M_{j,i}$ 
inherits a compression body structure from $C_i$
(see (3) of Remark~\ref{remark of compression body}), 
ie, 
$\partial_+ E = (E \cap \partial_+ C_i) \cup \text{Fr}_{C_i} (E)$. 
Then we can obtain a splitting, 
denoted by $C_{j,1} \cup_{P_j} C_{j,2}$, of $M_j$ as follows (\cite[Lemme 1.4]{BO}). 

\leftskip 28.5pt

Suppose that $M_j$ satisfies conclusion~1 (2 respectively) of 
Lemma~\ref{HS from weakly reducing collection of disks}. 
Recall that $M_{j,1}$ ($M_{j,2}$ respectively) inherits a compression body 
structure from $C_1$ ($C_2$ respectively). 
Then let\hfill\break 
$C_{j,1} = \text{cl}(M_{j,1} - N( \partial_+ M_{j,1}, M_{j,1}))$ 
($C_{j,2} = \text{cl}(M_{j,2} - N( \partial_+ M_{j,2}, M_{j,2}))$ 
respectively), 
and\hfill\break
$C_{j,2} = N( \partial_+ M_{j,1}, M_{j,1}) \cup M_{j,2}$ 
($C_{j,1} = N( \partial_+ M_{j,2}, M_{j,2}) \cup M_{j,1}$ 
respectively). 

\leftskip0pt

\begin{lemma}\label{C from M}
Suppose that each component of $\hat{P}( \Delta )$ is not a 2--sphere. 
Then each $C_{j,i}$ is a compression body such that, 
for each $j$, we have 
$\partial_+ C_{j,2} = \partial_+ C_{j,1} = C_{j,1} \cap C_{j,2}$, 
ie, 
$C_{j,1} \cup_{P_j} C_{j,2}$ 
is a Heegaard splitting of $M_j$. 
\end{lemma}

\begin{proof}
Since the argument is symmetric, we may suppose that 
$M_j$ satisfies conclusion~1 of 
Lemma~\ref{HS from weakly reducing collection of disks}. 
Since $M_{j,1}$ is a compression body, it is clear that 
$C_{j,1}$ is a compression body. 
Let ${\cal D}_1 = \text{Fr}_{C_2} M_{j,2}$. 
There is a union of mutually disjoint meridian disks, 
say ${\cal D}_2$, of $M_{j,2}$ such that 
${\cal D}_2 \cap {\cal D}_1 = \emptyset$, 
and each component of the manifold obtained from $M_{j,2}$ by 
cutting along ${\cal D}_2$ is homeomorphic to either a 3--ball 
or $G \times [0,1]$, where $G$ is a component of $\partial_- M_{j,2}$ 
with $G \times \{ 0 \}$ corresponding to $G$. 
Hence $C_{j,2}$ $(= N( \partial_+ M_{j,1}, M_{j,1}) \cup M_{j,2})$ 
is homeomorphic to a manifold obtained from
$N( \partial_+ M_{j,1}, M_{j,1})$ 
$( \cong \partial_+ M_{j,1} \times [0,1])$ 
by attaching 2--handles along the simple closed curves 
corresponding to $\partial ({\cal D}_1 \cup {\cal D}_2)$ 
in $\partial_+ M_{j,1} \times \{ 1 \}$, 
and capping off some of the resulting 2--sphere boundary components. 
By the definition of compression body (Section~2.A), 
this implies that $C_{j,2}$ is a compression body, unless 
there exists a 2--sphere component $S$ of $\partial C_{j,2}$, 
which is disjoint from 
$N( \partial_+ M_{j,1}, M_{j,1}) \cap C_{j,1}$
$(= \partial_+ M_{j,1} \times \{ 0 \})$. 
However such $S$ must be a component of $\hat{P}( \Delta )$, 
a contradiction. 
\end{proof}

Let $M$, $C_1 \cup_P C_2$, $\Delta$, $M_j$, $M_{j,i}$, and 
$C_{j,1} \cup_{P_j} C_{j,2}$ be as above. 

\begin{lemma}\label{non trivial compression body}
Suppose that each component of $\hat{P}( \Delta )$ is not a 2--sphere. 
If $\partial M$ is incompressible in $M$, 
then each compression body $C_{j,i}$ is not trivial. 
\end{lemma}

\begin{proof}
Suppose that some compression body is trivial. 
By changing subscripts if necessary, we may suppose that 
$C_{1,1}$ is trivial. 
Then we claim that $M_{1,2} \cap P \subset \text{Int}(M_{1,1} \cap P)$, 
ie, 
we have conclusion~1 of 
Lemma~\ref{HS from weakly reducing collection of disks}. 
In fact, if we have conclusion~2 of 
Lemma~\ref{HS from weakly reducing collection of disks}, 
then 
$C_{1,1} = N(\partial_+ M_{1,2}, M_{1,2}) \cup M_{1,1}$. 
However this expression obviously implies $C_{1,1}$ is not trivial, 
a contradiction. 
Hence 
$C_{1,1} = \text{cl}(M_{1,1} - N(\partial_+ M_{1,1}, M_{1,1}))$, 
and this implies that $M_{1,1}$ is a trivial compression body 
such that 
$\partial_- M_{1,1}$ is a component of $\partial M$. 
Let $D$ be any component of $\text{Fr}_{C_2} M_{1,2}$. 
Then by extending $D$ vertically to $M_{1,1}$, we obtain a disk 
$\tilde{D}$ properly embedded in $M$. 
Since each component of $\hat{P}( \Delta )$ 
is not a 2--sphere, $\partial \tilde{D}$ is not contractible in 
$\partial M$. 
Hence $\tilde{D}$ is a compressing disk of $\partial M$, 
a contradiction. 
\end{proof}

\begin{lemma}\label{reducible implies reducible}
If some $C_{j,1} \cup_{P_j} C_{j,2}$ is reducible, then 
$C_1 \cup_P C_2$ is reducible. 
\end{lemma}

\begin{proof}
We prove this by using an argument of C.Frohman \cite[Lemma 1.1]{F}. 
If $M$ is reducible, then by \cite[Lemma 1.1]{CG}, 
we see that any Heegaard splitting of $M$ is reducible. 
Hence we may suppose that $M$ is irreducible. 
If a component of $\hat{P}( \Delta )$ is a 2--sphere, 
then $C_1 \cup_P C_2$ is reducible (Lemma~\ref{2-sphere implies reducible}). 
Hence we may suppose that each component of $\hat{P}( \Delta )$ 
is not a 2--sphere 
(hence, $C_{j,1} \cup_{P_j} C_{j,2}$ is a Heegaard splitting of $M_j$). 
Since the argument is symmetric, we may suppose that 
the pair $M_{j,1}$, $M_{j,2}$ satisfies conclusion~1 
of Lemma~\ref{HS from weakly reducing collection of disks}.  
By \cite[Lemma 1.1]{CG}, there exists an incompressible 2--sphere $S$ 
in $M_j$ such that $S$ intersects $P_j$ in a circle. 
Let $D_1 = S \cap C_{j,1}$. 
Note that $D_1$ is a meridian disk of $C_{j,1}$. 
Since $M$ is irreducible, $S$ bounds a 3--ball $B^3$ in $M_j$. 
Let $C_{j,1}'$ be the closure of the component of $C_{j,1}-D_1$ 
such that $C_{j,1}' \subset B^3$. 
Since $\partial M \cap B^3 = \emptyset$, 
we see that $C_{j,1}'$ is a handlebody, 
ie, 
$\partial_- C_{j,1}' = \emptyset$. 
Let $X$ be a spine of $C_{j,1}'$, 
and 
$M_X = \text{cl}(M-N(X, C_1))$. 
It is clear that $S$ is an incompressible 2--sphere in $M_X$, 
and $P$ is a Heegaard surface of $M_X$. 
Hence, by \cite[Lemma 1.1]{CG}, 
there exists an incompressible 2--sphere $S_X$ in $M_X$ such that 
$S_X$ intersects $P$ in a circle. 
It is obvious that the 2--sphere $S_X$ gives a reducibility of 
$C_1 \cup_P C_2$. 
\end{proof}

\section{Reducing genus $g$, $n$--bridge positions}

\begin{figure}[ht!]
\begin{center}
\includegraphics[width=6cm, clip]{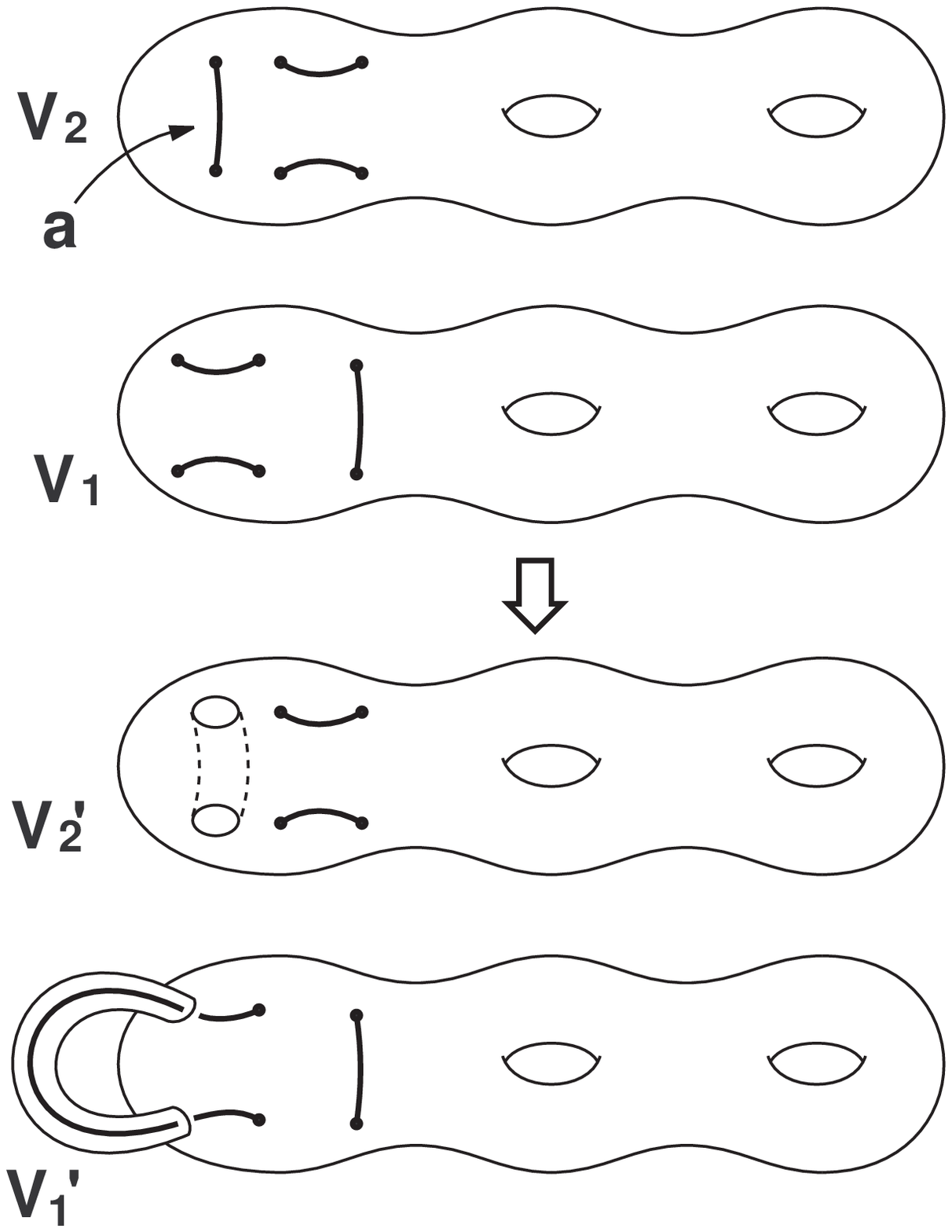}
\end{center}
\caption{}\label{5.1}
\end{figure}

Let $K$ be a knot in a closed, orientable 3--manifold $M$. 
Let $V_1 \cup V_2$ be a Heegaard splitting of $M$, 
which gives a genus $g$, $n$--bridge position of $K$ with 
$n \ge 1$. 
Let $a$ be a component of $K \cap V_i$ 
($i=1$ or $2$, say $2$). 
Let $V_1' = V_1 \cup N(a, V_2)$, 
and $V_2' = \text{cl} (V_2-N(a, V_2))$. 
By the definition of genus $g$, $n$--bridge positions, 
it is easy to see that 
$V_1' \cup V_2'$ gives a genus $(g+1)$, $(n-1)$--bridge position 
of $K$. 
We say that the Heegaard splitting 
$V_1' \cup V_2'$ is obtained from 
$V_1 \cup V_2$ by a {\it tubing} (along $a$). 
See Figure~\ref{5.1}. 

We say that a knot $K$ in $M$ is a {\it core knot} if there is a genus one 
Heegaard splitting $V \cup W$ of $M$ such that $K$ is a core curve 
of the solid torus $V$, 
ie, $K$ admits a genus one, $0$--bridge position. 
Note that if $M$ is a 3--sphere, then $K$ is a core knot if and only 
if $K$ is a trivial knot. 
We say that $K$ is {\it small} if the exterior $E(K)$ 
of $K$ does not contain a closed essential surface. 
We say that a surface $F$ properly embedded in $E(K)$ is 
{\it meridional} if $\partial F$ is a union of non-empty meridian loops. 
We note that \cite[Theorem 2.0.3]{CGLS} implies that 
if $M$ is a 3--sphere and $K$ is small, then $E(K)$ does not contain a 
meridional essential surface. 

\begin{proposition}\label{tubing}
Let $K$ be a knot in a closed orientable 3--manifold $M$ 
with the following properties. 

\begin{enumerate}

\item 
$M$ is $K$--irreducible. 

\item 
There exists a 2--fold branched covering space of $M$ with
branch set $K$. 

\item 
$K$ is not a core knot. 

\item 
$K$ is small and there does not exists a meridional essential surface 
in $E(K)$. 
\end{enumerate}

Let $C_1 \cup_P C_2$ be a Heegaard splitting of $M$, 
which gives a genus $g$, $n$--bridge position of $K$. 
Suppose that $C_1 \cup_P C_2$ is weakly $K$--reducible. 
Then we have either one of the following. 

\begin{enumerate}

\item 
There exists a 
weakly $K$--reducing pair of disks $E_1$, $E_2$ in $C_1$, $C_2$ 
respectively such that 
$E_1 \cap K = \emptyset$, and 
$E_2 \cap K = \emptyset$. 

\item 
There exists a Heegaard splitting $H_1 \cup H_2$ of 
$M$, which gives a genus $(g-1)$ $(n+1)$--bridge position 
of $K$ such that $C_1 \cup C_2$ is obtained from $H_1 \cup H_2$ 
by a tubing. 

\end{enumerate}
\end{proposition}

\begin{remark}\label{remark of tubing}
Note that, in Proposition~\ref{tubing}, 
if $g=0$, then we always have conclusion~1. 
\end{remark}

\begin{proof}
Let $D_1$, $D_2$ be a pair of $K$--essential disks 
in $C_1$, $C_2$ respectively, which gives a weak 
$K$--reducibility of $C_1 \cup C_2$. 
If $D_1 \cap K = \emptyset$ and $D_2 \cap K = \emptyset$, 
then we have conclusion~1. 
Hence in the rest of the proof, we may suppose that 
$D_1 \cap K \ne \emptyset$. 
We have the following two cases. 

\medskip
\noindent 
{\bf Case 1}\qua $D_2 \cap K = \emptyset$. 

\medskip
In this case, we first show the following. 

\medskip
\noindent
{\bf Claim 1}\qua If $D_1$ is separating in $C_1$, 
then we have conclusion~1. 

\begin{proof}
Let $C_1'$, $C_1''$ be the closures of the components of $C_1 - D_1$ 
such that $\partial D_2 \subset \partial C_1''$. 
Then $C_1'$ is a $K$--handlebody which is not a $K$--ball. 
Hence there exists a $K$--essential disk $D_1'$ in $C_1'$ such that 
$D_1' \cap K = \emptyset$, and $D_1' \cap D_1 = \emptyset$ 
(hence, $D_1'$ is properly embedded in $C_1$). 
See Figure~\ref{5.2}. 
It is clear that $D_1'$ is a $K$--meridian disk of $C_1$. 
Hence, by regarding, 
$E_1 = D_1'$, $E_2 = D_2$, we have conclusion~1. 
\end{proof}

\begin{figure}[ht!]
\begin{center}
\includegraphics[width=6cm, clip]{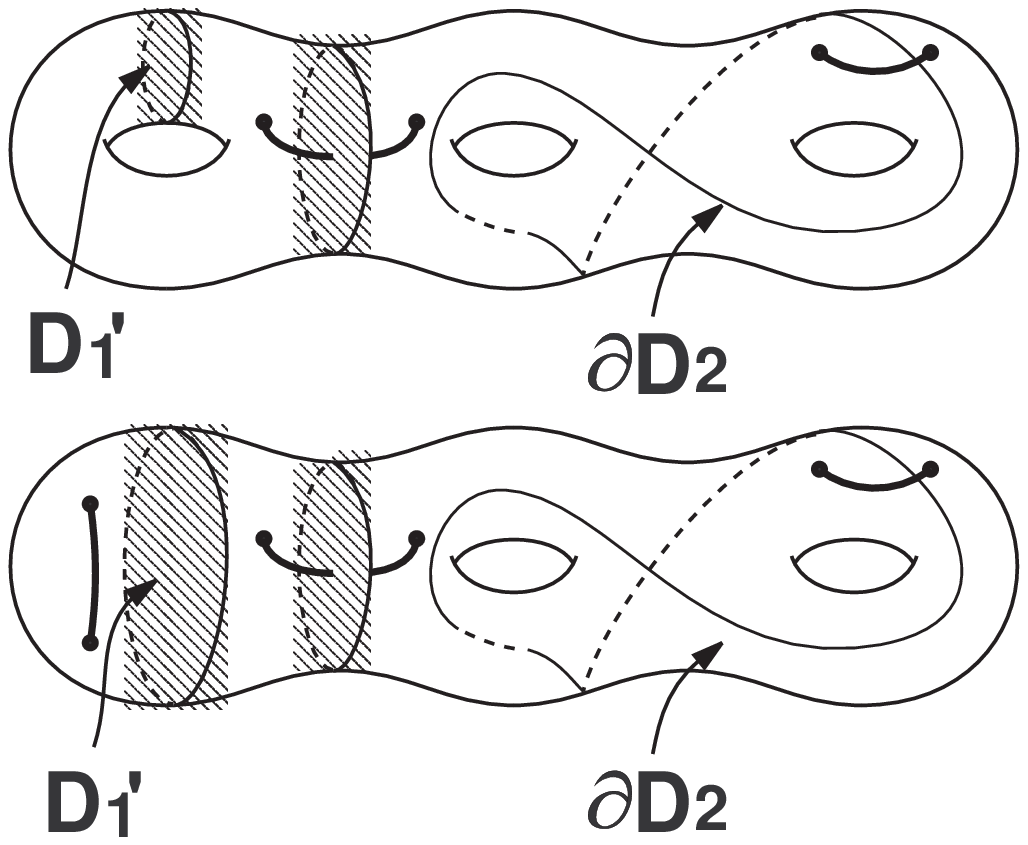}
\end{center}
\caption{}\label{5.2}
\end{figure}

By Claim~1, we may suppose that $D_1$ is non-separating in $C_1$. 
Let $P'$ be the surface obtained from $P$ by $K$--compressing along $D_1$, 
and 
$\hat{P'} = P' \cap E(K)$. 
We note that $P'$ separates $M$ into two components, 
say $C_1'$ and $C_2'$, where $C_1'$ is obtained from $C_1$ 
by cutting along $D_1$. 
Let $\hat{C}_i' = C_i' \cap E(K)$ $(i=1,2)$. 
Hence, by Section~3.B, we see that $\hat{C}_1'$ is a compression body 
with $\partial_+ \hat{C}_1' = \hat{P}'$. 
Let ${\cal D}$ be a union of maximal mutually disjoint, 
non parallel compressing disks for $\hat{P}'$ such that 
${\cal D} \subset \hat{C}_2'$. 
Note that since $D_2$ is a compressing disk for $\hat{P}'$ such that 
$D_2 \subset \hat{C}_2'$, there actually exists such ${\cal D}$. 
Let $\hat{C}_2'' = N( \hat{P}', \hat{C}_2') \cup N( {\cal D}, \hat{C}_2')$. 
Note that $\hat{C}_2'$ is homeomorphic to 
$C_2 \cap E(K) = \text{cl}(C_2-N(K))$, 
hence, $\hat{C}_2'$ is irreducible (Section~3.A). 
Hence the 2--sphere components ${\cal S}$ (possibly ${\cal S}=\emptyset$) of 
$\partial \hat{C}_2'' - \hat{P}'$ bounds mutually disjoint 3--balls 
in $\hat{C}_2'$. 
Let $C_2^* = \hat{C}_2'' \cup \text{(the 3--balls)}$. 
Then $C_2^*$ is a compression body such that 
$\partial_+ C_2^* = \hat{P}'$. 
Let $P^* = \partial_- C_2^*$. 

\medskip
\noindent
{\bf Claim 2}\qua If $P^*$ is compressible in $E(K)$, 
then we have conclusion~1. 

\begin{proof}
Suppose that there exists a compressing disk $E$ of $P^*$ in $E(K)$. 
Let $M^* = \hat{C}_1' \cup C_2^*$. 
By the maximality of ${\cal D}$, we see that $E$ is contained in $M^*$. 
Note that $\hat{C}_1' \cup_{\hat{P}'} C_2^*$ 
is a Heegaard splitting of $M^*$, 
and $\partial E \subset \partial_- C_2^*$. 
Hence, by \cite[Lemma 1.1]{CG}, we see that 
$\hat{C}_1' \cup_{\hat{P}'} C_2^*$ is weakly reducible, and this implies 
conclusion~1. 
\end{proof}

By Claim~2, we may suppose that $P^*$ is incompressible in $E(K)$. 
Note that each component of $\partial P^*$ is a meridian loop of $K$. 
Since $K$ is small and there does not exist a meridional essential surface 
in $E(K)$, we see that each component of $P^*$ is a 
boundary parallel annulus properly embedded in $E(K)$. 
Recall that ${\cal S}$ is the union of the 2--sphere components of 
$\partial \hat{C}_2'' - \hat{P}'$. 
Note that we can assign labels 
$C_1$ and $C_2$ to the components of $E(K) - (P^* \cup {\cal S})$  alternately 
so that the $C_2$ region are contained in $\hat{C}_2'$, 
and that 
$\hat{P}'$ is recovered from $P^* \cup {\cal S}$ by adding tubes 
along mutually disjoint arcs in $C_1$--regions. 
Recall that 
$\hat{P}'$ is connected. 
Since each component of $P^* \cup {\cal S}$ is separating in
$E(K)$, this shows that exactly one component of 
$E(K) - (P^* \cup {\cal S})$ is a $C_1$--region. 
Let $\tilde{P}^*$ be a surface in $M$ obtained from $P^*$ 
by capping off the boundary components by mutually disjoint 
$K$--disks in $N(K)$ 
(hence, via isotopy, 
$P'$ is recovered from $\tilde{P}^* \cup {\cal S}$ by adding the tubes 
used for recovering $\hat{P}'$ from $P^* \cup {\cal S}$). 
Then each component of $\tilde{P}^* \cup {\cal S}$ is a $K$--sphere. 
Since $M$ is $K$--irreducible, the components of $\tilde{P}^* \cup {\cal S}$ 
bounds $K$--balls, say $B_1, \dots , B_m$, in $M$. 

\medskip
\noindent
{\bf Claim 3}\qua The $K$--balls $B_1, \dots , B_m$ are mutually disjoint. 

\begin{proof}
Suppose not. 
By exchanging the subscript if necessary, we may suppose that 
$B_2 \subset B_1$. 
Since there exists exactly one $C_1$--region, 
this implies that the $K$--balls 
$B_2, \dots , B_m$ are included in $B_1$ in a non-nested configuration. 
Hence $P'$ is contained is the $K$--ball $B_1$. 
See Figure~\ref{5.3}. 
Note that $P$ is recovered from $P'$ by adding a tube 
along the component of $K - P'$, 
which intersects $D_1$. 
Hence we see that $P$ is contained in a regular neighborhood of $K$, 
say $N_K$. 

\begin{figure}[ht!]
\begin{center}
\includegraphics[width=4cm, clip]{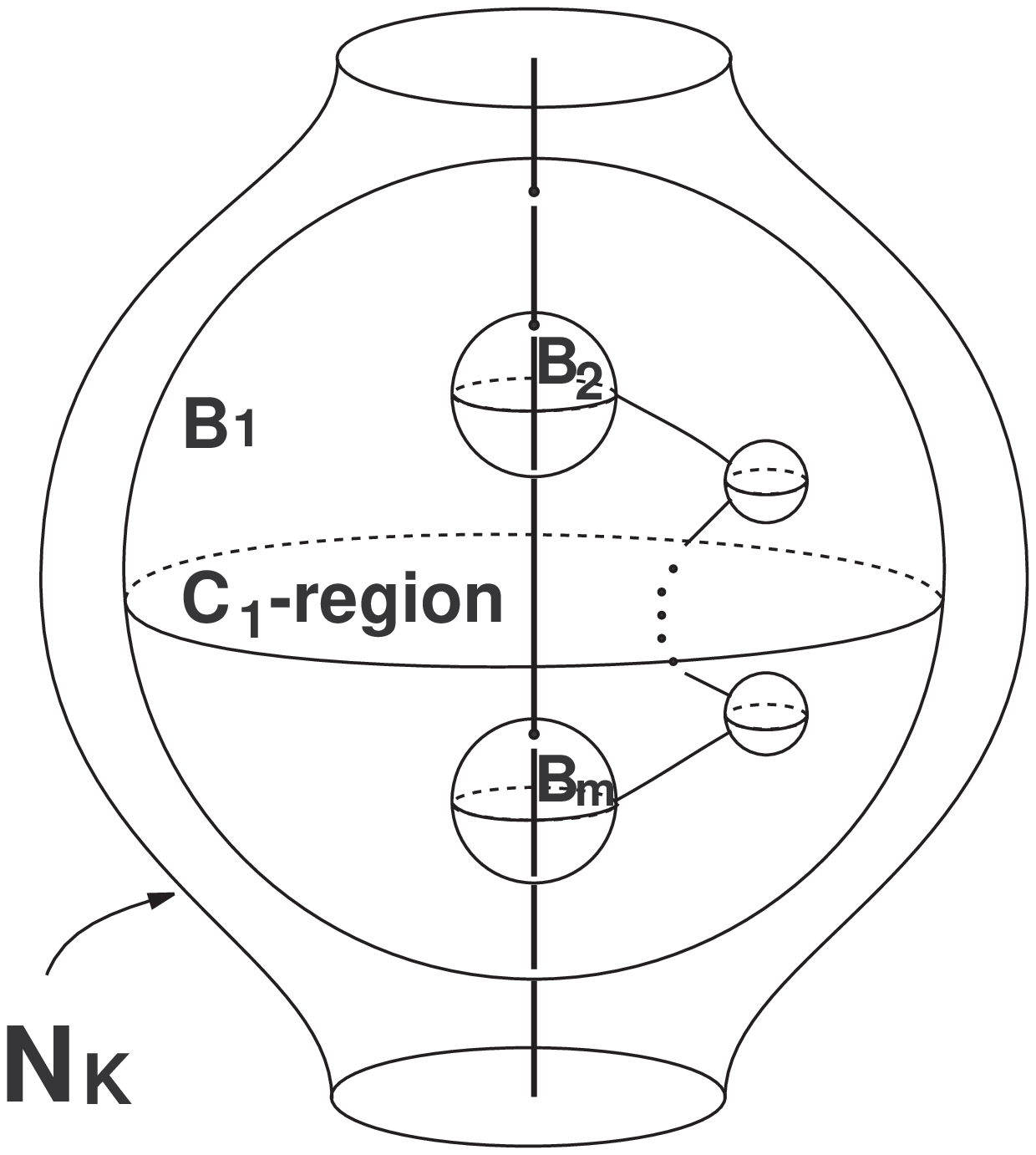}
\end{center}\vspace{5pt}
\caption{}\label{5.3}
\end{figure}

Note that $\text{cl}(M - N_K)$ is contained in $C_2$. 
Since $C_2$ is a $K$--compression body, there exists a 
$K$--compressing disk $D_N$ for $\partial N_K$ in $C_2$. 
Suppose that $D_N \subset N_K$. 
Since $N_K$ is a regular neighborhood of $K$, 
we see that $\text{cl} (N_K - N(D_N, N_K))$ is a $K$--ball. 
Since $C_2$ is $K$--irreducible, 
$\text{cl} (M-N_K) \cup N(D_N, N_K)$ is also a $K$--ball. 
These show that $M$ is the 3--sphere, and $K$ is a trivial knot, 
contradicting the condition~3 of the assumption of 
Proposition~\ref{tubing}. 
Suppose that $D_N \subset \text{cl} (M-N_K)$. 
Since $M$ is $K$--irreducible, 
we see that $\text{cl} (M-N_K)$ is irreducible. 
This shows that we obtain a 3--ball by cutting 
$\text{cl} (M-N_K)$ along $D_N$. 
This shows that $\text{cl} (M-N_K)$ is a solid torus. 
Hence $N_K \cup \text{cl} (M-N_K)$ is a genus one Heegaard splitting 
of $M$. 
Hence $K$ is a core knot, 
contradicting the condition~3 of the assumption of 
Proposition~\ref{tubing}. 

This completes the proof of Claim~3.
\end{proof}

Recall that $P'$ is the surface obtained from $P$ by $K$--compressing 
along $D_1$, 
and $C_1'$, $C_2'$ the closures of the components of $M-P'$, 
where $C_1'$ is obtained from $C_1$ by cutting along $D_1$. 
By Proposition~\ref{2-fold cover of compression body} and 
(3) of Remark~\ref{remark of compression body}, $C_1'$ is a $K$--handlebody. 
By Claim~3, we see that the $C_1$--region is 
$\text{cl}(M - (\cup_{i=1}^m B_i)) \cap E(K)$. 
Hence $P'$ is recovered from 
$\partial B_1 \cup \cdots \cup \partial B_m$ by adding 
tubes along arcs properly embedded in 
$\text{cl} (M - (\cup_{i=1}^m B_i))$. 
Hence, we see that $C_2'$ is obtained from the $K$--balls 
$B_1, \dots , B_m$ by adding 1--handles disjoint from $K$. 
Hence $C_2'$ is also a $K$--handlebody. 
These show that $P'$ is a Heegaard surface for $(M,K)$. 
It is clear that $C_1' \cup C_2'$ gives 
a genus $(g-1)$, $(n+1)$ bridge position 
of $K$, and $C_1 \cup C_2$ is obtained from $C_1' \cup C_2'$ 
by a tubing along a component of $K \cap (\cup B_i)$. 

Hence, by regarding $H_1 = C_1'$, $H_2 = C_2'$, we have 
conclusion~2 of Proposition~\ref{tubing}. 

\medskip
\noindent 
{\bf Case 2}\qua $D_2 \cap K \ne \emptyset$. 

\medskip
In this case, we first show the following. 

\medskip
\noindent
{\bf Claim 1}\qua If $\partial D_i$ $(i=1$ or $2)$ is separating in $P$, 
then we have conclusion of Proposition~\ref{tubing}. 

\begin{proof}
Since the argument is symmetric, 
we may suppose that $\partial D_2$ is separating in $P$. 
This implies that $D_2$ is separating in $C_2$. 
Let $C_2'$, $C_2''$ be the closures of the components of 
$C_2 - D_2$ such that $\partial D_1 \subset C_2''$. 
Then $C_2'$ is a $K$--handlebody which is not a $K$--ball. 
Hence there exists a $K$--essential disk $D_2'$ in $C_2'$ such that 
$D_2' \cap K = \emptyset$, and 
$D_2' \cap D_2 = \emptyset$ 
(hence, $D_2'$ is properly embedded in $C_2$). 
See Figure~\ref{5.2}. 
It is clear that $D_2'$ is a $K$--meridian disk of $C_2$. 
Hence by applying the arguments of Case~1 to the pair 
$D_1$, $D_2'$, we have conclusion of Proposition~\ref{tubing}. 
\end{proof}

Let $P'$ be the surface obtained from $P$ by $K$--compressing along 
$D_1 \cup D_2$, 
and $\hat{P}' = P' \cap E(K)$. 
Let $C_1'$, $C_2'$ be the closures of the components of $M-P'$ 
such that 
$C_1'$ is obtained from $C_1$ by cutting along $D_1$ and 
attaching $N(D_2,C_2)$, and 
$C_2'$ is obtained from $C_2$ by cutting along $D_2$ and 
attaching $N(D_1,C_1)$. 
Then let $\hat{C}_i' = C_i' \cap E(K)$ $(i=1,2)$. 

\medskip
\noindent
{\bf Claim 2}\qua 
If $\hat{P}'$ is compressible in $E(K)$, then we have
conclusion of Proposition~\ref{tubing}. 

\begin{proof}
Suppose that there is a compressing disk $D$ for $\hat{P}'$ 
in $E(K)$. 
Since the argument is symmetric, we may suppose that 
$D \subset \hat{C}_2'$. 
We may regard that $D$ is a compressing disk for $P'$. 
Since $P$ is recovered from $P'$ by adding two tubes along a 
component of $K \cap C_1'$ and a component of $K \cap C_2'$, 
we may suppose that $D \cap P = \partial D$. 
Hence $D$ is a $K$--meridian disk of $C_2$ 
such that $D \cap K = \emptyset$. 
Hence, by applying the arguments of Case~1 to the pair 
$D_1$, $D$, we have the conclusion of Proposition~\ref{tubing}. 
\end{proof}

By Claims~1 and 2, we see that, for the proof of Proposition~\ref{tubing}, 
it is enough to show that either 
(1) $\partial D_i$ $(i=1$ or $2)$ is separating in $P$, or 
(2) $\hat{P}'$ is compressible in $E(K)$. 
Suppose that $\partial D_i$ ($i=1, 2$) is non-separating in $P$, 
and that $\hat{P}'$ is incompressible in $E(K)$. 
Then, by the argument preceding Claim~3 of Case~1, 
we see that each component of $P'$ is a $K$--sphere, and 
$P$ is recovered from $P'$ by adding tubes along 
two arcs $a_1$, $a_2$ such that $a_i$ is a component of 
$K \cap C_i'$ $(i=1,2)$, and that $a_1 \cap a_2 = \emptyset$. 
Note that $P$ is connected. 
Since $\partial D_1$, $\partial D_2$ 
are non-separating in $P$, 
we see that $P'$ consists of one $K$--sphere, or 
two $K$--spheres, and 
this shows that $K \cap C_i'$ consists of one arc, or 
two arcs. 
But since $K$ is a knot, we have $a_1 \cap a_2 \ne \emptyset$ in 
either case, a contradiction. 
Hence we have the conclusion of Proposition~\ref{tubing} in Case~2. 

This completes the proof of Proposition~\ref{tubing}. 
\end{proof}

\section{Heegaard splittings of $(S^3, \text{two bridge knot})$}

In this section, we prove the following. 

\begin{proposition}\label{HS of 2-bridge knots are weakly K-reducible} 
Let $K$ be a non-trivial two bridge knot, and 
$X \cup_Q Y$ a Heegaard splitting of $S^3$, 
which gives a genus $g$, $n$--bridge position of $K$. 
Suppose that $(g, n) \ne (0,2)$. 
Then $X \cup_Q Y$ is weakly $K$--reducible. 
\end{proposition}

\begin{proposition}\label{HS of 2-bridge knots are weakly reducible} 
Let $K$ be a non-trivial two bridge knot. 
Then, for each $g \ge 3$, 
every genus $g$ Heegaard splitting of the exterior $E(K)$ of $K$ 
is weakly reducible. 
\end{proposition}

\subsection*{6.A\qua Comparing $X \cup_Q Y$ with a two bridge position}

Let $A \cup_P B$ be a genus $0$ Heegaard splitting of $S^3$, 
which gives a 2--bridge position of $K$. 
Then, by \cite[Corollary 6.22]{KS} (if $n \ge 1$) or 
by \cite[Corollary 3.2]{K} (if $n = 0$), we have the following. 

\begin{proposition}\label{comparing with a 2-bridge position}
Let $X \cup_Q Y$ be a Heegaard splitting of $S^3$, 
which gives a genus $g$, $n$--bridge position of $K$. 
If $X \cup_Q Y$ is strongly $K$--irreducible, 
then $Q$ is $K$--isotopic to a position such that 
$P\cap Q$ consists of non-empty collection of transverse 
simple closed curves which are $K$--essential in both 
$P$ and $Q$. 
\end{proposition}

In this subsection, we prove the following proposition. 

\begin{proposition}\label{essential-essential intersection}
Let $X \cup_Q Y$ be a Heegaard splitting of $S^3$, 
which gives a genus $g$, $n$--bridge position of $K$ 
with $(g, n) \ne (0,2)$. 
Suppose that $P \cap Q$ consists of non-empty collection of 
transverse simple closed curves which are 
$K$--essential in both $P$ and $Q$. 
Then $X \cup_Q Y$ is weakly $K$--reducible. 
\end{proposition}

We note that Proposition~\ref{HS of 2-bridge knots are weakly K-reducible} 
is a consequence of 
Propositions~\ref{comparing with a 2-bridge position} 
and \ref{essential-essential intersection}. 

\begin{proof}
[Proof of Proposition~\ref{HS of 2-bridge knots are weakly K-reducible}
from Propositions~\ref{comparing with a 2-bridge position} and \ref{essential-essential intersection}]

Let $X \cup_Q Y$\break be a Heegaard splitting of $S^3$, 
which gives a genus $g$, $n$--bridge position of $K$ 
with $(g, n) \ne (0,2)$. 
Suppose, for a contradiction, 
that $X \cup_Q Y$ is strongly $K$--irreducible. 
Then, by Propositions~\ref{comparing with a 2-bridge position}, 
we may suppose that $P\cap Q$ consists of non-empty collection of transverse 
simple closed curves which are $K$--essential in both $P$ and $Q$. 
By Propositions~\ref{essential-essential intersection}, 
we see that $X \cup_Q Y$ is weakly $K$--reducible, 
a contradiction. 
\end{proof}

\begin{proof}[Proof of Proposition~\ref{essential-essential intersection}]
First of all, we would like to remark that the proof given below 
is just an orbifold version of the proof of \cite[Corollary 6.4]{RS}. 
We suppose that $\vert P \cap Q \vert$ is minimal among 
all surfaces $P$ such that $P$ gives a two bridge position of $K$, 
and that 
$P \cap Q$ consists of non-empty collection of simple closed 
curves which are $K$--essential in both $P$ and $Q$. 
Note that the closure of each component of $P-Q$ is 
either an annulus which is disjoint from $K$, 
or a disk intersecting $K$ in two points. 
We divide the proof into several cases. 

\begin{figure}[ht!]
\begin{center}
\includegraphics[width=4cm, clip]{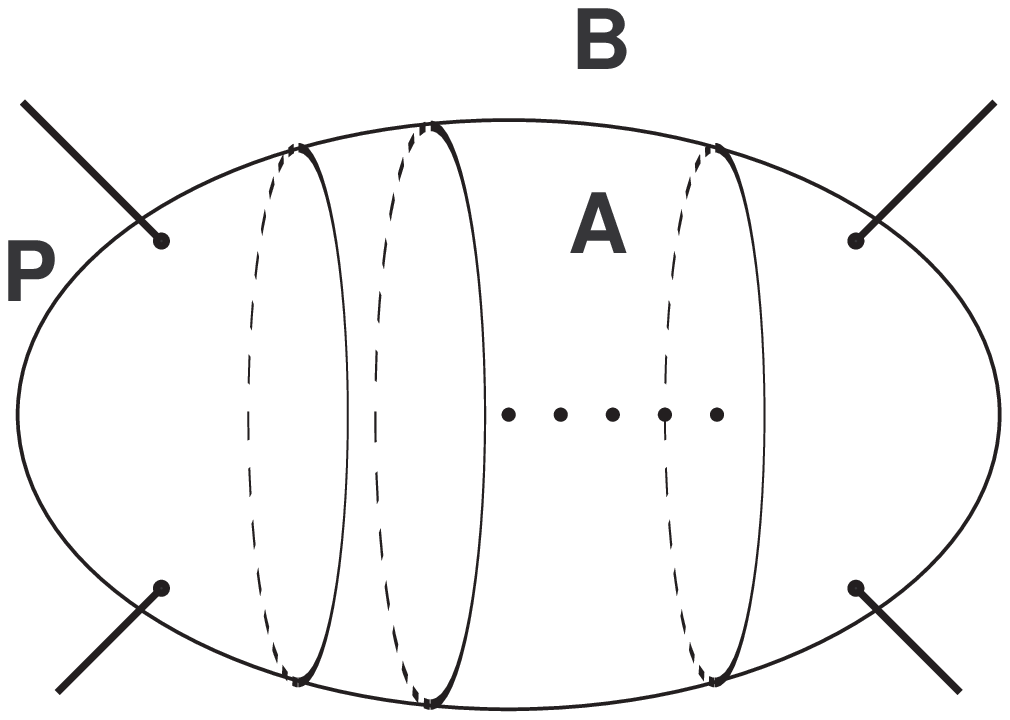}
\end{center}\vspace{5pt}
\caption{}\label{6.1.1}
\end{figure}

\medskip
\noindent
{\bf Case 1}\qua Each component of $P \cap X$ is not $K$--boundary parallel in $X$, 
and each component of $P \cap Y$ is not $K$--boundary parallel in $Y$. 

\medskip
Case~1 is divided into the following subcases. 

\medskip
\noindent
{\bf Case 1.1}\qua $P \cap X$ contains a component which is $K$--compressible in $X$, 
and $P \cap Y$ contains a component which is $K$--compressible in $Y$

\medskip
In this case, by $K$--compressing the components in $X$ and $Y$, 
we obtain $K$--meridian disks $D_X$, $D_Y$ in $X$, $Y$ respectively. 
By applying a slight $K$--isotopy if necessary, we may suppose that 
$D_X \cap D_Y = \emptyset$, and this shows that $X \cup_Q Y$ is 
weakly $K$--reducible. 

\medskip
\noindent
{\bf Case 1.2}\qua Either $P \cap X$ or $P\cap Y$, say $P \cap X$, 
contains a component which is $K$--compressible in $X$, 
and each component of $P \cap Y$ is $K$--incompressible in $Y$. 

\medskip
Let $D_X$ be a $K$--meridian disk obtained by 
$K$--compressing the component of $P \cap X$. 
Note that $P \cap Y$ is $K$--boundary compressible 
in $Y$ (Corollary~\ref{incomp., boundary incomp. surface in gamma-comp.}).   
Let $D_Y$ be a disk in $Y$ obtained by $K$--boundary 
compressing $P \cap Y$. 
By the minimality of $\vert P \cap Q \vert$, 
we see that $D_Y$ is a $K$--meridian disk in $Y$. 
Since $\partial D_X$ is a component of $P \cap Q$, we may suppose that 
$D_X \cap D_Y = \emptyset$ by applying a slight $K$--isotopy 
if necessary. 
This shows that $X \cup_Q Y$ is weakly $K$--reducible. 

\medskip
\noindent
{\bf Case 1.3}\qua Each component of $P \cap X$ is $K$--incompressible in $X$, 
and each component of $P \cap Y$ is $K$--incompressible in $Y$. 

\medskip
Let $D_X$ ($D_Y$ respectively) be a $K$--meridian disk obtained by 
$K$--boundary compressing $P \cap X$ ($P \cap Y$ respectively). 
By applying slight isotopies, we may suppose that 
$\partial D_X \cap P = \emptyset$, 
$\partial D_Y \cap P = \emptyset$ 
(hence, 
$\partial D_X \subset A$ or $B$, 
$\partial D_Y \subset A$ or $B$). 
If one of $\partial D_X$ or $\partial D_Y$ is contained in $A$, 
and the other in $B$, 
then $D_X \cap D_Y = \emptyset$, and this shows that 
$X \cup_Q Y$ is weakly $K$--reducible. 
Suppose that 
$\partial D_X \cup \partial D_Y$ is contained in $A$ or $B$, 
say $A$. 
Let $D_B$ be a $K$--meridian disk in $B$
(ie, 
$D_B$ is a disk properly embedded in $B$ such that 
$D_B \cap K = \emptyset$, and 
$D_B$ separates the components of $K \cap B$). 
Note that since each component of $P \cap X$, $P \cap Y$ 
is $K$--incompressible, $D_B \cap Q \ne \emptyset$. 
We take $D_B$ so that $\vert D_B \cap Q \vert$ 
is minimal among all $K$--essential disks $D'$ in $B$  
such that each component of 
$D' \cap (P \cap X)$ ($D' \cap (P \cap Y)$ respectively) 
is a $K$--essential arc properly embedded in 
$P \cap X$ ($P \cap Y$ respectively). 
Suppose that $D_B \cap Q$ contains a simple closed curve component. 
Let $D^* ( \subset D_B)$ be an innermost disk. 
Since the argument is symmetric, we may suppose that 
$D^* \subset X$. 
By the minimality of $\vert D_B \cap Q \vert$, 
we see that $D^*$ is a $K$--meridian disk in $X$. 
Since $D^* \subset B$, 
$\partial D^* \cap \partial D_Y = \emptyset$. 
Hence the pair $D^*$, $D_Y$ gives a weak $K$--reducibility 
of $X \cup_Q Y$. 
Suppose that each component of $D_B \cap Q$ is an arc. 
Let $\Delta (\subset D_B)$ be an outermost disk. 
Since the argument is symmetric, 
we may suppose that $\Delta \subset X$. 
Recall that $\Delta \cap (P \cap X)$ is a $K$--essential 
arc in $P \cap X$. 
By the minimality of $\vert D_B \cap Q \vert$, we see that at least 
one component, say $D^{**}$, of the surface obtained from 
$P \cap X$ by $K$--boundary compressing along $\Delta$ is a 
$K$--meridian disk in $X$. 
Since $\Delta \subset B$, 
$\partial D^{**} \cap \partial D_Y = \emptyset$. 
Hence the pair $D^{**}$, $D_Y$ gives a weak $K$--reducibility 
of $X \cup_Q Y$. 

\medskip
\noindent
{\bf Case 2}\qua A component of $P \cap X$ or $P \cap Y$, say $P \cap Y$, 
is $K$--boundary parallel in $Y$. 

\medskip
By the minimality of $\vert P \cap Q \vert$, 
we have either 
$\vert P \cap Q \vert = 1$ 
(and $P \cap Y$ ($P \cap X$ respectively) is a disk intersecting $K$ 
in two points) or, 
$\vert P \cap Q \vert = 2$
(and $P \cap Y$ is an annulus disjoint from $K$). 

\medskip
\noindent
{\bf Case 2a}\qua $\vert P \cap Q \vert = 1$. 

\medskip
Let 
$P_X = P \cap X$ and $P_Y = P \cap Y$. 
Let $E$ be the closure of the component of $Q-P$ such that 
$E$ and $P_Y$ are $K$--parallel in $Y$. 
Since the argument is symmetric, we may suppose that 
$E \subset A$. 
We have the following subcases. 

\medskip
\noindent
{\bf Case 2a.1}\qua $P_X$ is $K$--boundary parallel in $X$ 

\medskip
Since $(g, n) \ne (0,2)$, 
$P_X$ is parallel to $E$ in $A$, and cannot be parallel to 
$\text{cl} (Q-E)$. 

\begin{figure}[ht!]
\begin{center}
\includegraphics[width=5cm, clip]{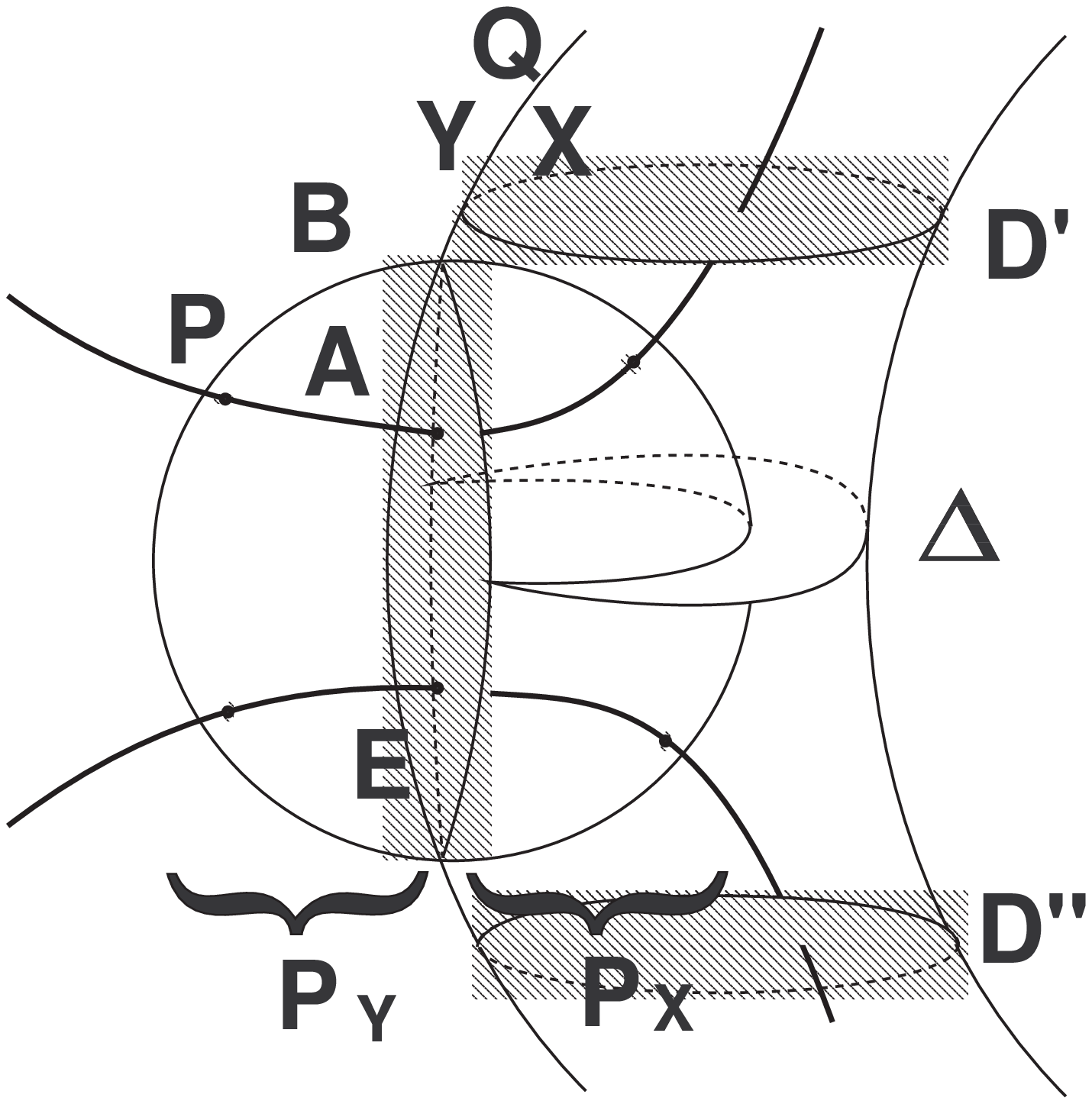}
\end{center}
\caption{}\label{6.1.2}
\end{figure}

Let $D_B$ be a $K$--meridian disk in $B$. 
Since $K$ is not a trivial knot, 
$\partial D_B$ and $\partial E$ are not isotopic in $P-K$. 
Hence $D_B \cap Q \ne \emptyset$. 
We suppose that 
$\vert D_B \cap Q \vert$ 
is minimal among all $K$--meridian disks $D'$ in $B$ 
such that each component of $D' \cap P_X$ ($D' \cap P_Y$ respectively) 
is a $K$--essential arc in $P_X$ ($P_Y$ respectively). 
Suppose that $D_B \cap Q$ contains a simple closed curve. 
Let $D^* ( \subset D_B)$ be an innermost disk. 
Since the argument is symmetric, we may suppose that 
$D^* \subset X$. 
By the minimality of $\vert D_B \cap Q \vert$, 
we see that $D^*$ is a $K$--meridian disk in $X$. 
Then by pushing $D_Y$ into $X$ along the parallelism through $E$, 
we can $K$--isotope $P$ to $P'$ such that 
$P' \subset \text{Int} X$, and 
$P' \cap D^* = \emptyset$. 
Hence, by Proposition~\ref{Rubinstein-Scharlemann's lemma}, 
we see that $X \cup_Q Y$ is weakly $K$--reducible. 
Suppose that each component of $D_B \cap Q$ is an arc. 
Let $\Delta (\subset D_B)$ be an outermost disk. 
Since the argument is symmetric, we may suppose that 
$\Delta \subset X$. 
See Figure~\ref{6.1.2}. 

\medskip
\noindent
{\bf Claim}\qua At least one component of the disks 
obtained from $P_X$ by $K$--boundary compressing along $\Delta$ 
is a $K$--meridian disk. 

\begin{proof}
Let $D'$, $D''$ be the disks obtained from 
$P_X$ by $K$--boundary compressing along $\Delta$. 
Suppose that $D'$ is $K$--boundary parallel, 
ie, there exists a $K$--disk $D_Q$ in $Q$ such that 
$\partial D_Q = \partial D'$. 
Note that since $D' \cup D''$ is obtained from 
$P_X$ by $K$--boundary compressing along $\Delta$, 
there is an annulus $A_Q$ is $Q$ such that 
$\partial A_Q = \partial D' \cup \partial D''$, and that 
$A_Q \cap K = E \cap K$: two points. 
Note also that $D_Q \cap K$ consists of one point. 
Hence $A_Q$ is not contained in $D_Q$, and this implies that 
$A_Q \cap D_Q = \partial D'$. 
Then $A_Q \cup D_Q$ is a disk intersecting $K$ in three points, 
whose boundary is $\partial D''$. 
Since $(g, n) \ne (0,2)$, $\text{cl}(Q-(A_Q \cup D_Q))$ is not 
a $K$--disk. 
Hence $D''$ is a $K$--meridian disk in $X$. 
\end{proof}

Let $D''$ be a $K$--meridian disk in $X$ obtained as in Claim. 
By applying a slight isotopy, we may suppose that $P \cap D'' = \emptyset$. 
Then by pushing $P_Y$ into $X$ along the parallelism through $E$, 
we can $K$--isotope $P$ to $P'$ such that 
$P' \subset \text{Int} X$, and 
$P' \cap D'' = \emptyset$. 
Hence, by Proposition~\ref{Rubinstein-Scharlemann's lemma}, 
we see that $X \cup_Q Y$ is weakly $K$--reducible. 

\medskip
\noindent
{\bf Case 2a.2}\qua $P_X$ is not $K$--boundary parallel in $X$, 
and $P_X$ is $K$--incompressible in $X$, 
ie, $P_X$ is $K$--essential in $X$. 

\medskip
Since $P_X$ is $K$--incompressible, there is a $K$--boundary 
compressing disk $\Delta$ for $P_X$ in $X$. 

\medskip
\noindent
{\bf Claim}\qua $\Delta \subset B$. 

\begin{proof}
Suppose that $\Delta \subset A$. 
Note that $K \cap E$ consists of two points in $\text{Int} E$, 
and $\Delta \cap E$ is an arc properly embedded in $E$, which 
separates the points. 
Then, by $K$--boundary compressing $P_X$ along $\Delta$, 
we obtain two $K$--disks.  
Since $X$ is $K$--irreducible, these $K$--disks are 
$K$--boundary parallel in $X$. 
This shows that $P_X$ is $K$--boundary parallel in $X$, 
contradicting the condition of Case~2a.2. 
\end{proof}

Then, by using the argument of the proof of Claim of Case~2a.1, 
we see that at least one component, say $D''$, 
of the $K$--disks obtained from $P_X$ by $K$--boundary compressing
along $\Delta$ is a $K$--meridian disk in $X$. 
By applying a slight isotopy, 
we may suppose that $D'' \cap P = \emptyset$. 
By Claim, we see that $\partial D'' \subset B$. 
Then by pushing $P_Y$ into $X$ along the parallelism through $E$, 
we can $K$--isotope $P$ to $P'$ such that 
$P' \subset \text{Int} X$, and 
$P' \cap D'' = \emptyset$. 
Hence, by Proposition~\ref{Rubinstein-Scharlemann's lemma}, 
we see that $X \cup_Q Y$ is weakly $K$--reducible. 

\medskip
\noindent
{\bf Case 2a.3}\qua $P_X$ is not $K$--boundary parallel in $X$, 
and $P_X$ is $K$--compressible in $X$. 

\medskip
Let $D$ be the $K$--compressing disk for $P_X$. 
Since there does not exist a 2--sphere $(\subset S^3)$ intersecting $K$ 
in three points, $D \cap K = \emptyset$. 
Let $D^*$ be the disk component of a surface obtained from $P_X$ 
by $K$--compressing along $D$. 
Since $(g,n) \ne (0,2)$, 
we see that $D^*$ is a $K$--meridian disk of $X$. 
By applying a slight isotopy, 
we may suppose that $D^* \cap P = \emptyset$. 
Suppose that $D^* \subset B$. 
Then by pushing $D_Y$ into $X$ along the parallelism through $E$, 
we can $K$--isotope $P$ to $P'$ such that 
$P' \subset \text{Int} X$, and 
$P' \cap D^{*} = \emptyset$. 
Hence, by Proposition~\ref{Rubinstein-Scharlemann's lemma}, 
we see that $X \cup_Q Y$ is weakly $K$--reducible. 
Hence, in the rest of this subcase, we suppose 
that $D^* \subset A$ (Figure~\ref{6.1.3}). 
Let $D_B$ be a $K$--meridian disk in $B$. 
Since $K$ is not a trivial two component link, 
$\partial D$ and $\partial D_B$ are not isotopic in $P-K$. 
Hence $D_B \cap Q \ne \emptyset$. 
We suppose that $\vert D_B \cap Q \vert$ is 
minimal among all $K$--meridian disks $D'$ in $B$ 
such that each component of $D' \cap P_X$ ($D' \cap P_Y$ respectively) 
is a $K$--essential arc in $P_X$ ($P_Y$ respectively). 

\begin{figure}[ht!]
\begin{center}
\includegraphics[width=4cm, clip]{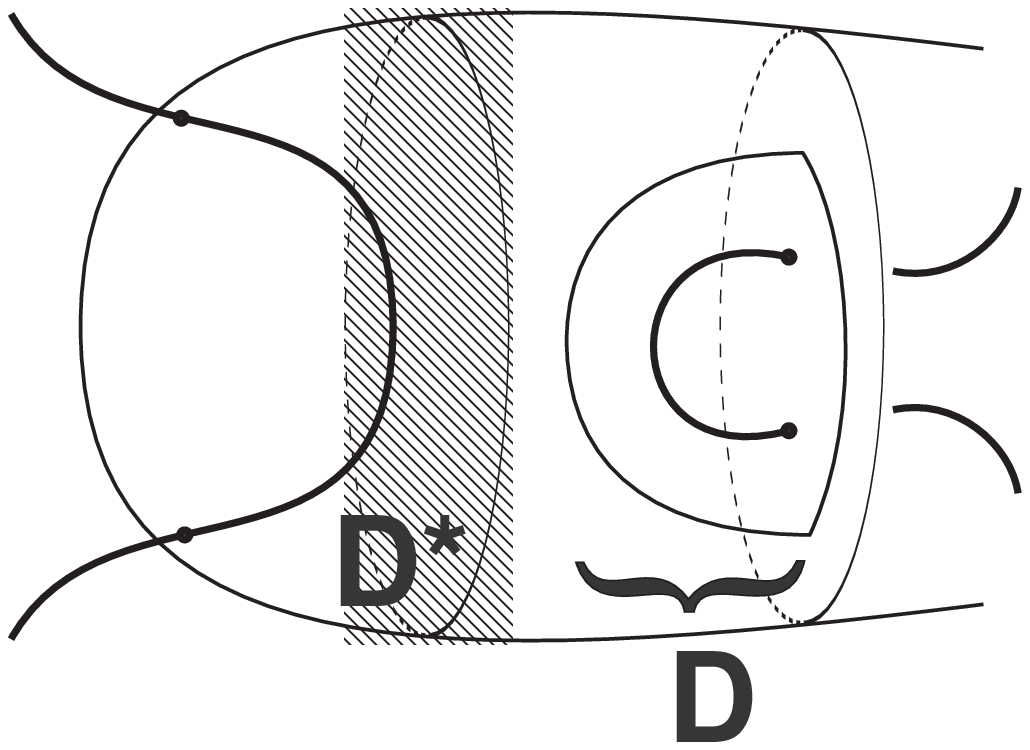}
\end{center}
\caption{}\label{6.1.3}
\end{figure}

Suppose that $D_B \cap Q$ contains a simple closed curve. 
Let $D^{**} ( \subset D_B)$ be an innermost disk. 
By the minimality of $\vert D_B \cap Q \vert$, 
we see that $\partial D^{**}$ is $K$--essential in $Q$. 
Note that $\partial D^{**} \subset B$. 
If $D^{**} \subset Y$, 
then the pair $D^{**}$, $D^*$ gives a weak $K$--reducibility of 
$X \cup_Q Y$. 
If $D^{**} \subset X$, 
then by pushing $P_Y$ into $X$ along the parallelism through $E$, 
we can $K$--isotope $P$ to $P'$ such that 
$P' \subset \text{Int} X$, and 
$P' \cap D^{**} = \emptyset$. 
Hence, by Proposition~\ref{Rubinstein-Scharlemann's lemma}, 
we see that $X \cup_Q Y$ is weakly $K$--reducible. 

Suppose that each component of $D_B \cap Q$ is an arc. 
Let $\Delta (\subset D_B)$ be an outermost disk. 
If $\Delta \subset X$, 
then by using the argument as in the proof of Case~2a.1, 
we see that $X \cup_Q Y$ is weakly $K$--reducible. 
Suppose that $\Delta \subset Y$. 
Then, by using the argument as in the proof of Claim of 
Case~2a.1, we can show that at least one component , say $D''$, 
of the $K$--disks obtained from $P_Y$ by $K$--boundary compressing 
along $\Delta$ is a $K$--meridian disk in $Y$. 
By applying slight $K$--isotopy, 
we may suppose that $D'' \subset B$. 
Hence the pair $D^*$, $D''$ gives a weak $K$--reducibility 
of $X \cup_Q Y$. 

\medskip
\noindent
{\bf Case 2b}\qua $\vert P \cap Q \vert = 2$. 

\medskip
Let $D_1$, $D_2$ be the components of $P\cap X$, 
and $A_1 = P \cap Y$. 
Recall that $A_1$ is a $K$--boundary parallel annulus in $Y$ 
such that $A_1 \cap K = \emptyset$, 
and that 
$D_1$, $D_2$ are not $K$--boundary parallel. 
We also note that $\partial D_1 \cup \partial D_2$ bounds 
an annulus $A'$ in $Q$ such that 
$A_1$ and $A'$ are $K$--parallel in $Y$. 
Without loss of generality, we may suppose that $A'$ is 
contained in the 3--ball $A$. 

\medskip
\noindent
{\bf Case 2b.1}\qua $D_1 \cup D_2$ is $K$--incompressible in $X$. 

\medskip
Since $D_1 \cup D_2$ is $K$--incompressible, 
there is a $K$--boundary compressing disk $\Delta$ for 
$D_1 \cup D_2$. 
Without loss of generality, we may suppose that 
$\Delta \cap D_1 \ne \emptyset$, 
$\Delta \cap D_2 = \emptyset$. 
Since $D_1$ is not $K$--boundary parallel, 
at least one component, say $D^{**}$, of the $K$--disks obtained 
from $D_1$ by $K$--boundary compressing along $\Delta$ 
is a $K$--meridian disk in $X$. 
By applying a slight $K$--isotopy, 
we may suppose that $D^{**} \cap P = \emptyset$. 

\medskip
\noindent
{\bf Claim}\qua $D^{**} \subset B$. 

\begin{proof}
Suppose, for a contradiction, that $D^{**} \subset A$. 
Then $\partial D^{**}$ is contained in the annulus $A'$ 
bounded by $\partial D_1 \cup \partial D_2$. 
We note that $D^{**}$ intersects $K$ in one point. 
Hence $\partial D^{**}$ is not contractible in $Q$. 
This shows that $\partial D^{**}$ is a core curve of $A'$. 
Let $A''$ be the annulus in $A'$ bounded by 
$\partial D^{**} \cup \partial D_1$. 
Then the 2--sphere $D_1 \cup A'' \cup D^{**}$ 
intersects $K$ in three points, a contradiction. 
\end{proof}

By Claim we see that, 
by pushing $A_1$ into $X$ along the parallelism through $A'$, 
we can $K$--isotope $P$ to $P'$ such that 
$P' \subset \text{Int} X$. 
By the above claim, we may suppose that $P' \cap D^{**} = \emptyset$. 
Hence, by Proposition~\ref{Rubinstein-Scharlemann's lemma}, 
we see that $X \cup_Q Y$ is weakly $K$--reducible. 

\medskip
\noindent
{\bf Case 2b.2}\qua $D_1 \cup D_2$ is $K$--compressible. 

\medskip
Let $D$ be the $K$--compressing disk for $D_1 \cup D_2$. 
Without loss of generality, we may suppose that 
$D \cap D_1 \ne \emptyset$, 
$D \cap D_2 = \emptyset$. 
Let $D^*$ be a $K$--meridian disk of $X$ obtained from 
$D_1$ by $K$--compressing along $D$. 
By applying slight isotopy, we may suppose that $D^* \cap P = \emptyset$. 
Suppose that $D^* \subset B$. 
By pushing $A_1$ into $X$ along the parallelism through $A'$, 
we can $K$--isotope $P$ to $P'$ such that 
$P' \subset \text{Int} X$, and $P' \cap D^* = \emptyset$. 
Hence, by Proposition~\ref{Rubinstein-Scharlemann's lemma}, 
we see that $X \cup_Q Y$ is weakly $K$--reducible. 

\begin{figure}[ht!]
\begin{center}
\includegraphics[width=5cm, clip]{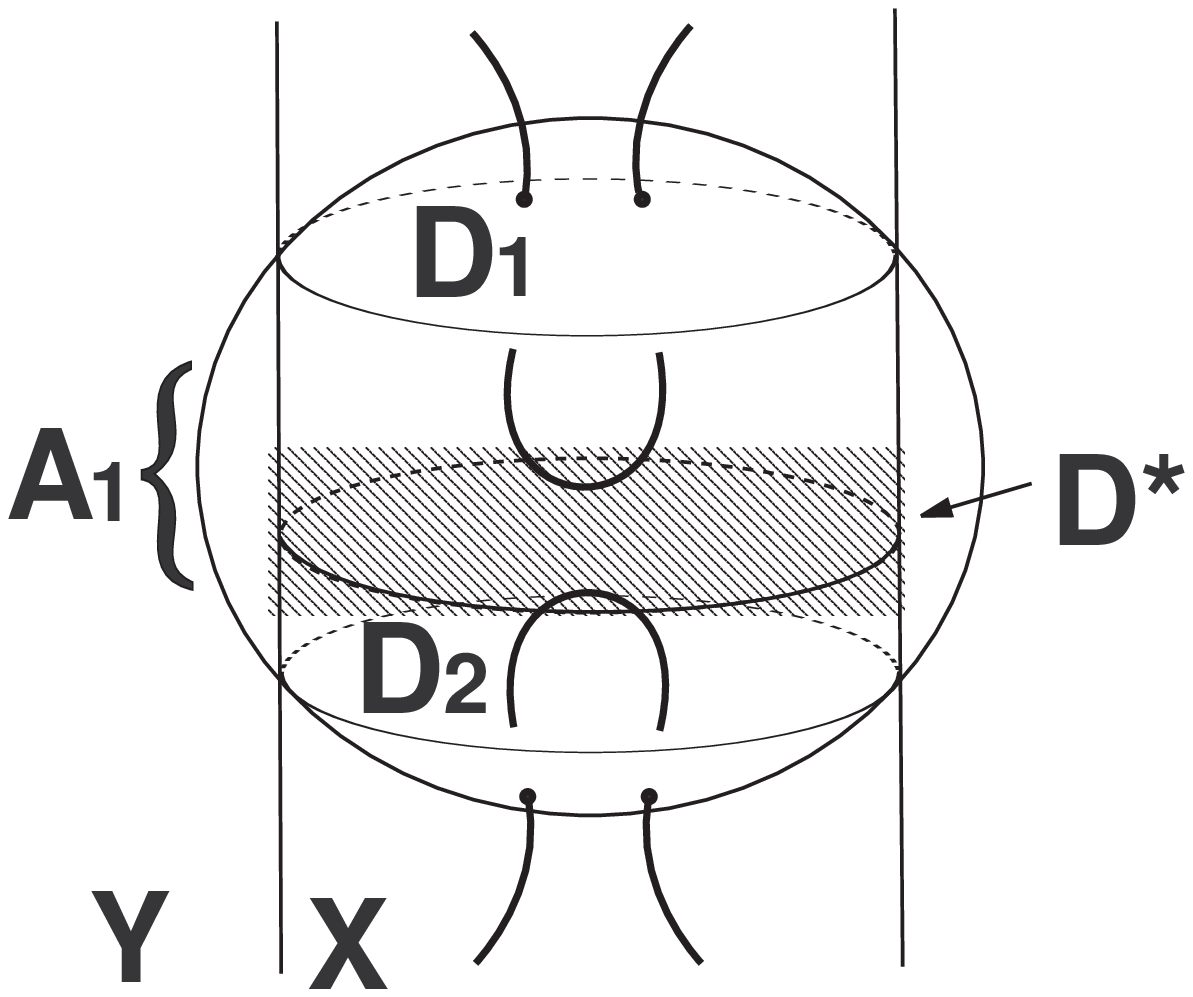}
\end{center}
\caption{}\label{6.1.4}
\end{figure}

Suppose that $D^* \subset A$ (Figure~\ref{6.1.4}). 
Let $D_B$ be a $K$--meridian disk in $B$. 
Since $K$ is not a trivial two component link, 
$\partial D$ and $\partial D_B$ are not isotopic in $P-K$. 
Hence $D_B \cap Q \ne \emptyset$. 
We suppose that $\vert D_B \cap Q \vert$ is minimal among all 
$K$--essential disks $D'$ in $B$ such that each component of 
$D' \cap D_1$ ($D' \cap D_2$, $D' \cap A_1$ respectively) 
is a $K$--essential arc in $D_1$ ($D_2$, $A_1$ respectively). 
Suppose that $D_B \cap Q$ contains a simple closed curve component. 
Let $D'$ be an innermost disk. 
By the minimality of $\vert D_B \cap Q \vert$, 
we see that $\partial D'$ is $K$--essential in $Q$. 
Note that $\partial D' \subset B$. 
If $D' \subset Y$, then the pair $D'$, $D^*$ gives a weak $K$--reducibility 
of $X \cup_Q Y$. 
If $D' \subset X$, then 
by pushing $A_1$ into $X$ along the parallelism through $A'$, 
we can $K$--isotope $P$ to $P'$ such that 
$P' \subset \text{Int} X$, and $P' \cap D' = \emptyset$. 
Hence, by Proposition~\ref{Rubinstein-Scharlemann's lemma}, 
we see that $X \cup_Q Y$ is weakly $K$--reducible. 
Suppose that each component of $D_B \cap Q$ is an arc. 
Let $\Delta (\subset D_B)$ be an outermost disk. 
If $\Delta \subset X$, 
then by using the argument as in the proof of Case~2b.1, 
we see that $X \cup_Q Y$ is weakly $K$--reducible. 
Suppose that $\Delta \subset Y$. 
Let $D^{**}$ be the disk obtained from $A_1$ by $K$--boundary compressing 
along $\Delta$. 

\medskip
\noindent
{\bf Claim}\qua $D^{**}$ is a $K$--meridian disk of $Y$. 

\begin{proof}
Suppose that $D^{**}$ is not a $K$--meridian disk of $Y$, 
ie, 
$D^{**}$ is $K$--parallel to a disk, say $D''$, 
in $\partial Y (= Q)$. 
Since $\Delta \subset B$, we see that $D'' \subset \text{cl} (Q-A')$. 
Note that $\text{cl} (Q-A')$ is recovered from $D''$ by adding a band 
along an arc intersecting $\Delta \cap Q$ in one point. 
This shows that $\text{cl} (Q-A')$ is an annulus not intersecting $K$. 
Hence $Q$ is a torus, and $X$ is a solid torus such that 
$Q \cap K = \emptyset$. 
However, since $D^*$ is a meridian disk of $X$, this implies that 
$X$ is $K$--reducible, 
contradicting Corollary~\ref{gamma-compression body is gamma-irreducible}. 
\end{proof}

By Claim, we see that, by applying a slight isotopy, 
we may suppose that $D^{**} \cap P = \emptyset$, 
and $D^{**} \subset B$. 
Hence the pair $D^*$, $D^{**}$ gives a weak $K$--reducibility 
of $X \cup_Q Y$. 

This completes the proof of 
Proposition~\ref{essential-essential intersection}
\end{proof}

\subsection*{6.B\qua Proof of 
Proposition~\ref{HS of 2-bridge knots are weakly reducible}}

Let $K$ be a non-trivial two bridge knot, and 
$C \cup_P V_2$ a genus~$g$ Heegaard splitting of $E(K)$ with $g \ge 3$. 
Note that $K$ satisfies the conditions of the assumption of 
Proposition~\ref{tubing}. 
Let $V_1$ be the handlebody in $S^3$ such that 
$\partial V_1 = P$, and 
$C \subset V_1$. 
Then $V_1 \cup V_2$ is a Heegaard splitting of $S^3$ which 
gives a genus~$g$, $0$--bridge position of $K$. 
By Propositions~\ref{HS of 2-bridge knots are weakly K-reducible} 
and \ref{tubing}, we have either one of the following. 

\medskip
\noindent
(1.1) 
There exists a weakly $K$--reducing pair of disks $D_1$, $D_2$ 
for $V_1 \cup V_2$ such that 
$D_1 \cap K = \emptyset$, and $D_2 \cap K = \emptyset$. 

\medskip
\noindent
(1.2) 
There exists a Heegaard splitting 
$V_{1,1} \cup_{P_1} V_{1,2}$ of $(S^3, K)$ which gives 
a genus $(g-1)$, $1$--bridge position of $K$ such that 
$V_1 \cup V_2$ is obtained from $V_{1,1} \cup_{P_1} V_{1,2}$ 
by a tubing. 

\medskip
If (1.1) holds,
then we immediately have the conclusion of 
Propositions~\ref{HS of 2-bridge knots are weakly reducible}. 
If (1.2) holds, then we further apply 
Propositions~\ref{HS of 2-bridge knots are weakly K-reducible} 
and \ref{tubing}, and we have either one of the following. 

\medskip
\noindent
(2.1) 
There exists a weakly $K$--reducing pair of disks $D_1$, $D_2$ for 
$V_{1,1} \cup_{P_1} V_{1,2}$ such that 
$D_1 \cap K = \emptyset$, and $D_2 \cap K = \emptyset$. 

\medskip
\noindent
(2.2) 
There exists a Heegaard splitting 
$V_{2,1} \cup_{P_2} V_{2,2}$ of $(S^3, K)$ which gives 
a genus $(g-2)$, $2$--bridge position of $K$ such that 
$V_{1,1} \cup_{P_1} V_{1,2}$ is obtained from 
$V_{2,1} \cup_{P_2} V_{2,2}$ by a tubing. 

\medskip
We claim that if (2.1) holds,
then we have the conclusion of 
Propositions~\ref{HS of 2-bridge knots are weakly reducible}. 
In fact, since $D_1 \cap K = \emptyset$, and $D_2 \cap K = \emptyset$, 
and tubing operations are performed in a small neighborhood 
of $K$, the pair $D_1$, $D_2$ survives in $V_1 \cup V_2$ 
to give a weak reducibility. 
If (2.2) holds, then we further apply 
Propositions~\ref{HS of 2-bridge knots are weakly K-reducible} 
and \ref{tubing}, and we have either one of the following. 

\medskip
\noindent
(3.1) 
There exists a weakly $K$--reducing pair of disks $D_1$, $D_2$ for 
$V_{2,1} \cup_{P_2} V_{2,2}$ such that 
$D_1 \cap K = \emptyset$, and $D_2 \cap K = \emptyset$. 

\medskip
\noindent
(3.2) 
There exists a Heegaard splitting 
$V_{3,1} \cup_{P_3} V_{3,2}$ of $(S^3, K)$ which gives 
a genus $(g-3)$, $3$--bridge position of $K$ such that 
$V_{2,1} \cup_{P_2} V_{2,2}$ is obtained from 
$V_{3,1} \cup_{P_3} V_{3,2}$ by a tubing. 

\medskip
Then we apply the same argument as above, and so on. 
Then either we have the conclusion of 
Propositions~\ref{HS of 2-bridge knots are weakly reducible}, 
or the procedures are repeated $(g-1)$ times to give the following. 

\medskip
\noindent
($g$.1) 
There exists a weakly $K$--reducing pair of disks $D_1$, $D_2$ for 
$V_{g-1,1} \cup_{P_{g-1}} V_{g-1,2}$ such that 
$D_1 \cap K = \emptyset$, and $D_2 \cap K = \emptyset$. 

\medskip
\noindent
($g$.2) 
There exists a Heegaard splitting 
$V_{g,1} \cup_{P_g} V_{g,2}$ of $(S^3, K)$ which gives 
a genus $0$, $g$--bridge position of $K$ such that 
$V_{g-1,1} \cup_{P_{g-1}} V_{g-1,2}$ is obtained from 
$V_{g,1} \cup_{P_g} V_{g,2}$ by a tubing. 

\medskip
If ($g$.1) holds, then by using the arguments as above, 
we see that we have the conclusion of 
Propositions~\ref{HS of 2-bridge knots are weakly reducible}. 
Suppose that ($g$.2) holds. 
Then we see that 
there exists a weakly reducing pair of disks $D_1$, $D_2$ for 
$V_{g,1} \cup_{P_g} V_{g,2}$ such that 
$D_1 \cap K = \emptyset$, and $D_2 \cap K = \emptyset$
(see Remark~\ref{remark of tubing}), 
and this together with the arguments as for the case ($g$.1), 
we see that we have the conclusion of 
Propositions~\ref{HS of 2-bridge knots are weakly reducible}. 

This completes the proof of 
Propositions~\ref{HS of 2-bridge knots are weakly reducible}. 

\section{Proof of Theorem~\ref{main theorem}}

Let $K$ be a knot in a closed 3--manifold $M$. 

\begin{definition}
A {\it tunnel} for $K$ is an embedded arc $\sigma$ 
in $S^3$ such that $\sigma \cap K = \partial \sigma$.  
We say that a tunnel $\sigma$ for $K$ is {\it unknotting} 
if $S^3-$Int~$N(K \cup \sigma ,S^3)$ 
is a genus two handlebody.  
\end{definition} 

For a tunnel $\sigma$ for $K$, let $\hat{\sigma} = \sigma \cap E(K)$. 
Then $\hat{\sigma}$ is an arc properly embedded in $E(K)$, 
and we may regard that $N(K \cup \sigma)$ is obtained from 
$N(K)$ by attaching $N(\hat{\sigma}, E(K))$, 
where $N(\hat{\sigma}, E(K)) \cap N(K)$ consists of two disks, 
ie, 
$N(\hat{\sigma}, E(K))$ is a 1--handle attached to $N(K)$. 

\begin{definition}
Let $\sigma_1$, $\sigma_2$ be tunnels for $K$. 
We say that $\sigma_1$ is {\it isotopic} to $\sigma_2$ 
if there is an ambient isotopy $h_t$ $(0 \le t \le 1)$ of 
$E(K)$ such that $h_0 = \text{id}_{E(K)}$, and 
$h_1 ( \hat{\sigma}_1 ) = \hat{\sigma}_2$. 
\end{definition}

\begin{remark}\label{unknotting tunnel}
Let $\sigma$ be an unknotting tunnel for $K$, 
and let $V = N(K \cup \sigma , M)$, and $W = \text{cl}(M-V)$. 
Note that $V \cup W$ is a Heegaard splitting of $(M,K)$, 
which gives a genus two, $0$--bridge position of $K$. 
Let $\sigma_1$, $\sigma_2$ be unknotting tunnels for $K$, 
and $V_1 \cup_{P_1} W_1$, $V_2 \cup_{P_2} W_2$ Heegaard splittings 
obtained from $\sigma_1$, $\sigma_2$ respectively as above. 
Then it is known that 
$\sigma_1$ is isotopic to $\sigma_2$ 
if and only if $P_1$ is $K$--isotopic to $P_2$. 
\end{remark}

\begin{figure}[ht!]
\begin{center}
\includegraphics[width=5cm, clip]{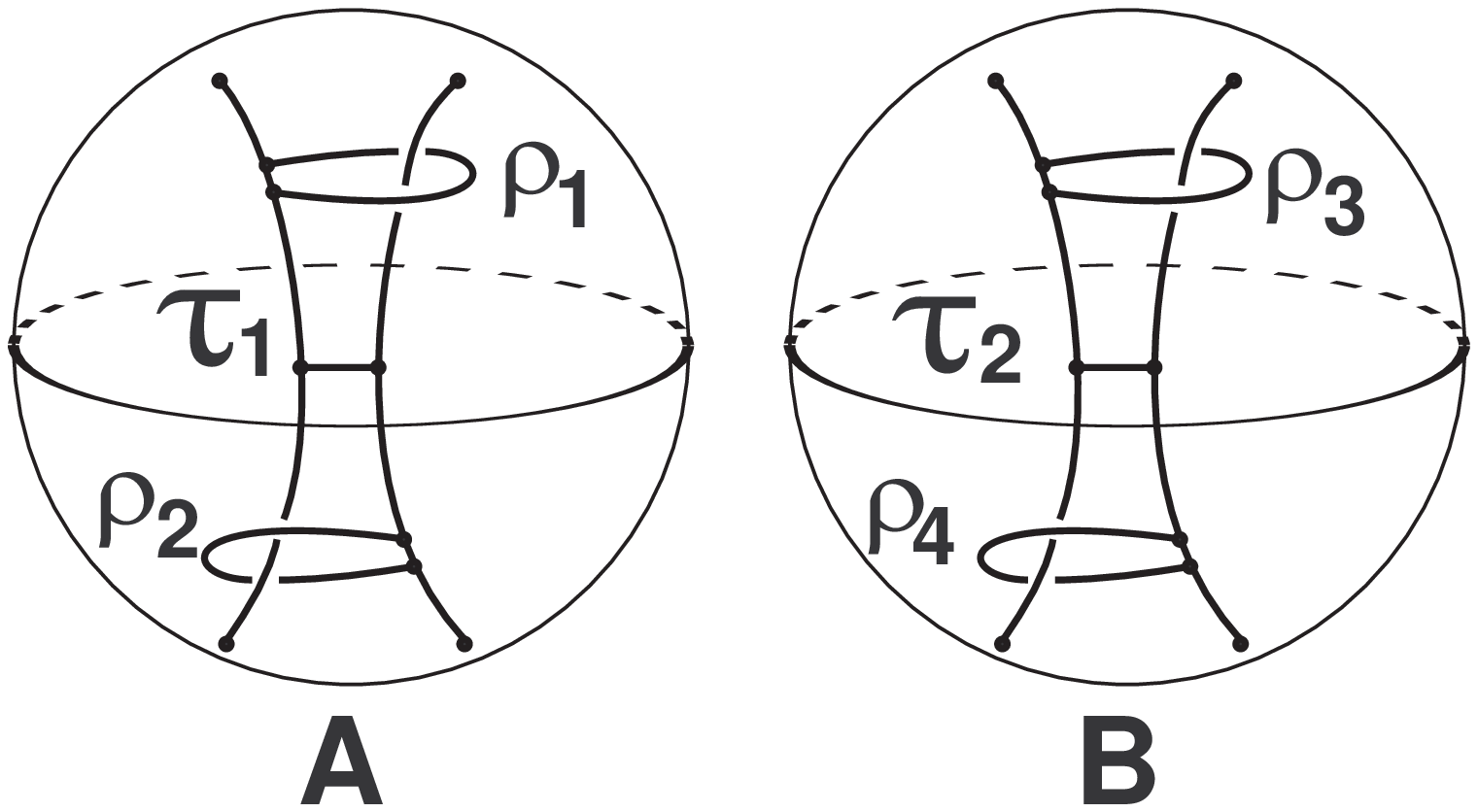}
\end{center}\vspace{5pt}
\caption{}\label{7.1}
\end{figure}

Now, in the rest of this paper, 
let $K$ be a non-trivial 2--bridge knot, and $A \cup_P B$ a 
genus 0 Heegaard splitting of $S^3$, which gives a two bridge 
position of $K$ (Figure~\ref{7.1}). 

\subsection*{7.A\qua Genus two Heegaard splittings of $E(K)$} 

Here we show the next lemma on unknotting tunnels of $K$, 
which is used in the proof of Theorem~\ref{main theorem}. 

\begin{lemma}\label{pair of meridian disks for unknotting tunnels}
Let $\sigma$ be an unknotting tunnel for $K$, and 
$V \cup W$ a Heegaard splitting obtained from $\sigma$ as in 
Remark~\ref{unknotting tunnel}. 
Then there exist meridian disks $D_1$, $D_2$ of $V$, $W$ respectively 
such that $D_1$ intersects $K$ transversely in one point, 
$D_1 \cap N(\hat{\sigma}, E(K)) = \emptyset$, 
and 
$\partial D_1$ intersects $\partial D_2$ transversely in one point. 
\end{lemma}

\begin{figure}[ht!]
\begin{center}
\includegraphics[width=5cm, clip]{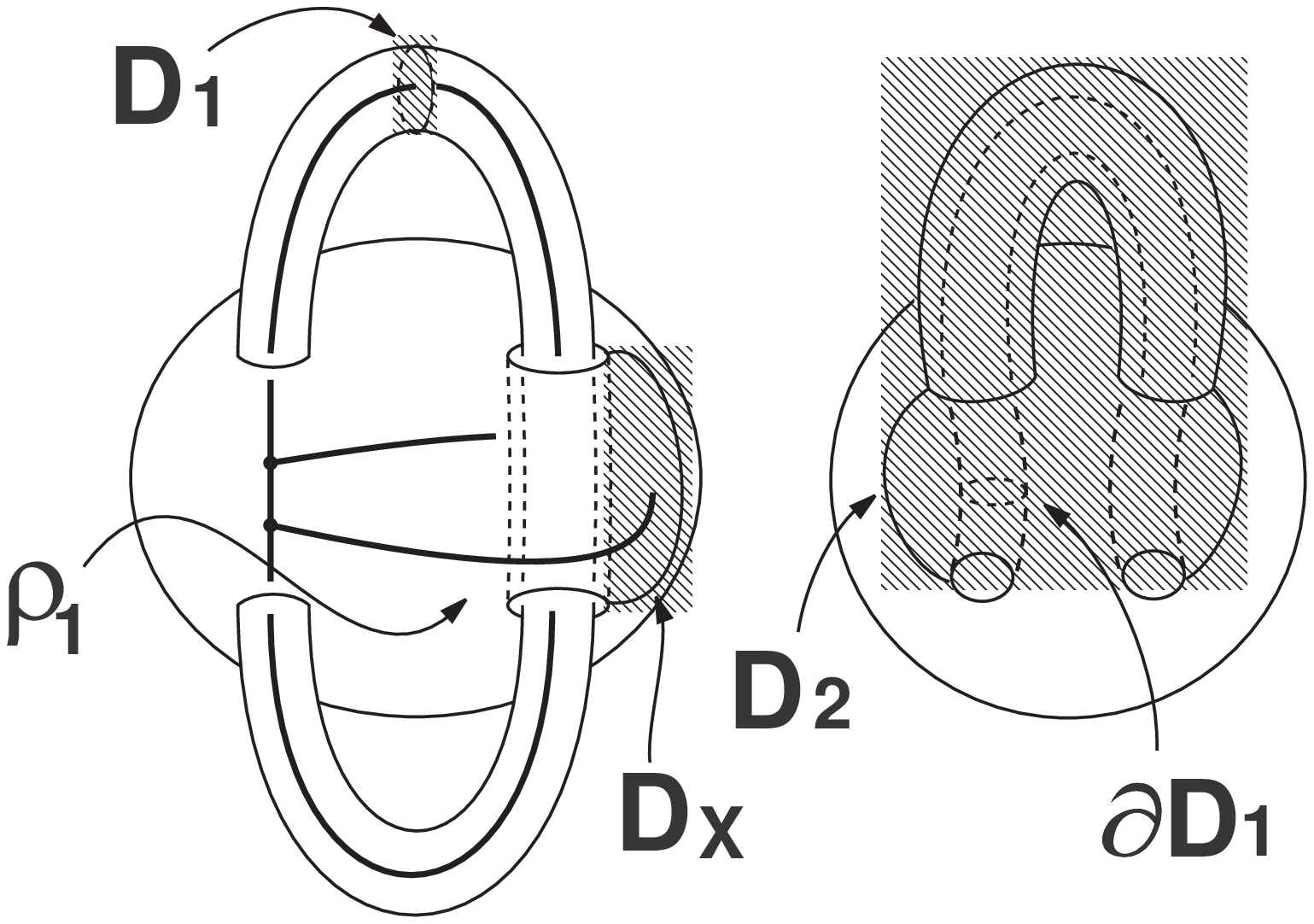}
\end{center}
\caption{}\label{7.1.2}
\end{figure}

\begin{proof}
We note that $\sigma$ is isotopic to either one of the six unknotting 
tunnels $\tau_1$, $\tau_2$, $\rho_1$, $\rho_2$, $\rho_3$, or $\rho_4$ 
in Figure~\ref{7.1} (see \cite{GST} or \cite{K}). 
Suppose that $\sigma$ is isotopic to $\tau_i$, 
$i=1$ or $2$, say $1$. 
Then we may regard that $V = A \cup N(K \cap B, B)$ 
(Figure~\ref{7.1.1}). 
Here $N(\hat{\sigma}, E(K)) = N(D_A, A)$, 
where $D_A$ is a disk properly embedded in $A$, 
such that $D_A$ separates the components of $K \cap A$, 
and $N(D_A, A) \cap N(K \cap B, B) = \emptyset$ 
(hence, $D_A$ is properly embedded in $V$). 
Then we can take a pair $D_1$, $D_2$ satisfying the conclusion of 
Lemma~\ref{pair of meridian disks for unknotting tunnels} 
as in Figure~\ref{7.1.1}.

\begin{figure}[ht!]
\begin{center}
\includegraphics[width=5cm, clip]{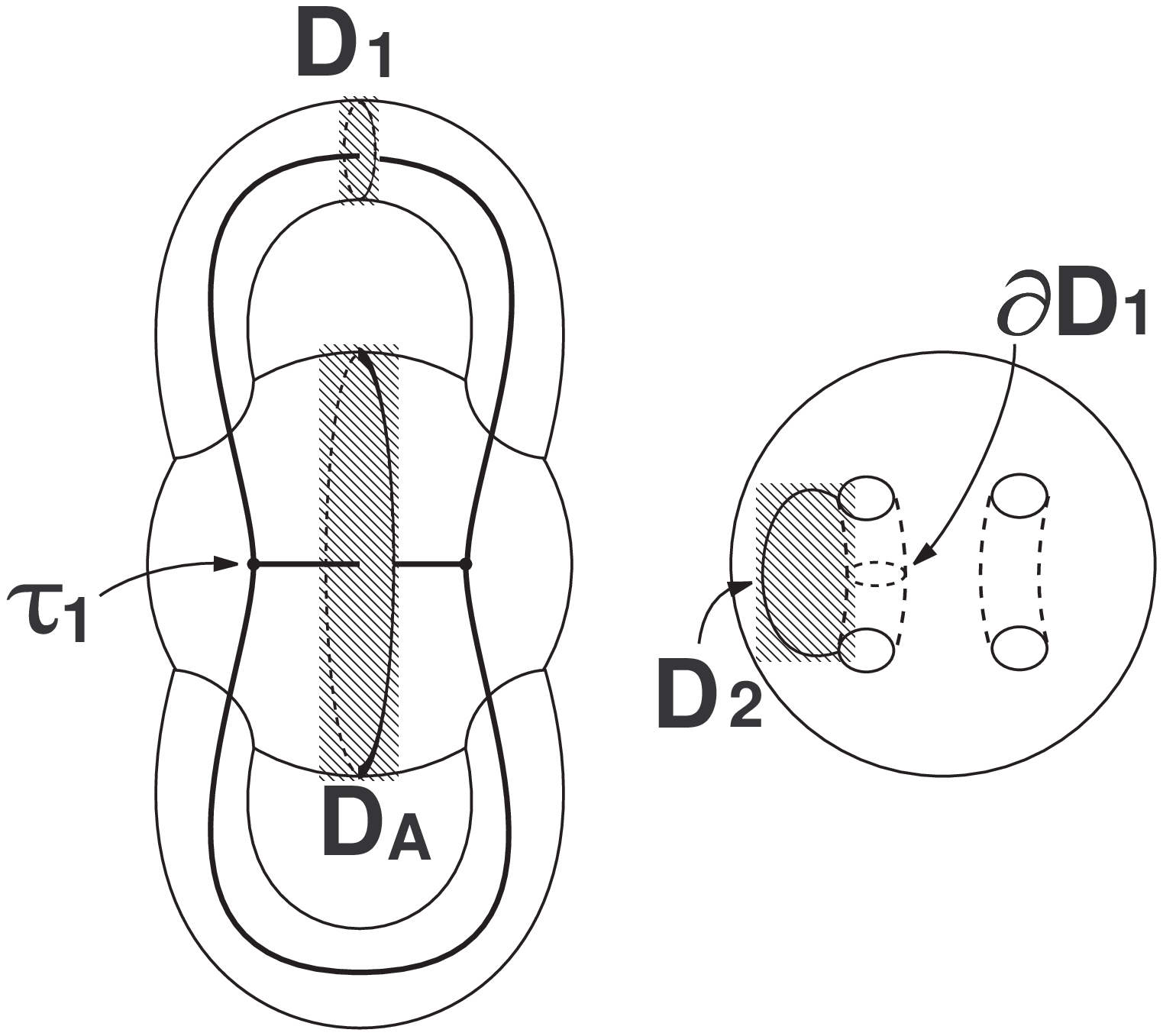}
\end{center}\vspace{-10pt}
\caption{}\label{7.1.1}
\end{figure}

Suppose that $\sigma$ is isotopic to $\rho_i$, 
$i=1$, $2$, $3$, or $4$, say $1$.  
Then we may regard that $V$ is obtained from the Heegaard 
splitting $A \cup_P B$ of $(S^3, K)$ as follows. 

Let $a$ be the component of $K \cap A$, which is disjoint from $\sigma$, 
and 
$V' = \text{cl}(A - N(a,A))$, 
$W' = B \cup N(a,A)$. 
Let $a' = a \cup (K \cap B)$. 
Note that $a'$ is an arc properly embedded in $W'$. 
Then $V = V' \cup N(a', W')$. 
See Figure~\ref{7.1.2}. 
That is, $V \cup W$ is obtained from $A \cup_P B$ 
by successively tubing along $a$, and $a'$. 
We can take a pair $D_1$, $D_2$ satisfying the conclusion of 
Lemma~\ref{pair of meridian disks for unknotting tunnels}, 
as in Figure~\ref{7.1.2}.\end{proof}

\subsection*{7.B\qua Irreducible Heegaard splittings of 
$\text{(torus)} \times [0,1]$} 

In \cite{BO}, M Boileau, and J-P Otal gave a classification 
of Heegaard splittings of $\text{(torus)} \times [0,1]$, 
and M.Scharlemann, and A.Thompson \cite{ST} proved that the same kind 
of results hold for $F \times [0,1]$, where $F$ is any closed orientable 
surface. 
The result of Boileau--Otal will be used for 
the proof of Theorem~\ref{main theorem}, 
and in this section we quickly state it. 

Let $T$ be a torus. 
Let $Q_1$ be the surface $T \times \{ 1/2 \}$ in $T \times [0,1]$. 
It is clear that $Q_1$ separates $T \times [0,1]$ into two trivial 
compression bodies. 
Hence $Q_1$ is a Heegaard surface of $T \times [0,1]$. 
We call this Heegaard splitting {\it type~1}. 
Let $a$ be a vertical arc in $T \times [0,1]$. 
Let $V_1 = N( (T \times \{0, 1 \}) \cup a, T \times [0,1])$, 
and 
$V_2 = \text{cl}(T \times [0,1] - V_1)$. 
It is easy to see that $V_1$ is a compression body, 
$V_2$ is a genus two handlebody, and 
$V_1 \cap V_2 = \partial_+ V_1 = \partial_+ V_2 (= \partial V_2)$. 
Hence $V_1 \cup V_2$ is a Heegaard splitting of $T \times I$. 
We call this Heegaard splitting {\it type~2}. 
Then in \cite[Th\'eor\`eme 1.5]{BO}, or \cite[Main theorem 2.11]{ST}, 
the following is shown. 

\begin{theorem}\label{HS of torus times [0,1]}
Every irreducible Heegaard splitting of $T \times [0,1]$ is isotopic
to either a Heegaard splitting of type~1 or type~2. 
\end{theorem}

\subsection*{7.C\qua Proof of Theorem~\ref{main theorem}} 

Let $C_1 \cup_P C_2$ be a genus $g$ Heegaard splitting of 
the exterior of $K$, $E(K)=\text{cl}(S^3-N(K))$, 
with $g \ge 3$ and $\partial_- C_1 = \partial E(K)$. 
Then, by Proposition~\ref{HS of 2-bridge knots are weakly reducible}, 
we see that $C_1 \cup_P C_2$ is weakly reducible. 
By Proposition~\ref{CG Theorem 3.1}, 
either $C_1 \cup_P C_2$ is reducible, or 
there is a weakly reducing collection of disks 
$\Delta$ for $P$ such that each component of $\hat{P}( \Delta )$ 
is an incompressible surface in $E(K)$, 
which is not a 2--sphere. 
Suppose that the second conclusion holds and let 
$M_j$ $(j=1, \dots , n)$, $M_{j,i}$ $(i=1,2)$, and 
$C_{1,1} \cup_{P_1} C_{1,2}, \cdots , C_{n,1} \cup_{P_n} C_{n,2}$ 
be as in Section~4. 
Note that each component of $\partial_- C_{i,j}$ is either 
$\partial E(K)$ or a closed incompressible surface in
$\text{Int} E(K)$. 
Since every closed incompressible surface in $\text{Int} E(K)$ 
is a $\partial$--parallel torus, 
we see that the submanifolds $M_1, \dots , M_n$ 
lie in $E(K)$ in a linear configuration, 
ie, 
by exchanging the subscripts if necessary, we may suppose that 

\begin{enumerate}
\item 
$\partial_- C_{1,1} = \partial E(K)$, 

\item 
For each $i$ $(1 \le i \le n-1)$, 
$M_i$ is homeomorphic to (torus)$\times [0,1]$, and 
$M_i \cap M_{i+1} = F_i$: a $\partial$--parallel torus in $E(K)$. 
\end{enumerate}

\medskip
\noindent
{\bf Claim 1}\qua If $n>2$, then $C_1 \cup_P C_2$ is reducible. 

\begin{proof}
Let 
$M_1' = \text{cl}(C_1 - M_{n,1})$, and 
$M_2' = \text{cl}(C_2 - M_{n,2})$. 
Then from the pair $M_1'$, $M_2'$ we can obtain, as in Section~4, 
a Heegaard splitting, 
say $C_1' \cup_{P'} C_2'$, of the product region between 
$F_{n-1}$ and $\partial E(K)$. 
Since $n>2$, we see, by \cite[Remark 2.7]{Schu'}, that 
$\text{genus}(P') > 2$. 
Hence by Theorem~\ref{HS of torus times [0,1]}, 
$C_1' \cup_{P'} C_2'$ is reducible. 
Hence, by Lemma~\ref{reducible implies reducible}, 
$C_1 \cup_P C_2$ is reducible. 
\end{proof}

By Claim~1, we may suppose, in the rest of the proof, that $n=2$. 
Now we prove Theorem~\ref{main theorem} by the induction on $g$. 

Suppose that $g=3$. 
By Lemma~\ref{reducible implies reducible}, 
we may suppose that both $C_{1,1} \cup_{P_1} C_{1,2}$, and 
$C_{2,1} \cup_{P_1} C_{2,2}$ are irreducible. 
By Lemma~\ref{non trivial compression body} and 
Theorem~\ref{HS of torus times [0,1]}, 
we see that $C_{1,1}$ is a genus 2 compression body 
with $\partial_-C_{1,1} = \partial E(K) \cup F_1$, 
and $C_{1,2}$ is a genus 2 handlebody. 

\medskip
\noindent
{\bf Claim 2}\qua $(M_{1,1} \cap P) \subset (M_{1,2} \cap P)$. 

\begin{proof}
Suppose not. 
Then, by Lemma~\ref{HS from weakly reducing collection of disks}, 
we see that 
$(M_{1,1} \cap P) \supset (M_{1,2} \cap P)$. 
Recall that 
$C_{1,1} = \text{cl}(M_{1,1}-N(\partial_+ M_{1,1}, M_{1,1}))$. 
This implies that 
$\partial_- M_{1,1} = \partial_- C_{1,1}$. 
Note that $C_{1,1} \cup_{P_1} C_{1,2}$ is a Heegaard splitting of 
type~2 in Section~7.B. 
These show that 
$\partial_- M_{1,1} =\partial E(K) \cup F_1$. 
However, this is impossible since 
$\partial_- M_{1,1} \subset \partial E(K)$. 
\end{proof}

\begin{figure}[ht!]
\begin{center}
\includegraphics[width=6cm, clip]{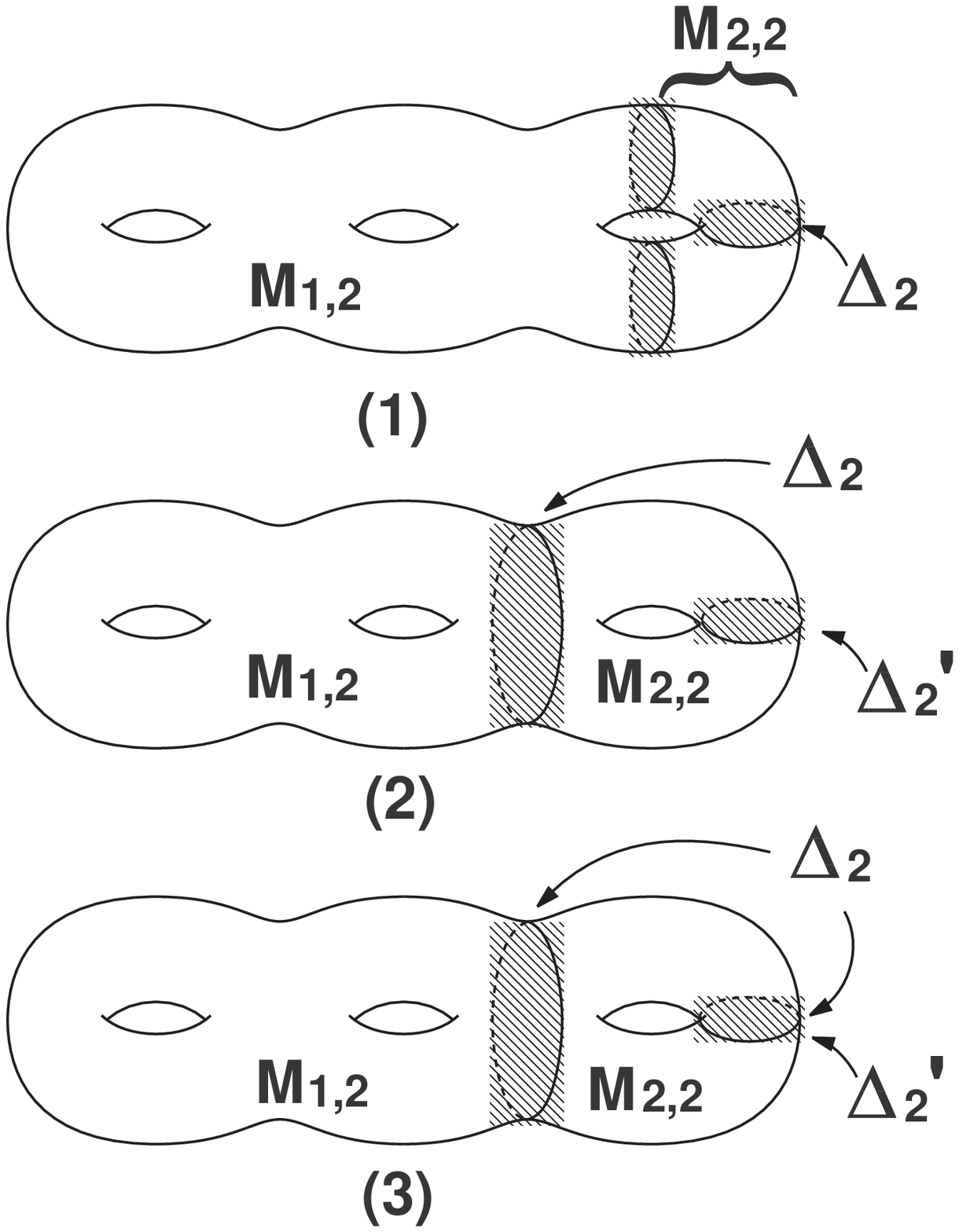}
\end{center}\vspace{10pt}
\caption{}\label{7.3.1}
\end{figure}

By Claim~2, we see that $M_{1,2}$ is a genus two handlebody. 
Hence $\Delta_2$ is either one of Figure~\ref{7.3.1}, 
ie, 
either (1) $\Delta_2$ consists of a non-separating disk in $C_2$, 
(2) $\Delta_2$ consists of a separating disk in $C_2$, or 
(3) $\Delta_2$ consists of two disks, 
one of which is a separating disk, and the other is 
a non-separating disk in $C_2$. 

\begin{figure}[ht!]
\begin{center}
\includegraphics[width=6cm, clip]{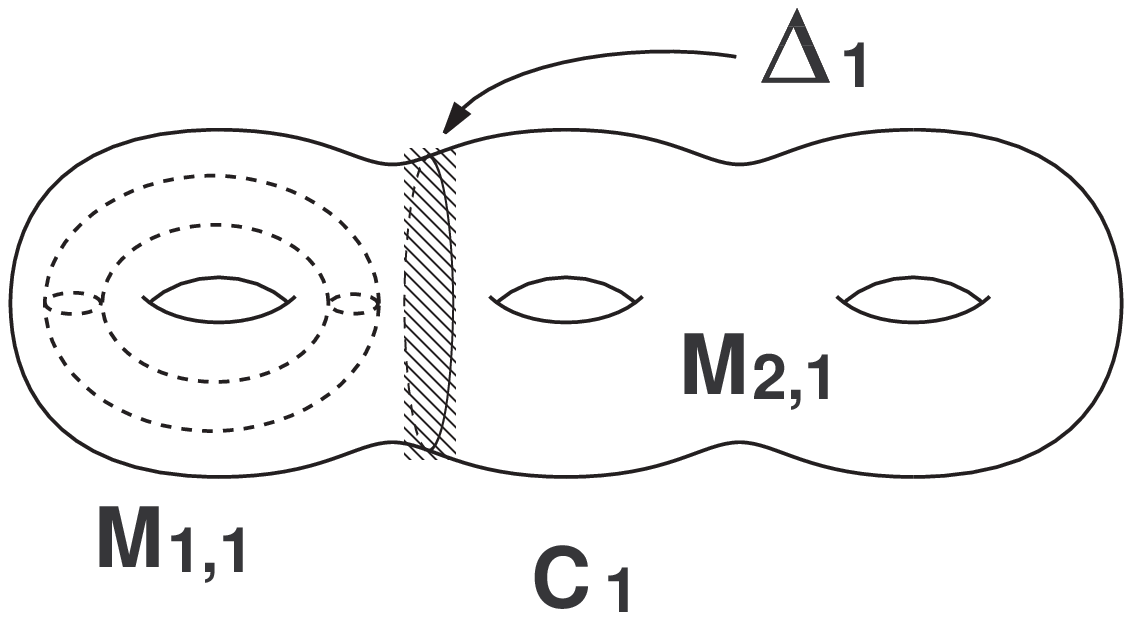}
\end{center}\vspace{10pt}
\caption{}\label{7.3.2}
\end{figure}

Suppose that $\Delta_2$ is of type~(1) in Figure~\ref{7.3.1}. 
Since no component of $\hat{P} ( \Delta )$ is a 2--sphere, 
we see that $\partial \Delta_1 \subset M_{1,2}$. 
By Claim~2, we see that 
$(M_{2,1} \cap P) \supset (M_{2,2} \cap P)$. 
Since 
$\partial (M_{2,1} \cup M_{2,2}) = \partial M_2 = F_1$: a torus, 
we see that $M_{2,1}$ is a genus two handlebody, and $\Delta_1$ 
consists of a separating disk in $C_1$ (Figure~\ref{7.3.2}). 

Let $N_K = \text{cl} (S^3-M_2)$. 
Since $F_1$ is a $\partial$--parallel torus in $E(K)$, 
we see that $N_K$ is a regular neighborhood of $K$, 
hence $M_2$ is an exterior of $K$. 
Note that $M_{2,2}$ is a 1--handle attached to $N_K$ such that 
$\text{cl}(S^3-(N_K \cup M_{2,2})) = M_{2,1}$, a genus 
two handlebody. 
This shows that $M_{2,2}$ is a regular neighborhood of an arc 
properly embedded in $M_2$, which comes from an unknotting tunnel of $K$. 
Hence, by Lemma~\ref{pair of meridian disks for unknotting tunnels}, 
we see that there is a pair of disks $D_1$, $D_2$ in 
$N_K \cup M_{2,2}$, $M_{2,1}$ respectively such that 
$D_1$ intersects $K$ transversely in one point, 
$D_1 \cap M_{2,2} = \emptyset$, and 
$\partial D_1$ intersects $\partial D_2$ transversely 
in one point. 
Here, by deforming $D_2$ by an ambient isotopy of $M_{2,1}$ 
if necessary, we may suppose that $D_2 \cap \Delta_1 = \emptyset$ 
(hence, $D_2$ is a meridian disk of $C_1$). 
Since $D_1$ and $K$ intersect transversely in one point, 
we may suppose that $D_1 \cap E(K)$ $(= D_1 \cap M_1)$ is 
a vertical annulus, say $A_1$, properly embedded in 
$M_1$ $( \cong T^2 \times [0,1])$. 
Recall that $C_{1,1} \cup_{P_1} C_{1,2}$ is a type~2 
Heegaard splitting of $M_1$. 
This implies that there exists a vertical arc $a$ in $M_1$ 
such that 
$M_{1,1} = N(\partial E(K) \cup a, M_1)$. 
Since $a$ is vertical, we may suppose, by isotopy, that $a \subset A_1$, 
ie, 
$a$ is an essential arc properly embedded in $A_1$. 
Let $\ell$ be the component of $\partial A_1$ contained in $\partial E(K)$. 
Hence $A_1 \cap C_2$ $= A_1 \cap M_{1,2}$ 
$= \text{cl}(A_1 - N(\ell \cup a,M_1))$, and 
this is a disk, say $D_1'$, properly embedded in $C_2$. 
Obviously $\partial D_1'$ and $\partial D_2$ intersect 
transversely in one point. 
Recall that $D_2$ ($D_1'$ respectively) is a disk properly embedded 
in $C_1$ ($C_2$ respectively). 
Hence $C_1 \cup_P C_2$ is stabilized and this shows that 
$C_1 \cup_P C_2$ is reducible if $g=3$
(see 2 of Remark~\ref{reducible Heegaard splitting}). 

Suppose that $\Delta_2$ is of type (2) or (3) in Figure~\ref{7.3.1}.
Then we take $\Delta_2'$ as in Figure~\ref{7.3.1}, and  let 
$\Delta' = \Delta_1 \cup \Delta_2'$. 
We note that $\Delta'$ is a weakly reducing collection of disks for $P$, 
where $\Delta'$ is of type (1) in  Figure~\ref{7.3.1}. 
Let $F_1'$ be the torus obtained from $\Delta'$, 
which is corresponding to $F_1$. 
It is directly observed from Figure~\ref{7.3.1} that $F_1'$ is 
isotopic to $F_1$. 
Hence we can apply the argument for type~1 weakly reducing 
collection of disks to $\Delta'$, and we can show that 
$C_1 \cup_P C_2$ is reducible. 

Suppose that $g \ge 4$. 
If $\text{genus}(P_1) > 2$, then 
by Theorem~\ref{HS of torus times [0,1]} 
and Lemma~\ref{reducible implies reducible}, 
we see that $C_1 \cup_P C_2$ is reducible. 
Suppose that $\text{genus}(P_1) = 2$. 
Then, by \cite[Remark 2.7]{Schu'}, 
we see that $\text{genus}(P_2) = g-1$. 
Hence, by the assumption of the induction, 
we see that $C_{2,1} \cup_{P_2} C_{2,2}$ is reducible. 
Hence, by Lemma~\ref{reducible implies reducible}, 
$C_1 \cup_P C_2$ is reducible. 

This completes the proof of Theorem~\ref{main theorem}.

\end{document}